\documentclass[10pt]{article}
\usepackage{graphicx} 
\usepackage{authblk}

\usepackage{bbm}
\usepackage{amsmath}
\usepackage{amssymb}
\usepackage{graphicx}

\usepackage{import}
\usepackage{stackrel}
\usepackage{amsthm}
\usepackage{caption}
\usepackage{authblk}
\usepackage{placeins}
\usepackage{adjustbox}
\usepackage[nottoc]{tocbibind}

\usepackage{ragged2e}
\usepackage{etoolbox}
\usepackage{subcaption}
\usepackage{float}
\usepackage[english]{babel}
\usepackage[a4paper, total={6in, 10in}]{geometry}
\usepackage{booktabs}
\usepackage{colortbl}
\usepackage{xcolor}
\usepackage{xfrac}
\usepackage{constants}
\usepackage{mathrsfs}
\usepackage{bbm}
\usepackage{amssymb}
\usepackage{hyperref}

\newtheorem{maintheorem}{Theorem}

\newtheorem{theorem}{Theorem}[section]
\newtheorem{proposition}[theorem]{Proposition}
\newtheorem{lemma}[theorem]{Lemma}
\newtheorem{sublemma}[theorem]{Sublemma}
\newtheorem{corollary}[theorem]{Corollary}

\theoremstyle{definition}
\newtheorem{definition}[theorem]{Definition}
\newtheorem{question}{Question}
\newtheorem*{conjecture}{Conjecture}
\newtheorem{remark}[theorem]{Remark}
\newtheorem{note}[theorem]{\textbf{{Note}}}

\newcommand{\todelete}[1]{{\color{red}#1}} 
\newcommand{\toadd}[1]{{\color{blue}#1}}

\title{Doubly Intermittent Singular Skew Products}
\author[1]{Stefano Luzzatto}
\author[2]{Rizwan Rizwanullah}
\affil[1]{Abdus Salam International Centre for Theoretical Physics (ICTP), Trieste, Italy.  \emph{luzzatto@ictp.it}}
\affil[2]{Kungliga Tekniska H\"ogskolan (KTH), Stockholm, Sweden. \emph{rizwanullahssgm@gmail.com}}

\date{\today}

\begin{document}
\maketitle


\begin{abstract}
We study a large class of skew product systems, which can also be thought of as random dynamical systems, coupling a uniformly expanding circle map with a full branch interval map with two indifferent fixed points and unbounded derivatives. This significantly generalizes the setting of several skew products considered in the literature, in particular the Gou\"ezel skew products  \cite{Gou07}, in which the fibre maps are intermittency maps with a single indifferent fixed point and bounded derivatives.
Through the construction of an induced full branch map, we  prove various statistical properties of the skew products, including existence of an invariant probability measure equivalent to Lebesgue and decay of correlations. 
\end{abstract}

\section{Introduction and Statement of Results}

Over the last several decades, there has been a substantial amount of research on skew-product systems, i.e. maps of the form \( \widehat F(\omega, x) = (F(\omega), T_{\omega}(x)) \) in which a \emph{base} map \( F \) drives a family of \emph{fibre} maps \( T_{\omega} \). Even in the two-dimensional case, these can exhibit a level of complexity and technical challenges that go well beyond the study of the one-dimensional systems that form them while having additional structure which makes them often easier to study than general two-dimensional maps. 

Skew products also admit a natural interpretation as \emph{random} dynamical systems. Indeed, for a fixed \( \omega \), the second coordinate of the orbit of \( (\omega, x) \) under \( \widehat F \) is given by the composition \( T_{F^{n-1}(\omega)} \circ \cdots \circ T_{F(\omega)} \circ T_{\omega}(x) \), which we will denote by \( T^{n}_{\omega}(x) \), i.e. by iterating a sequence of maps which is not fixed in advance but is chosen at each step according to the orbit of \( \omega \) under \( F \). When the base map is chaotic, as in our setting where \( F \) is a uniformly expanding circle map and is therefore isomorphic to a Bernoulli shift with respect to its unique absolutely continuous invariant probability measure, we can think of the sequence \( \omega, F(\omega), F^{2}(\omega), \dots \) as a (correlated) sequence of random variables, and thus of the fibre dynamics as the iteration of a random sequence of maps, where the ``noise'' driving the choice of map at each step is generated by the base. Conversely, the classical setting of random maps, in which the maps are chosen independently at each step from some finite family according to a fixed probability distribution, is exactly a skew product over a full shift with a Bernoulli measure, see e.g. \cite{Kif86, Arn98, KL06} for the general theory of random dynamical systems and Section \ref{sec:context} for further references. From this point of view the skew products we consider here correspond to a situation in which the maps are chosen from a continuous family, depend smoothly on the noise, and the noise itself is generated deterministically by a chaotic map rather than by independent trials. The two points of view are complementary. The random maps formulation suggests the natural questions, namely the existence of a physical measure and its statistical properties, which in the language of random dynamical systems correspond to the \emph{annealed} statistics of the fibre dynamics (see Section \ref{sec:main_results} for a more precise discussion), whereas the skew product formulation makes explicit the geometry of the problem, and it is precisely this geometry which we exploit in our arguments.

The family of skew products we study in this paper have uniformly expanding maps on the base and fibre maps which are a generalization of the \emph{doubly intermittent} maps introduced and studied in  \cite{CLM23, CL24}, which  themselves significantly generalize and extend classical intermittency maps by allowing \emph{both} fixed points to be neutral and the derivative at the discontinuity to be unbounded. Generalizations of doubly intermittent maps were subsequently studied in \cite{MS25,CMT26}.

While our general strategy, involving the construction of an induced map, is similar to that used for the skew-product maps with classical intermittency maps on the fibres, studied in  \cite{Gou07}, the characteristics of our fibre maps introduce some significant technical challenges and we introduce some interesting technical innovations which might be relevant in other more general settings. 

In Section \ref{sec:def} we give the precise definition of the class of skew product maps we consider. In Section~\ref{sec:main_results} we state our main results. In Section \ref{sec:context} we give some context and review of existing related results. Finally, in Section \ref{sec:overview} we give an overview of our proofs, and indeed prove our main results modulo two technical propositions which will then be proved in the rest of the paper.

\paragraph{Acknowledgments.} This project was started when R.R. was at the Scuola Internazionale Superiore di Studi Avanzati (SISSA) in Trieste, Italy, and completed with the generous support and hospitality of  Kungliga Tekniska H\"ogskolan (KTH) in Stockholm, Sweden, the Abdus Salam International Centre for Theoretical Physics (ICTP) in Trieste, Italy, and funding from the Swedish Research Council VR grant SFB 67578 \emph{Renormalization in Dynamics} of Liviana Palmisano.  We thank Asad Ullah for introducing us to \cite{Gou07} and for many initial discussions, and Sebastien Gou\"ezel for pointing out a minor but important  gap in our arguments in a preliminary version.  AI was used at the final stages of writing to proofread the text and draw the figures. It helped to improve the presentation  and suggested  ways to optimize some of the arguments. All resulting suggestions were  checked and edited by the authors, who take full responsibility for the contents. 

\subsection{Doubly Intermittent Singular Skew Products}\label{sec:def}
 Let $I= [-1,1]$ and \( \mathbb S^1\) be the unit circle. We  define a class of maps 
\(
\widehat{F}: \mathbb{S}^{1}\times I \rightarrow \mathbb{S}^{1} \times I
\)
 of the form 
\begin{equation}\label{eq:definnition_of_That}
   \widehat{F}(\omega,x) = (F(\omega), T(\omega, x)) = (F(\omega), T_{\omega}(x)).
\end{equation}
Notice that to simplify the notation, as in \eqref{eq:definnition_of_That}, we will often switch between writing  \(T(\omega, x)\), to emphasize that \( T \) is a function of two variables, and writing \( T_{\omega}(x)\), to emphasize that for many estimates we can consider \( \omega\) fixed and see \( T_{\omega}(x) \) as a function of \( x \).  We will apply the same principle also to other functions to be introduced below. 
We assume that \( \widehat F \) satisfies the following conditions. 

\medskip
\noindent
\textbf{A1)} \( F: \mathbb S^1 \to \mathbb S^1 \) is a $C^{2}$ uniformly expanding map, i.e. \( |F'(\omega)|\geq \theta \) for some \( \theta>1\) and all \( \omega \in \mathbb S^{1}\);

\medskip
\noindent
\textbf{A2)}  
The restrictions  \( T:  \mathbb{S}^{1}\times (-1,0)  \to (-1,1) \) and \( T:  \mathbb{S}^{1}\times (0,1)  \to (-1,1) \) are \( C^{2}\). Moreover, 
for each $\omega \in \mathbb{S}^{1}$  the corresponding restrictions  $T_{\omega}: (-1,0) \rightarrow (-1,1)$ and $T_{\omega}: (0,1) \rightarrow (-1,1)$ are orientation preserving $C^{2}$ diffeomorphisms such that the only fixed points are the endpoints of $I$ and such that \( \frac{\partial^{2}}{\partial x^{2}}T_{\omega}(x) \geq 0\) for \( x\in (-1,0) \) and \(\frac{\partial^{2}}{\partial x^{2}} T_{\omega}(x) \leq 0\) for \( x\in (0,1) \). 

\begin{figure}[h]
    \centering
    \includegraphics[width=4cm]{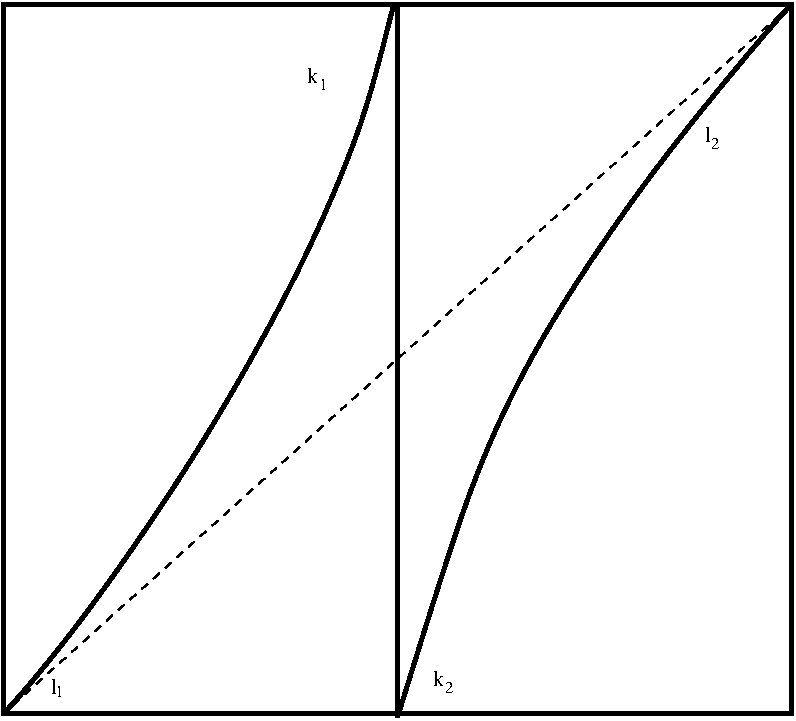}
  \caption{For each \( \omega \in \mathbb S^{1}\), the map \(T_{\omega}\) is a full branch map with two monotone branches}\label{fig:fig}
\end{figure}

\noindent
\textbf{A3)} There exist (one-sided) neighbourhoods $U_{-1}, U_{-0}, U_{+0},  U_{+1}$ of the endpoints of \( [-1,0]\) and \( [0,1]\) and two-variables \(C^{3}\) functions
 \[
 \eta^{(1)} : \mathbb S^{1} \times U_{-1}^{\circ}\to I,  \quad   
\eta^{(2)} : \mathbb{S}^{1} \times U_{+1}^{\circ}\to I, \quad 
 \xi^{(1)}: \mathbb{S}^{1} \times U_{-0}^{\circ}\to I, \quad 
 \xi^{(2)}: \mathbb{S}^{1} \times U_{+0}^{\circ}\to I,
 \]
with a \(C^{1}\) extension to the boundary of their respective domains, and $C^{2}$ functions 
\[
 \ell_{1},\ell_{2}: \mathbb{S}^{1}\rightarrow [0,1), \quad  k_{1}, k_{2} : \mathbb{S}^{1} \rightarrow (0, 1] \text{ and } a_{1}, a_{2}, b_{1}, b_{2}: \mathbb{S}^{1} \to (0,\infty),
\]
such that for all \(\omega\),  $T_{\omega}(U_{+0}) \subseteq U_{-1} $ and $T_{\omega}(U_{-0}) \subseteq U_{+1}$ and
\begin{equation}\label{eq:definition_of_Tomega}
T_{\omega}(x) = 
    \begin{cases}
         x+b_{1}(\omega)\big(x+1\big)^{\ell_{1}(\omega)+1}+ \eta_{\omega}^{(1)}(x) & \text{if }  x \in U_{-1}
         \\
       1-a_{1}(\omega)\big|x\big|^{k_{1}(\omega)} - \xi_{\omega}^{(1)}(x)& \text{if }   x \in U_{-0} \\
       -1+a_{2}(\omega) x^{k_{2}(\omega)}+ \xi_{\omega}^{(2)}(x) & \text{if } x\in U_{+0} \\
       x-b_{2}(\omega)\big(1-x\big)^{\ell_{2}(\omega)+1}-\eta_{\omega}^{(2)}(x) & \text{if } x \in U_{+1}.
    \end{cases}
\end{equation}

\medskip
\noindent
\textbf{A4)} 
For every \( \omega\in \mathbb S^{1}\), and for every \( x\in U_{-1}\) and \( x\in U_{+1}\) respectively, we have 
\begin{equation}\label{eq:neg_schwar_eta}
   \frac{\partial^{3} \eta_{\omega}^{(1)}}{\partial x^{3}}(x) \leq  0 \quad \text{ and } \quad  \frac{\partial ^{3} \eta_{\omega}^{(2)}}{\partial x^{3}}(x) \geq 0 
\end{equation}
Moreover, if \(\eta^{(1)}_{\omega}\) is not identically zero, and \(\ell_{1}(\omega) = 0\) then the first inequality in \eqref{eq:neg_schwar_eta} is strict, similarly if \(\eta^{(2)}_{\omega}\) is not identically zero, and \(\ell_{2}(\omega) = 0\) then the second inequality is strict.

\medskip
\noindent 
\textbf{A5)} There exist \(1>\delta>0\) such that for \(i=1,2\) and every \( \omega \in \mathbb S^{1}\) we have 
\begin{equation}\label{eq:expans_cond}
\delta<b_{i}(\omega) 
\quad \text{ and } \quad 
1+b_{i}(\omega)+\delta <\theta
\end{equation}
and,  letting, \( \bar\eta_{\omega}^{(1)}(x)= (1+x)^{\ell_{1}(\omega)+1}  \),  for 
every \( x \in U_{-1}^{\circ}\), and \(\alpha \in \{(i,j); i +j \leq 2\}\) we have 
\begin{equation}\label{eq:higherorder}
|\partial^{\alpha} \eta_{\omega}^{(1)}(x)| \leq \delta|\partial^{\alpha}  \bar\eta_{\omega}^{(1)}(x)|
\quad \text{ and } \quad 
 \partial^{\alpha} \eta_{\omega}^{(1)}(x)  \partial^{\alpha} \bar\eta_{\omega}^{(1)} (x)\geq 0,
 \end{equation}
where \(\partial^{\alpha}\) denote the mixed higher order partial derivatives with respect to \( \omega \) and \( x \) respectively. Similarly,  letting 
 \(
   \bar\eta_{\omega}^{(2)}(x)=(1-x)^{\ell_{2}(\omega)+1}, 
  \bar\xi_{\omega}^{(1)}(x) = |x|^{k_{1}(\omega)}, \) and \( 
    \bar \xi_{\omega}^{(2)} (x)= x^{k_{2}(\omega)},
 \) we suppose that the corresponding inequality \eqref{eq:higherorder} holds for each of these functions in relation to the corresponding functions  \( \eta_{\omega}^{(2)}, \xi_{\omega}^{(1)},  \xi_{\omega}^{(2)} \) for all \( x \) in their corresponding domains.

\medskip
 This completes our definition of the skew product system \eqref{eq:definnition_of_That} and we let \begin{equation*} 
 \mathcal{G} = \{  \widehat{F}: \mathbb{S}^{1} \times  I \rightarrow  \mathbb{S}^{1} \times I \text{ satisfying conditions (A1)-(A5)}\}.
\end{equation*}

 The \emph{doubly intermittent} interval maps satisfying (A2)-(A3) were introduced  in \cite{CLM23} (without the higher order terms \( \eta^{(i)}, \xi^{(i)}\)) as a generalization of the well-known \emph{intermittency maps}, such as those considered in  \cite{LSV99}. The \emph{doubly intermittent singular skew products} we introduce here are a corresponding generalization of the skew products considered in \cite{Gou07}. Our methods are also a generalization of those used in \cite{Gou07} albeit with some significant differences and some key technical innovations which we discuss in Section \ref{sec:context} together with  detailed historical background.
 
\begin{remark}\label{rem:derivative_ge1} 
Notice that \( \ell_{1}(\omega) \in (0,1) \) and/or \( \ell_{2}(\omega) \in (0,1) \) implies that the corresponding fixed points \( \pm 1 \)  for \( T_{\omega}\) are \emph{indifferent}, in the sense that \(\frac{\partial}{\partial x}T_{\omega}(-1)=1\) and/or \(\frac{\partial}{\partial x}T_{\omega}(+1)=1\). If 
$\ell_{1}(\omega)=0$ and/or $\ell_{2}(\omega)=0$ then the corresponding fixed points \( -1 \) and/or \( +1\) are expanding with derivative \( \frac{\partial}{\partial x}T_{\omega}(-1) \geq 1+b_{1}(\omega)\) and/or \( \frac{\partial}{\partial x}T_{\omega}(1) \geq 1+b_{2}(\omega)\). Also, \( k_{1}(\omega) \in (0,1) \) and/or \( k_{2}(\omega) \in (0,1) \) implies that the derivative on the left and/or the right of 0 respectively, is unbounded, see Figure \ref{fig:fig}. Notice also that (A2) and (A3) together imply that \( \frac{\partial}{\partial x}T_{\omega}(x)\geq 1 \) for all \( \omega \) and all \( x\in(-1,0)\cup(0,1) \). Indeed, by (A2) the derivative \( \frac{\partial}{\partial x}T_{\omega} \) is non-decreasing on \( (-1,0) \) and non-increasing on \( (0,1) \), and by (A3) and \eqref{eq:higherorder} it equals \( 1 \) or \( 1+b_{i}(\omega)\) at the fixed points \( \mp 1 \). We will use this fact repeatedly. \end{remark}

\begin{remark}
The assumption \eqref{eq:higherorder} says that \( \eta_{\omega}^{(1)} \) consist of \emph{higher order terms} in the first line of the expression for \( T_{\omega}\) in \eqref{eq:definition_of_Tomega} and that these higher order terms have the same sign as those of the main term \(  \bar\eta_{\omega}^{(1)}(x)= (1+x)^{\ell_{1}(\omega)+1}  \). 
The addition of the higher order terms significantly complicates some of the computations in the proofs but allows us to formally include many more maps, including the standard intermittency maps, as special cases. 
\end{remark}

 \begin{remark}\label{rem:expans_cond}
In (A2), the assumption that \( \frac{\partial^{2}}{\partial x^{2}}T_{\omega}(x) \geq 0\) for \( x\in (-1,0) \) and \( \frac{\partial^{2}}{\partial x^{2}}T_{\omega}(x) \leq 0\) for \( x\in (0,1) \) is not strictly necessary. It is  assumed here for simplicity and convenience, but could be replaced by a number of weaker but more technical assumptions.
Similarly, in (A5), the condition  \(1+b_{i}(\omega)+\delta<\theta\) is used only in the proofs of Lemma \ref{lem:expansion_contol_by_horizontal} and Proposition \ref{pro:hdistortion}, and there only in the case in which the corresponding exponent \( \ell_{i}\) vanishes for some \( \omega \), i.e. \( \ell_{i}^{-}=0 \). However, making this distinction would have required us to consider these as special cases in a number of lemmas and so we have decided to consider it a general assumption. 
\end{remark}

\subsection{Main Results}\label{sec:main_results}
Our interest in this paper is the statistical properties of skew products in the family \(  \mathcal{G} \). We begin with the first fundamental result concerned with the existence of a physical measure. 

\begin{maintheorem}\label{the:existence}
Every $\widehat{F} \in  \mathcal{G}$ admits a unique absolutely  continuous invariant probability measure \(\mu\). Moreover \(\mu\) is ergodic and equivalent to Lebesgue.
\end{maintheorem}

We then consider some of the dynamical and statistical properties with respect to the measure~\( \mu \). For  H\"older continuous observables $\varphi,\psi: \mathbb S^{1}\times I \to \mathbb R$ and \( n \geq 1 \) we define the 
 correlation function 
\begin{equation}\label{eq:def_correlation}
     Cor_{\mu}(\varphi, \psi \circ \widehat{F} ^{n}) = \left| \int \varphi(\psi \circ \widehat{F} ^{n})d\mu - \int\varphi d\mu \int \psi d\mu \right|.
\end{equation}
We say that $\mu$ has the property of \emph{decay of correlations} if $Cor_{\mu}(\varphi,\psi \circ \widehat{F}^{n}) \to 0$ as $n \to \infty$ and that the correlation function decays at a \emph{rate} \( \gamma_{n}\to 0\) (or that \( \mu \) is \emph{mixing} at a rate \( \gamma_{n}\to 0\)) if, for all H\"older continuous functions $\varphi,\psi$ there exists a constant $C_{\varphi,\psi}>0$ such that 
\begin{equation}\label{eq:corr}
Cor_{\mu}(\varphi,\psi \circ \widehat{F}^{n}) \leq C_{\varphi,\psi} \gamma_{n}
\end{equation}
 for every \( n \geq 1\). We then say that \( \mu \) is, for example, \emph{polynomially} or \emph{exponentially} mixing depending on whether \( \gamma_{n}\to 0 \) polynomially or exponentially. 
We will show that \( \mu \) is indeed always mixing, at a rate that depends on the parameters of the fibre maps \(T_{\omega}\) in \eqref{eq:definition_of_Tomega}. 
For $ \widehat{F} \in  \mathcal{G}$ and \( i=1,2\), we let 
\[
\ell_{i}^{+} := \max_{\omega\in \mathbb S^{1}}\{\ell_{i}(\omega)\} \quad \ell_{i}^{-} := \min_{\omega\in \mathbb S^{1}}\{\ell_{i}(\omega)\}, \quad \ell^{+} := \max\{\ell_{1}^{+},\ell_{2}^{+}\}, \quad \ell^{-} := \min\{\ell_{1}^{-},\ell_{2}^{-}\}.
\]
Similarly, we define \(k_{i}^{\pm},k^{\pm},a_{i}^{\pm},a^{\pm},b_{i}^{\pm}, b^{\pm}\). Finally,  we define
\begin{equation}\label{eq:beta}
\beta_{1} := \ell_{1}^{+}k_{2}^{+},
 \qquad \beta_{2} :=\ell_{2}^{+}k_{1}^{+},
\quad \text{ and } \quad  
\beta := \max\{ \beta_{1}, \beta_{2} \}.
\end{equation}
Notice that \(0\leq \beta<1\), moreover \( \beta=0\) implies that  \( \ell_{1}^{+}= \ell_{2}^{+}=0\), which in turn implies \( \ell_{1}(\omega)= \ell_{2}(\omega)=0\) for every \( \omega \in \mathbb S^{1}\), in which case \( T_{\omega}\) is uniformly expanding for every \( \omega \in \mathbb S^{1}\). On the other hand,  \( \beta > 0 \) implies that either \( \ell_{1}^{+}>0\) and/or \( \ell_{2}^{+}>0\) which means that there exists at least some \( \omega\in \mathbb S^{1}\) for which at least one of the two fixed points of the map \( T_{\omega}\) is neutral.   
It is not surprising therefore that we get different rates of decay of correlations in these two cases. 

\begin{maintheorem}\label{the:mixing}
Let $ \widehat{F} \in  \mathcal{G}$ and \( \mu \) the invariant probability measure equivalent to Lebesgue. 
If  \( \beta = 0 \), then \( \mu \) is exponentially mixing and 
if  \( 1>\beta > 0 \), then $\mu$ is polynomially mixing with rate $\gamma_{n}=n^{-\frac{1}{\beta}+1}$. 
\end{maintheorem}

It is worth saying a few words, albeit informally, about the interpretation of Theorems \ref{the:existence} and \ref{the:mixing} from the point of view of random dynamical systems discussed at the beginning of the introduction, referring to \cite{Arn98, KL06} and to Section \ref{sec:context} for the precise framework, which we do not need here. The measure \( \mu \) is a measure on the skew product, i.e. on pairs \( (\omega, x) \) consisting of a realization \( \omega \) of the noise and a point \( x \) on the fibre, and the correlations \eqref{eq:def_correlation} are integrals over all such pairs, i.e. over all realizations of the noise at the same time as over the fibre. In the terminology of random dynamical systems, Theorems \ref{the:existence} and \ref{the:mixing} are therefore statements about the \emph{annealed} (or averaged) statistical properties of the random compositions \( T^{n}_{\omega} \): Theorem \ref{the:existence} says that there is a unique absolutely continuous measure which describes the statistics of the fibre dynamics for Lebesgue almost every initial condition and almost every realization of the noise, and Theorem \ref{the:mixing} says that, averaging over the noise, the fibre dynamics ``forgets'' its initial condition at a rate which depends on the parameters of the fibre maps. Notice that the rate in Theorem \ref{the:mixing} depends only on the extremal values \( \ell^{+}_{i}, k^{+}_{i} \) of the parameter functions, i.e. on the ``slowest'' fibre maps which occur in the family. We do not know whether this rate is sharp, and there are in fact some reasons to believe that it is not: for random compositions of intermittent maps chosen independently at each step from a finite family it was shown in \cite{BBD14} that the annealed rate of decay of correlations is governed by the \emph{fastest} map in the family rather than by the slowest, and it would be interesting to know whether a similar phenomenon occurs in our setting, where the parameters vary continuously and the noise is generated by the base map \( F \), see Section \ref{sec:context}. Annealed statistical properties are the most basic and natural ones from the random point of view, but not the only ones. The corresponding \emph{quenched} statements, concerning the statistical properties of the compositions \( T^{n}_{\omega} \) for a \emph{fixed} typical realization \( \omega \) of the noise, are of a somewhat different nature and are not addressed in this paper, see Section \ref{sec:context} for some references and for a precise question in this direction.

\subsection{Innovations, and Open Questions}\label{sec:context}

One fundamental example of skew product system is the Alves-Viana family \cite{Via97} obtained by coupling a uniformly expanding circle map with a quadratic map. It was shown in \cite{Alv00} that these maps admit an absolutely continuous (with respect to the bi-dimensional Lebesgue measure) invariant probability measure. Another class of skew products that serves as the basis for the family consider in this paper was introduced in \cite{Gou07} and is obtained by coupling uniformly expanding maps with a particular class of intermittent maps \cite{LSV99}. More precisely, they are of the form   $T(\omega,x) = (F(\omega), T_{\alpha(\omega)}(x))$, where $F(\omega)= 4 \omega ~mod ~1$, $\alpha: \mathbb{S}^{1}: \to (0,1)$ is a $C^{2}(\mathbb{S}^{1})$ map, and  
\begin{equation}\label{eq:Talpha}
    T_{\alpha(\omega)}(x)=  \begin{cases}
         x(1+2^{\alpha(\omega)}x^{\alpha(\omega)}) & \text{if }  0\leq x \leq \frac{1}{2}
         \\
       2x-1 & \text{if }   \frac{1}{2} < x \leq 1.
    \end{cases}
\end{equation}
It was shown in \cite{Gou07} that the map \emph{T} admits a unique absolutely continuous ergodic invariant probability measure, which also satisfies a number of limit theorems, and in \cite{UV25} it was concluded that the invariant measure has polynomial decay of correlations. Several other papers study the same family or relatively minor generalizations. These include \cite{BB16} where \( \alpha \) is piecewise constant but can take any positive values and not just values in \( (0,1)\). A number of sophisticated limit theorems are proved in  \cite{EMV21} which also apply to some slightly more general classes of maps \eqref{eq:Talpha}.
In \cite{BQT21}, the base is replaced by the shift map on a symbolic space and the fibre maps are of the form \eqref{eq:Talpha}. In \cite{LV18} the base map is a shift and the fibre maps satisfy some non-uniform expansivity conditions.  A rather different class of skew products is that in which the fibre maps are rotations of the circle by an angle which depends on the base point, so that the fibre dynamics is isometric and the mixing properties of the skew product depend entirely on the base map and on the ``twist'' produced by the variation of the rotation angle. Exponential mixing and a local limit theorem for such systems, over a non-uniformly expanding base and with the rotation angle depending smoothly on the base point, were obtained in \cite{Gou09}, and exponential mixing was obtained in \cite{ButEsl17} when both the base map and  rotation angle are only piecewise smooth, and in \cite{But24} when the rotation angle is allowed to have singularities, with unbounded derivative, which produce an unbounded twist in the fibre direction. 

As mentioned above, the family of skew products we study in this paper have fibre maps as defined in  \eqref{eq:definition_of_Tomega} above, which  extend the maps \eqref{eq:Talpha}  (it is easy to see that the maps \eqref{eq:Talpha} can be obtained in the family \eqref{eq:definition_of_Tomega} by an appropriate choice of parameters).  We highlight here some of the key technical challenges and innovations which we will introduce in our proofs. 

\medskip
Notice first of all that the circles \(\mathbb S^{1} \times \{1,-1\}\) can be neutral curves for our maps, and therefore the inducing region \(\mathbf{H}\) needs to be  chosen away from these two curves. This makes our inducing scheme more sophisticated (see Section \ref{sec:markov_partition}) than that of \cite{Gou07}, where there was only one neutral curve and the inducing region was the upper half of the cylinder. The induced map \(T_{\mathbf{H}}\) is a sort of double inducing: first from the left neighbourhood (\(\mathbf{H}\)) of the singular curve \(\mathbb S^{1} \times \{0\}\) to the right, and then  back to the left. This is an adaptation of the one-dimensional construction in \cite{CLM23}.

\label{page:stair}
Secondly, as we shall see in Sections \ref{sec:eventually_uniformly_expansion} and \ref{sec:distortion}, the expansion and distortion estimates between two arbitrary points are obtained by introducing a sequence of intermediate points within the same partition element such that each consecutive pair shares either the same horizontal or the same vertical coordinate,  see in particular the ``\emph{Monotone Staircase Lemma}'' \ref{lem:stairs_between_points}. We first establish the desired estimates for each pair of consecutive points and then combine them to obtain the corresponding estimates for the original pair of points. This approach is the key innovation that allows us to construct a full-branch induced  Gibbs-Markov map, whereas the induced map in \cite{Gou07} is not full-branch but instead satisfies the weaker ``big image'' property.

Thirdly, in  \cite{Gou07} the uniformly expanding map \(F\) was an affine expanding map of degree \(4\) and the degree of the map was a rather strict requirement for technical reasons. In our setting, our counterpart of this assumption is the second condition in \eqref{eq:expans_cond}. Moreover, as we mentioned in Remark~\ref{rem:expans_cond}, this can be relaxed further if both fixed points are neutral for all  fibre maps (\(\ell^{-}_{i}>0\)).

\medskip
Notwithstanding these technical innovations, there are still several open questions and generalizations which would be very interesting to address. A first very natural question is how the invariant probability  measures given by Theorem \ref{the:existence} depend on the parameters.
\begin{question}
    Does the density of the absolutely continuous invariant measure depend continuously, or even differentiably, on the parameter functions?
\end{question}
Secondly, with some restriction on the parameter function  \(\alpha\), see \eqref{eq:Talpha}, the convergence of the limit laws of Birkhoff sums to a normal law or a stable law were established in \cite{Gou07}\todelete{ }. This requires some better estimates, particularly on the tail of the return time, which is not part of our study.

\begin{question}
Do our maps satisfy  limit laws analogous to those of \cite{Gou07}?
\end{question}

Thirdly, as we mentioned above, the interval maps satisfying (A2)-(A3), without the higher order terms  \(\eta^{(i)}\) and \(\xi^{(i)}\), were introduced in \cite{CLM23} as part of a larger family which included the possibility of having critical points instead of singularities, corresponding to having   \( k_{1}>1 \) and/or \( k_{2}>1\) in~\eqref{eq:definition_of_Tomega}. This case is not included in our study because we are unable to establish  the uniform expansion of the induced map, required to prove the  separation and distortion properties. In the one-dimensional case (even with the higher order terms), the expansion can be obtained relatively easily but the argument does not extend, at least not easily,  to the skew-product setting. 

\begin{question}
Can our techniques and results be generalized to  fibre maps \eqref{eq:definition_of_Tomega} with critical points? 
\end{question}

This may not be possible in full generality but should be true with some additional assumptions. 

\begin{conjecture}
The result remains true under the following assumptions. If \( \max_{\omega}k_{1}(\omega) > 1 \) then \emph{either} \( k_{1}\) \emph{or} \( \ell_{2}\) is constant, and if \( \max_{\omega}k_{2}(\omega) > 1 \) then \emph{either} \( k_{2}\) \emph{or} \( \ell_{1}\) is constant.  
\end{conjecture}

Another natural assumption, but possibly significantly more challenging from a technical point of view, might be to impose some averaging condition such as \( \int k_{1} dm < 1\) and \( \int k_{2} dm < 1\). 

\medskip
As discussed in the introduction, skew products over a chaotic base are closely related to random dynamical systems, a very large field with a long history and an extensive literature to which we cannot do justice here. The general theory, including precise definitions of the annealed and quenched statistical properties which we mentioned only informally in Section \ref{sec:main_results}, is developed in the monographs \cite{Kif86, Arn98} and in the survey \cite{KL06}. Random compositions of interval maps, and the existence of absolutely continuous invariant measures for them, were studied already in \cite{Oh83, Pel84, Mor85}, and the thermodynamic formalism and the spectral theory of transfer operators for random expanding and hyperbolic systems, including annealed and quenched decay of correlations, were developed in \cite{Kif93, BY93, BG95, Bal97, Buz99}, see also \cite{AA03} for random perturbations of non-uniformly expanding maps. More recently, a very effective approach to the quenched statistical properties of random uniformly expanding and hyperbolic systems has been developed on the basis of multiplicative ergodic theorems for transfer operator cocycles \cite{FLQ13, GQ14, DFGV18, AFGV21}, see also \cite{BRS20, DS23} for statistical stability and linear response, and a variety of annealed and quenched limit theorems and mixing results for random and sequential systems have been obtained in \cite{CR07, ALS09, ANV15, HNTV17, DH20, Su23}. We mention in particular \cite{Haf21} which, as in our setting, considers skew products whose base map is mixing rather than a Bernoulli shift. Closest to the present paper are the results on random compositions of intermittent maps of the form \eqref{eq:Talpha}, both with independent choices of the parameter \( \alpha \) at each step and in the more general sequential (non-stationary) setting, for which decay of correlations, loss of memory and limit theorems have been obtained in \cite{BBD14, AHNTV15, LS16, NTV18, BBR19, BQT21, NPT21, KL21}. As mentioned in Section \ref{sec:main_results}, a remarkable feature of these results is that both the annealed \cite{BBD14} and the quenched \cite{BBR19} rates of decay of correlations turn out to be governed by the fastest map in the family, which suggests the following.
\begin{question}\label{q:sharp}
Is the rate of decay of correlations in Theorem \ref{the:mixing} sharp? In particular, is the actual rate determined by the ``fastest'' fibre maps rather than by the slowest ones?
\end{question}

A further natural question, which is apparently more delicate, concerns the \emph{quenched} statistical properties of the fibre dynamics, i.e. the statistical properties of the compositions \( T^{n}_{\omega} \) for a fixed typical \( \omega \). There are several ways to formulate such properties precisely, and perhaps the simplest, which is the one used for example in \cite{CR07, AHNTV15, KL21} for sequential compositions of maps and which does not require any additional structure, is the following. Letting \( m \) denote Lebesgue measure on \( I \), we say that the compositions \( T^{n}_{\omega} \) have \emph{quenched decay of correlations}  at rate \( \gamma_{n} \to 0 \) if for all H\"older continuous \( \varphi, \psi: I \to \mathbb R \) there exists a constant \( C_{\omega, \varphi, \psi}>0 \) such that
\[
\left| \int \varphi \, (\psi \circ T^{n}_{\omega}) \, dm - \int \varphi \, dm \int \psi \circ T^{n}_{\omega} \, dm \right| \leq C_{\omega, \varphi, \psi} \, \gamma_{n}
\]
for every \( n \geq 1 \). Results of this kind for sequential and random compositions of intermittent maps were obtained in \cite{AHNTV15, KL21}, and quenched decay of correlations for random compositions of intermittent maps, with respect to a family of measures on the fibres which depends on the realization of the noise, was obtained in \cite{BBR19}.
\begin{question}\label{q:quenched}
Do the compositions \( T^{n}_{\omega} \) have, for Lebesgue almost every \( \omega \), quenched decay of correlations, and if so at what rate? In particular, is the rate the same as the rate of decay of the annealed correlations in Theorem \ref{the:mixing}?
\end{question}
Finally, we mention that a rather different mechanism of intermittency, in which the intermittent behaviour is produced by the random alternation of maps with different behaviour near a common fixed point (for example a contracting and an expanding map, or a map with a critical point and an expanding map) rather than by neutral fixed points of the individual maps, has been studied in \cite{GH17, HKRVZ22, KZ23, HPR24}.

\subsection{Overview of the proofs of the main results}\label{sec:overview}
We discuss here our overall strategy and prove our Theorems modulo two key technical propositions, which we then prove in the rest of the paper.   We assume throughout that $ \widehat{F} \in \mathcal{G}$. 
We start with the abstract definition of Gibbs-Markov map.
Let  $M$ be a measure space with a reference probability measure $m$. Let  $T$ : $M\rightarrow M$ be a measurable map and $\mathscr{R}$ be a countable  partition of $M$  (\emph{mod}\ \(0\)) into pairwise disjoint subsets and, for $n\geq 1$, let $\mathscr{R}_{n}:= \vee_{k=0}^{n-1}T^{-k}(\mathscr{R})$. The sequence $\{\mathscr{R}_{n}\}_{n \in \mathbb{N}}$ is a \emph{basis} of  \emph{M} if  each \( \mathscr R_{n}\) is a partition of \( M \) (\emph{mod}\ \(0\)), $\mathscr{R}_{\infty}:= $ $\vee_{n=0}^{\infty} T^{-n} (\mathscr{R})$ is a partition into single points, and $\{\mathscr{R}_{n}\}_{n\in \mathbb{N}}$ generates the Borel $\sigma$-algebra \( \mathcal{B}\) on  \( M \). 
A measurable function $J_{T} : M \to [0,\infty)$ is a \emph{Jacobian} of $T$ (with respect to \(m\)) if for all \(A \in \mathcal{B}\),
\( 
    m(T(A))= \int_{A}J_{T}~ dm
\). For any  \( x, y \in M\)  we define the \emph{separation time} 
\begin{equation}\label{eq:separation_time}
    s(x,y) := \min \{ n\geq 0:  \ T^{n}(x)\text{ and}  ~ T^{n}(y)~ \text{lie in different elements of}~ \mathscr{R} \}.
\end{equation}

\begin{definition}\label{def:gibbs_markov}
\emph{T} is a \emph{Gibbs-Markov} map (with respect to $\mathscr{R}$) if it satisfies the following properties:

\noindent 
\emph{Full Branch}:  \emph{T} maps each $\alpha \in \mathscr{R}$ bijectively to $M$ (\emph{mod} 0).\\
 \emph{Separability}: The sequence $\{\mathscr{R}_{n}\}_{n \in \mathbb{N}}$ is a basis of $M$.\newline
 \emph{Non-singularity}: \emph{T} has a strictly positive Jacobian $\emph{J}_{T}$.\newline
 \emph{Gibbs/Bounded distortion}: There are $C> 0$ and 0 $<\gamma< 1$ such that, for all $\alpha\in\mathscr{R}$ and $x,y\in \alpha$, 
\begin{equation}\label{eq:jacobian}
\ln \frac{{\emph{J}_{T}(x)}}{\emph{J}_{T}(y)} \leq C\gamma^{s(x,y)}.
\end{equation}
\end{definition}
By  standard results, see \cite{You98}, or \cite[Theorem 3.13]{Alv2020} , a Gibbs-Markov map $T$ admits a unique invariant, ergodic absolutely continuous (with respect to \emph{m}) probability measure, with density bounded above and below.

We cannot apply these classical results directly to our maps  $\widehat{F}$ defined in \eqref{eq:definnition_of_That} because they are not Gibbs-Markov,  since the neutral fixed points and the infinite derivatives mean that the bounded distortion conditions cannot hold. Therefore  we construct an \emph{induced} map that is Gibbs-Markov. The technique of constructing an induced Gibbs-Markov map is widely used but in  our setting, the existence of the two neutral curves makes the inducing highly nontrivial. To state our result, first of all  we consider a metric  on \(\mathbb{S}^{1}\times I\) by letting
\begin{equation}\label{eq:metricd}
    d\left((\omega_{1},x_{1}),(\omega_{2},x_{2})\right)=|\omega_{1}-\omega_{2}|+|x_{1}-x_{2}|.
\end{equation}
We then fix a region \(\mathbf{H} \subset \mathbb{S}^{1}\times I\), which we call the \emph{inducing region} on which we define the \emph{first return time} function, $\tau: \mathbf{H} \to \mathbb{N}$ by 
 \[ 
\tau(\omega,x):= \min \{ n>0: \ \widehat{F}^{n}(\omega,x) \in \mathbf{H} \} 
\]
 and the corresponding \emph{first return induced map}, $T_{\mathbf{H}}: \mathbf{H} \to \mathbf{H}$ by 
\begin{equation*}
 T_{\mathbf{H}}(\omega,x):=  \widehat{F}^{\tau(\omega,x)}(\omega,x).
\end{equation*}
Given such a first return induced map, we define two important properties of such maps. 
\begin{definition}\label{def:expanding}
A Gibbs-Markov  (with respect to a partition \(\mathscr{R}\)) induced map \(T_{\mathbf{H}}\), is \emph{expanding} if there are constants \(\Cl{Cexpanding}>0\) and \(0<\gamma<1\) such that for all \(\alpha \in \mathscr{R}\) and \((\omega_{1},x_{1}), (\omega_{2},x_{2}) \in \alpha \) we have
\begin{itemize}
      \item \(d\left(T_{\mathbf{H}}(\omega_{1},x_{1}),T_{\mathbf{H}}(\omega_{2},x_{2})\right)\leq \Cr{Cexpanding} \gamma^{s\left((\omega_{1},x_{1}),(\omega_{2},x_{2})\right)}\)
      \item \(d(\widehat{F}^{j}(\omega_{1},x_{1}),\widehat{F}^{j}(\omega_{2},x_{2}))\leq \Cr{Cexpanding} d(T_{\mathbf{H}}(\omega_{1},x_{1}),T_{\mathbf{H}}(\omega_{2},x_{2}))\),  for all \(0\leq j\leq \tau(\omega_{1},x_{1})=\tau(\omega_{2},x_{2})\).
\end{itemize}
\end{definition}

\begin{definition}\label{def:saturation}
A first return Gibbs-Markov induced map \( T_{\mathbf{H}}\)  \emph{saturates} the cylinder $\mathbb{S}^{1}\times I$ if 
\begin{equation*}
    \bigcup_{n \geq 0}\bigcup_{i=0}^{n-1}  \widehat{F}^{i} \big(\{\tau=n\} \big)= \bigcup_{n\geq 0} \widehat{F}^{n} \big(\{\tau>n\}\big) = \mathbb{S}^{1} \times I \quad (mod~0).
\end{equation*}
\end{definition}
Saturation ensures that before returning to the inducing region, the base dynamics goes through all the regions in the base space.  
We can now state our first main technical result.

\begin{proposition}\label{pro:gibbs_markov_structure}
Let
\begin{equation}\label{eq:inducing_region}
 \mathbf{H}:= \bigcup_{\omega \in \mathbb{S}^{1}}\{(\omega,y): \    x(\omega)<y <0 \text{ where } 
   x (\omega) \in(-1,0) \text{ satisfies }~T_{\omega}(x(\omega)) =0\}.
\end{equation}
Then the induced map $T_{\mathbf{H}} $ is a first return expanding Gibbs-Markov map that saturates $\mathbb{S}^{1} \times I$. Moreover, for every \( n \geq 2\), we have \( m\{\tau =n\}>0\), in particular the lowest common divisor  of the first return times is 1. 
\end{proposition}

By  Proposition \ref{pro:gibbs_markov_structure},  $T_{\mathbf{H}}$ admits a unique invariant ergodic probability measure $\hat{\mu}$ supported on~$\mathbf{H}$, which is absolutely continuous with respect to (indeed \emph{equivalent} to) Lebesgue measure $m$, and \(d\hat{\mu}/dm\) is bounded above and below. We \emph{spread} the measure in the whole cylinder by letting 
\begin{equation}\label{eq:spredingmu}
    \widetilde{\mu}: = \sum_{n=0}^{\infty}\widehat{F}^{n}_{*}\big(\Hat{\mu}|_{\{\tau> n\}} \big).
\end{equation}
Again, by standard results, $\widetilde{\mu}$ is a \emph{sigma-finite} measure, \(\widehat{F}\)-invariant, and ergodic. Using the non-singularity of $ \widehat{F}$, we  get that it is absolutely continuous with respect to Lebesgue and indeed,  since $T_{\mathbf{H}}$ saturates $\mathbb{S}^{1} \times I$,  it is in fact equivalent to Lebesgue. 
Notice that $ \widehat{F}^{-n}(\mathbb{S}^{1}\times I) = \mathbb{S}^{1} \times I$, therefore
        \begin{equation}\label{eq:tail}
     \widetilde{\mu}(\mathbb{S}^{1}\times I) =  \sum_{n=0}^{\infty} \Hat{\mu}( \widehat{F}^{-n}(\mathbb{S}^{1}\times I) \cap \{\tau > n\}) = \sum_{n=0}^{\infty} \Hat{\mu}((\mathbb{S}^{1}\times I) \cap \{\tau > n\}) = \sum_{n=0}^{\infty} \Hat{\mu}(\{\tau > n\}).
    \end{equation}
Notice, it is not obvious that the sum on the right side of \eqref{eq:tail} converges, and therefore we can not declare directly that measure $\widetilde{\mu}$ is a finite measure. To determine that it can be normalized to get a probability measure we use the following estimate on the tail of the return time.
\begin{proposition}\label{pro:tail_lebesgue}
For all $n \geq 1$, if \(\ell_{1}^{+} \neq 0\) and/or \(\ell_{2}^{+} \neq 0\), then
    \begin{equation}\label{eq:tail_lebesgue1}
        \hat{\mu}\big(\{\tau> n\}) \lesssim n^{-\frac{1}{\beta}},
    \end{equation}
and, if  $\ell_{1}^{+}=\ell_{2}^{+}=0$, then
\begin{equation}\label{eq:tail_lebesgue2}
    \hat{\mu}\big(\{\tau> n\}) \lesssim n(1+b^{-})^{\frac{-n}{k^{+}}} \lesssim e^{- \frac{\ln(1+b^{-})}{2k^{+}}n}.
\end{equation}
\end{proposition}

\begin{remark}
The notation \( \lesssim \) in  equations \eqref{eq:tail_lebesgue1} and \eqref{eq:tail_lebesgue2} means that the left hand side is bounded by a constant multiple of the right hand side for some constant \emph{independent of \( n\)}. Since our interest is in the asymptotic \emph{rates} of decay,  this is sufficient and significantly simplifies the notation, especially in the course of the technical proofs to be given below. 
In some situations, the specific value of the constants matters for the computations and in those cases we keep it explicit. 
\end{remark}

We now prove our main theorems assuming 
Propositions \ref{pro:gibbs_markov_structure}  and \ref{pro:tail_lebesgue}.

\begin{proof}[Proof of the Theorem \ref{the:existence} assuming Propositions \ref{pro:gibbs_markov_structure} and \ref{pro:tail_lebesgue}]
It follows from Proposition \ref{pro:gibbs_markov_structure} and the discussion above that 
the measure $\widetilde{\mu}$ defined in \eqref{eq:spredingmu}  is invariant, ergodic, and equivalent to Lebesgue. It just therefore remains to show that it is a finite measure so that it can be normalized to a probability measure. By plugging \eqref{eq:tail_lebesgue1} into \eqref{eq:tail}, we have if $\ell_{1}^{+}\neq 0$ or $\ell_{2}^{+} \neq 0$, then
\begin{equation*}
   \widetilde{\mu}(\mathbb{S}^{1} \times I) \lesssim 1+ \sum_{n=1}^{\infty}\frac{1}{n^{\frac{1}{\beta}}} <\infty,
\end{equation*}
and, if  $\ell_{1}^{+}=\ell_{2}^{+}=0$, then by \eqref{eq:tail_lebesgue2} and \eqref{eq:tail} we also similarly have $\widetilde{\mu}(\mathbb{S}^{1} \times I) <\infty$. 
In both cases we define \(\mu:= {\widetilde{\mu}}/{\widetilde{\mu}(\mathbb{S}^{1} \times I)}\) to get the required probability measure.
\end{proof}

\begin{proof}[Proof of Theorem \ref{the:mixing} assuming Propositions \ref{pro:gibbs_markov_structure} and \ref{pro:tail_lebesgue}]
By the classical results \cite{You99}, see also  \cite[Theorem 3.54]{Alv2020},  the decay of correlation follows from the fact that the lowest commond divisor of the return time is 1 and from decay of the tail of the inducing time. By  \eqref{eq:tail_lebesgue2}  if $\ell_{1}^{+}= \ell_{2}^{+}=0$ then \(\beta=0\) and we have an exponential decay of return times leading to exponential decay of correlations, and by~\eqref{eq:tail_lebesgue1} if $\ell_{1}^{+} \neq 0 $ and/or $\ell_{2}^{+} \neq 0$  there is a $1/\beta$-\emph{polynomial decay} on the return time leading to the corresponding correlation decay~rate. 
\end{proof}

In the rest of the paper we prove Propositions \ref{pro:gibbs_markov_structure} and \ref{pro:tail_lebesgue}.
In Section \ref{sec:markov_partition_and_size_estemats} we show that the first return map $T_{\mathbf{H}}$  is a full branch map, in Section \ref{sec:eventually_uniformly_expansion} we prove the non-singularity and separability, and  in Section~\ref{sec:distortion} we prove bounded distortion, thus completing the proof that $T_{\mathbf{H}}$ is a Gibbs-Markov map, as stated  in Proposition \ref{pro:gibbs_markov_structure}. Finally, in Section \ref{sec:tail}  we prove the tail estimates and   Proposition~\ref{pro:tail_lebesgue}.

\section{Markov Partition and their Size Estimates}\label{sec:markov_partition_and_size_estemats}

The first step in proving  $T_{\mathbf{H}}$ is a Gibbs-Markov map is the existence of the \emph{mod} \(0\)  partition \( \mathscr{R}\) so that the full branch property is satisfied.
In this section we will construct such a partition 
\begin{equation}\label{eq:induc_partition}
    \mathscr{R} := \{ \mathscr{A}_{n,k,r}; \ \text{for all integers }  n,k \geq 1 \text{ and }  1\leq r\leq q^{n+k}\}.
\end{equation}
In Section~\ref{sec:markov_partition} we construct the partition and prove  the following result.
\begin{proposition}\label{pro:full_branch_saturation}
$T_{\mathbf{H}}$ is a full branch  map on  \( \mathscr{R}\)  that saturates $\mathbb{S}^{1} \times I$ with \( \tau|_{\mathscr{A}_{n,k,r}} = n+k\).
\end{proposition}

 Topologically, the partition elements are ``rectangles''. We first partition the cylinder vertically into horizontal \emph{strips}, and then divide the strips horizontally into sub-regions.   
In Section~\ref{sec:size_estemats} we give some estimates on the sizes  of the partition elements which we will need in the following sections.

As we shall see from the construction, the indexes \( n, k \) have a very precise geometrical meaning related to the orbit of the partition element \(  \mathscr{A}_{n,k,r}\) in that they correspond to the time that the orbit spends on the ``left'' or ``right'' of some regions in the space before they return to the inducing region, and so in particular the return time is given exactly by the sum \( n+k\). For each pair \( n, k \geq 1 \) there are exactly \( q^{n+k}\) partition elements which share not only the same return time \( n+k \) but also the same ``combinatorics'', where \(q\) is the degree of the uniformly expanding map \( F \). Note that this allows us to count exactly the number of partition elements with a given return time and also with specific combinatorics, this will play an important role in our estimates.

\subsection{Topological construction of the Markov partition}\label{sec:markov_partition}
To construct the Markov partition \(\mathscr{R}\), we begin by vertically partitioning the cylinder using families of curves. Specifically, in Section \ref{sec:Auxi_partition}, we introduce auxiliary sequences of positive \(\{Y_{n,0}\}\), and negative \(\{Z_{0,k}\}\) functions. We show that the \emph{graphs} of the functions in each sequence are \emph{pairwise disjoint} and  respectively converge to the constant functions \(1\) and \(-1\), hence partition the upper and lower halves of the cylinder into vertical ``strips''. In particular the inducing region \(\mathbf{H}\) is the strip bounded by the graph of \(Z_{0,0}\) and the circle \(\mathbb{S}^{1} \times \{0\}\). The partition is further refined by the graphs of the two-parameter families of functions \({Z_{n,k}}\) and \({Y_{n,k}}\) with disjoint graphs, constructed as pullbacks of the auxiliary sequences. The graphs of the \(\{Z_{n,k}\}\) form the vertical partition of the inducing region into strips which we then further chop horizontally to get the Markov partition.

\begin{note}[Notation]
For simplicity, we will denote the inverse of the restriction of the intermittency map to  $I^{+}=[0,1]$ and $I^{-}=[-1,0]$ by $T^{-1}_{\omega}|_{I^{+}}$ and $T^{-1}_{\omega}|_{I^{-}}$ respectively. Moreover, for all $n\geq 1$ and $\omega \in \mathbb{S}^{1}$, we will write $T^{n}_{\omega}$ instead of $T_{F^{n-1}(\omega)}\circ T_{F^{n-2}(\omega)} \circ \dots \circ T_{\omega}$, and $T^{-n}_{\omega}|_{I^{\pm}}$ for the corresponding composition of the inverse branches $T^{-1}_{\omega}|_{I^{\pm}}\circ T^{-1}_{F(\omega)}|_{I^{\pm}}\circ\dots\circ T^{-1}_{F^{n-1}(\omega)}|_{I^{\pm}}$. 
\end{note}

\subsubsection{Auxiliary Partition}\label{sec:Auxi_partition}
We start by defining sequences of maps $\{Y_{n,0}\}$ and $\{Z_{0,k}\}$ from \(\mathbb S^{1}\) to \(I\), in the following way. First of all for \(\omega \in \mathbb S^{1}\) we define
\begin{equation*}
   Y_{0,0}(\omega):= T^{-1}_{\omega}|_{I^{+}}(0)  \quad \text{and} \quad Z_{0,0}(\omega) := T^{-1}_{\omega}|_{I^{-}}(0).
   \end{equation*}
Recall that \(T_\omega\) is smooth in both the variable \(\omega\) and \(x\) with non vanishing gradient, and therefore the level curves \(Y_{0,0}\) and $Z_{0,0}$ are smooth and $\widehat{F}(\omega, Y_{0,0}(\omega))=\widehat{F}(\omega, Z_{0,0}(\omega)) = (F(\omega),0)$. 
Inductively, for all \(n,k\geq 1\) and $ \omega \in \mathbb{S}^{1}$, we define
\begin{equation}\label{eq:Xn}
   Y_{n,0}(\omega) := T^{-1}_{\omega}|_{I^{+}}(Y_{n-1,0}(F(\omega))) \quad \text{and} \quad Z_{0,k}(\omega) := T^{-1}_{\omega}|_{I^{-}}(Z_{0,k-1}(F(\omega))).    
\end{equation}

In order to establish a well-defined partition of $\mathbb{S}\times I^{+}$ and $\mathbb{S}^{1} \times I^{-}$, we need to make sure that the graphs of different elements of the sequence $\{Y_{n,0}\}$ and $\{Z_{0,k}\}$ do not overlap. 

\begin{lemma}\label{lem:graphs_are_disjoint}
    The graphs of the maps in the sequence $\{Z_{0,k}\}$ are pairwise disjoint. Moreover, the sequence is pointwise monotone decreasing, and converges to the constant function \(-1\). Similarly, the graphs of $\{Y_{n,0}\}$ are disjoint, pointwise monotone increasing, and converge to \(1\).
 \end{lemma}
\begin{proof}
We prove the Lemma for the sequence $\{Z_{0,k}\}$ and the proof in case of the sequence $\{Y_{n,0}\}$ follows by an almost identical argument.
Let us first observe that the graphs of $Z_{0,0}$ and $Z_{0,1}$ are disjoint, by definition, for all \(\omega\) we have
\(
   Z_{0,0}(\omega) = T^{-1}_{\omega}|_{I^{-}}(0) 
   \) and 
   \(  Z_{0,1}(\omega)= T^{-1}_{\omega}|_{I^{-}} \circ T^{-1}_{F(\omega)}|_{I^{-}}(0).
\) Therefore, $Z_{0,0}(\omega)$ and $Z_{0,1}(\omega)$ are the images of two distinct points 0 and $T^{-1}_{F(\omega)}|_{I^{-}}(0)$ respectively, under a strictly increasing function $T^{-1}_{\omega}|_{I^{-}}$ and hence \(Z_{0,0}(\omega)>Z_{0,1}(\omega)\). Similarly
\(Z_{0,1}(\omega)= T^{-1}_{\omega}|_{I^{-}}(Z_{0,0}(F(\omega)) \) and \(Z_{0,2}(\omega)= T^{-1}_{\omega}|_{I^{-}}(Z_{0,1}(F(\omega))\). 
Since $Z_{0,0}(F(\omega))>Z_{0,1}(F(\omega))$, therefore $Z_{0,1}(\omega)>Z_{0,2}(\omega)$ . Inductively $Z_{0,k-2}(F(\omega))>Z_{0,k-1}(F(\omega))$, and therefore applying the monotone function $T_{\omega}^{-1}|_{I^{-}}$ to the last inequality we have $Z_{0,k-1}(\omega)>Z_{0,k}(\omega)$. 

Now let us establish the convergence, define  \(T_{\ell_{1}^{+}}(x):= x+b_{1}^{-}(1+x)^{\ell_{1}^{+}+1}\) in $U_{-1}$, and extend it continuously to \(I^{-}\), such that $x<T_{\ell_{1}^{+}}(x)<T_{\omega}(x)$ for all $\omega \in \mathbb{S}^{1}$ and $x\leq \widetilde{x}<0$, where $\widetilde{x}$ is such that $T_{\ell_{1}^{+}}(\widetilde{x})=0$. 
We inductively define a sequence of real numbers \(\{P_{k}\}\), by letting \(P_{0}=0\) and \(T_{\ell_{1}^{+}}(P_{k+1})=P_{k}\). Notice that for any \(x,y \in [-1,0]\) with \(y \leq x\) and for all \(\omega \in \mathbb{S}^{1}\), \(T^{-1}_{\omega}(y) \leq T^{-1}_{\ell_{1}^{+}}(x)\) which leads to \(Z_{0,k-1}(\omega)=T^{-k}_{\omega}(0) \leq T^{-k}_{\ell_{1}^{+}}(0)=P_{k}\). Therefore, it is enough to show that \(P_{k}\) converges to \(-1\). Again notice that \(P_{k}\) is monotone decreasing and bounded below by \(-1\). By the continuity of \(T_{\ell_{1}^{+}}\) and the definition of the sequence we get that the limit of \(P_{k}\) is a fixed point of \(T_{\ell_{1}^{+}}\), but the only fixed point of \(T_{\ell_{1}^{+}}\) in \([-1,0]\) is \(-1\), hence the limit is \(-1\).
\end{proof}

We refer to the region between the graphs of any two consecutive curves in the above sequences as \emph{strips}, precisely, 
\begin{equation*}
    \Hat{Y}_{0,0}~ \text{is the region bounded by the graph of} ~Y_{0,0} ~ \text{and} ~\mathbb{S}^{1} \times \{0\},
    \end{equation*}
and for all $n\geq 1$,
\begin{equation*}
\Hat{Y}_{n,0}~ \text{is the region between}~ Y_{n,0} ~\text{and}~Y_{n-1,0}.
\end{equation*}

\noindent Similarly,
\begin{equation*}
    \Hat{Z}_{0,0}~ \text{is the region between} ~Z_{0,0}~ \text{and}~ \mathbb{S}^{1} \times \{0\},
\end{equation*}
and for all $k\geq 1$,
\begin{equation*}
     \Hat{Z}_{0,k}~ \text{is the region between}~ Z_{0,k} ~\text{and}~Z_{0,k-1}.
\end{equation*}
Notice that \(\widehat{F}\) maps each strip to the one before, that is,
\begin{equation}\label{eq:pull_back1}
    \widehat{F}:\Hat{Y}_{0,0} \to \mathbb{S}^{1} \times I^{-} \quad\text{and} \quad \widehat{F}:\Hat{Z}_{0,0} \to \mathbb{S}^{1} \times I^{+},
\end{equation}
and  for all $n,k\geq 1$,
\begin{equation}\label{eq:pull_back2}
    \widehat{F}:\Hat{Y}_{n,0} \to \Hat{Y}_{n-1,0} \quad \text{and} \quad \widehat{F}:\Hat{Z}_{0,k} \to \Hat{Z}_{0,k-1}.
\end{equation}
These are surjective maps where the image wraps around the range $q$ times. Notice that for all \(k\geq 1\) the strips $\Hat{Z}_{0,k} \subset \mathbb{S}^{1} \times I^{-}$, therefore, by the first observation in \eqref{eq:pull_back1} we can pull the strips $\Hat{Z}_{0,k}$ back to  $\Hat
Y_{0,0}$, that is, for all $k\geq 1$ we define,
\begin{equation*}
    \Hat{Y}_{0,k}:= \widehat{F}^{-1}(\Hat{Z}_{0,k-1}) \cap\Hat{Y}_{0,0}.
\end{equation*}
Notice that $\Hat{Y}_{0,k} \subset \Hat{Y}_{0,0}$ and $\{\Hat{Y}_{0,k}\}$ form a partition of $\Hat{Y}_{0,0}.$ By the repeated application of  \eqref{eq:pull_back2}, we have $\widehat{F}^{n}: \Hat{Y}_{n,0} \to \Hat{Y}_{0,0}$ is surjective. Therefore, we can pull back the strips $\Hat{Y}_{0,k}$ to each $\Hat{Y}_{n,0}$, that is, for all $n,k\geq 1$ we define,
\begin{equation*}
    \Hat{Y}_{n,k}:= \widehat{F}^{-n}(\Hat{Y}_{0,k}) \cap\Hat{Y}_{n,0}.
\end{equation*}
Therefore, we have a refined \(mod~0\) partition of $\mathbb{S}^{1} \times I^{+}$, 
\begin{equation}\label{eq:partition_of_the_upper_half}
    \mathcal{Y}:=\{\Hat{Y}_{n,k}: \ n\geq 0,k\geq1\}.
\end{equation} 

\begin{remark}\label{rem:Ynk}
The strips $\Hat{Y}_{0,k}$ can also be seen as the closure of regions between the graphs of functions. First, for all \(n\geq 1\) and \(\omega\) we let
\begin{equation}\label{eq:Y_Z}
     Y_{0,1}(\omega):=T_{\omega}^{-1}|_{{I}^{+}}(Z_{0,0}(F(\omega))) \quad \text{and} \quad
    Y_{0,n}(\omega) := T^{-1}_{\omega}|_{I^{+}}(Z_{0,n-1}(F(\omega))),
\end{equation}
Notice that the sequence \(\{Y_{0,n}\}\) is pointwise decreasing and converges to the constant function \(0\). Then for all $ n, k \geq 1,$ we let
\begin{equation*}
    Y_{n,k}(\omega):=T^{-n}_{\omega}|_{I^{+}}(Y_{0,k}(F^{n}(\omega)).
\end{equation*}
By the monotonicity of \(T^{-1}_{\omega}\) and Lemma \ref{lem:graphs_are_disjoint} we conclude that the graphs of $Y_{n,k}$'s are mutually disjoint and the strip $\Hat{Y}_{n,k}$ is given by the region between the graphs of $Y_{n,k-1}$ and~$Y_{n,k}.$
\end{remark}

\subsubsection{Vertical partition of the Inducing Region}\label{sec:vertical_part_ind}
Notice that $\Hat{Z}_{0,0}$ is the inducing region $\mathbf{H}$, and by the second observation in \eqref{eq:pull_back1}, we can pull back the partition $\mathcal{Y}$ in \eqref{eq:partition_of_the_upper_half} to the inducing region $\mathbf{H}$, we define the following strips,
\begin{equation*}
    \Hat{Z}_{n,k}:= \widehat{F}^{-1}(\Hat{Y}_{n-1,k}) \cap\mathbf{H}.
\end{equation*}
Therefore, we get a \(mod~0\) partition of $\mathbf{H}$

\[\mathcal{Z}=\{\Hat{Z}_{n,k}: \quad n,k\geq1\}.\]
For all $n, k\geq 1$,  the iterates $1,\dots,n$ of the strip \(\hat Z_{n,k}\) lie in $\mathbb S^1\times I^+$, the iterates $n+1,\dots,n+k-1$ lie in $\hat Z_{0,k-1},\dots,\hat Z_{0,1}$, all disjoint from $\mathbf H$ and finally \(\hat{Z}_{0,1}\) is mapped to \(\mathbf{H}\). Therefore, we  conclude that the first return time of every $(\omega,x) \in \Hat{Z}_{n,k}$ is $n+k$, that is, for all $ n, k \geq 1$,
$$\tau|_{\Hat{Z}_{n,k}}:= n+k.$$
In each iteration, the strip $\Hat{Z}_{n,k}$ wraps around its image $q$ times. Therefore, the map $T_{\mathbf{H}}: \Hat{Z}_{n,k} \to \mathbf{H}$ is $q^{n+k}$ to $1$.
\begin{remark}\label{rem:Zn}
It will be useful for computation throughout the paper to group the strips $\Hat{Z}_{n,k}$ for a fixed value of $n$, that is 
\begin{equation*}
    \Hat{Z}_{n}:= \bigcup_{k\geq1} \Hat{Z}_{n,k}.
\end{equation*}
\end{remark}

\begin{remark}\label{rem:Znk}
Instead of pulling back the strips, we can pull back the graphs of $Y_{n,k}$, defined in Remark \ref{rem:Ynk}, that is, for all $n\geq 1, k\geq 0$
    \begin{equation}\label{eq:Z_Y}
        Z_{n,k}(\omega) := T^{-1}_{\omega}|_{I^{-}}(Y_{n-1,k}(F(\omega))).
    \end{equation}
Again, the sequence \(\{Z_{n,0}\}\) pointwise increases and converges to the constant function \(0\). The graphs of \(Z_{n,k}\)'s are disjoint and   $\Hat{Z}_{n,k}$ is the region between the graph of $Z_{n,k}$ and~$Z_{n,k-1}$.
\end{remark}

\subsubsection[The constant \( \bar n\)]{The constant \( \bar n\)}\label{sec:n_bar}

We now define a number \( \bar n \) after which  the curves defined above lie in the neighbourhoods \( \mathbb S^{1}\times U_{\pm 0}\) and \( \mathbb S^{1}\times U_{\pm 1} \), in which \( T_{\omega} \) has the explicit form \eqref{eq:definition_of_Tomega}. Throughout the proofs, we use this to separate the cases where an explicit expression is available from those where only bounds are known. 

Recall the map \(T_{\ell_{1}^{+}}(x)= x+b_{1}^{-}(1+x)^{\ell_{1}^{+}+1}\) and the sequence \(\{P_{n}\}\), the backward orbit of \( 0 \) under \( T_{\ell_{1}^{+}} \), defined in the proof of Lemma \ref{lem:graphs_are_disjoint}. Similarly, we define \(T_{\ell_{2}^{+}}(x) :=x-b_{2}^{-}(1-x)^{\ell_{2}^{+}+1}\) in \(U_{+1}\), extend it to \(I^{+}\), and define the sequence \(\{\tilde P_{n}\}\) as the backward orbit of \(0\) under \( T_{\ell_{2}^{+}} \). Letting \(T_{k_{1}^{+}}(x)=1-a_{1}^{-}|x|^{k_{1}^{+}}\) and \(T_{k_{2}^{+}}(x)= -1+a_{2}^{-}x^{k_{2}^{+}}\), we fix a natural number
\begin{equation}\label{eq:n_bar}
  \bar{n}:=\min\left\{i\geq 1;\quad T^{-1}_{k_{2}^{+}}(P_{i}) \in U_{+0} ~~\text{and}~~ T^{-1}_{k_{1}^{+}}(\tilde P_{i}) \in U_{-0} \right\}.
\end{equation}
Notice that for  \(n \ge \bar{n}\) and \(\omega \in \mathbb{S}^{1}\) we have \(Z_{0,n-1}(\omega) \in U_{-1},Y_{n-1,0}(\omega) \in U_{+1}, Y_{0,n}(\omega) \in U_{+0}\) and \(Z_{n,0}(\omega) \in U_{-0}\). Indeed, first notice that, \(Z_{0,n-1}(\omega) \leq P_{\bar{n}}\), then by the first relation in the definition of \(\bar{n}\) we have \( Z_{0,n-1}(\omega) \in  [-1,P_{\bar n}] \subset U_{-1}\). Similarly, by Remark \ref{rem:Ynk} we have \(Y_{0,n}(\omega)\leq Y_{0,\bar{n}}(\omega) =T^{-1}_{\omega}(Z_{0,\bar{n}-1}(F(\omega)) \leq T^{-1}_{k_{2}^{+}}(P_{\bar{n}}) \in U_{+0} \), therefore \(Y_{0,n}(\omega) \in U_{+0}\). An identical argument applies to  \(Y_{n-1,0}\) and \(Z_{n,0}\).

\subsubsection{Horizontal partition of the inducing region}
As established in the previous subsection, the induced map $T_{\mathbf{H}}$ is surjective from each strip $\Hat{Z}_{n,k} \in \mathcal{Z}$ to the inducing region $\mathbf{H}$, but it is not injective. We therefore subdivide the strip $\Hat{Z}_{n,k}$ into $q^{n+k}$ pieces (topological rectangles) such that $T_{\mathbf{H}}$ is injective, and hence bijective,  on each rectangle. 

Let $ \mathscr{P}_{0}$ be the \emph{mod} 0 partitions of the circle associated with the expanding map $F$, that is, for every \(\Omega \in \mathscr{P}_{0}\), the map \(F:\Omega \to \mathbb{S}^{1}\) is a smooth bijection (except at the boundary points). We then define the standard dynamical refinements of this partition, for every \( n, k\geq 1 \) we let 
\begin{equation}\label{eq:refinement}
\mathscr{P}_{n+k}=\bigvee_{i=0}^{n+k-1} F^{-i}(\mathscr{P}_{0}).
\end{equation}
We are now ready to define the sets $\mathscr{A}_{n,k,r}$ which form the partition \(  \mathscr{R}  \) as given in \eqref{eq:induc_partition}. 
To do this, we   divide the strip  $\Hat{Z}_{n,k}$ horizontally into sub-regions in such a way that for each partition element $\Omega \in \mathscr{P}_{n+k}$, there exists a topological rectangle in $\Hat{Z}_{n,k}$ such that its projection on the circle is $\Omega$ and vice versa. We represent a rectangle in this strip by $\mathscr{A}_{n,k,r}$, where $r$ corresponds to a partition element $\Omega_{r} \in \mathscr{P}_{n+k},$ and it ranges between $1$ and $q^{n+k}.$
This completes the construction of the partition \(  \mathscr{R}  \) and we can now  prove that it satisfies the required properties.

\subsubsection{Proof of Proposition \ref{pro:full_branch_saturation}}
\begin{proof}[Proof of Proposition \ref{pro:full_branch_saturation}]
It is clear from the above discussion that for all $n,k\geq 1$ and $1\leq r \leq q^{n+k}$, the map $\widehat{F}^{\tau}: \mathscr{A}_{n,k,r} \to \mathbf{H}$ is a bijection. For the saturation, first of all notice that
$$\{\tau >m\} = \bigcup_{n+k>m, 1\leq r \leq q^{n+k}} \mathscr{A}_{n,k,r} = \bigcup_{n+k> m} \Hat{Z}_{n,k} = \bigg(\bigcup_{k \geq m}\Hat{Z}_{1,k}\bigg) \bigcup \bigg(\bigcup_{n+k> m, ~n> 1} \Hat{Z}_{n,k}\bigg).$$ 
Notice that for all \(m>1\) we have $\widehat{F}^{m}(\Hat{Z}_{1,k})= \Hat{Z}_{0,k-m+1} $, therefore
\begin{equation*}
    \bigcup_{k> m}\Hat{Z}_{0,k-m} \subset \widehat{F}^{m}\bigg(\bigcup_{k \geq m}\Hat{Z}_{1,k}\bigg) \subset \widehat{F}^{m}(\{\tau>m\}).
\end{equation*}
Moreover, \( \mathbb{S}^{1} \times (0,1)= \widehat{F}(\{\tau>1\})\), and \(\hat{Z}_{0,0}=\widehat{F}^{0}(\{\tau>0\})\), therefore
\begin{equation*}
\mathbb S^{1} \times I~ (mod~ 0)= \bigcup_{k> m>1}\Hat{Z}_{0,k-m} \cup \hat{Z}_{0,0} \cup \left(\mathbb{S}^{1} \times (0,1)\right) \subset \bigcup_{m\geq 0} \widehat{F}^{m}(\{\tau >m\}).
\end{equation*}
\end{proof}

\subsection{Estimates for the size of the horizontal strips}\label{sec:size_estemats}

In this subsection, we estimate the size of the strips defined in Section \ref{sec:markov_partition}. In particular, we estimate the size of the region in which the graph of an element of the sequences defined in \eqref{eq:Xn},\eqref{eq:Y_Z} and \eqref{eq:Z_Y}, oscillates. In Lemma \ref{lem:graphs_are_disjoint}  we establish that the sequence of functions \( \{ Z_{0,n} \} \) pointwise decreases and converge to $-1$. Here we give some estimates for the convergence rate. Similarly, we give corresponding estimates for the rest of the sequences.
\begin{lemma}\label{lem:estimats_X}
  If $\ell_{1}^{\pm}, \ell_{2}^{\pm} > 0$,  then for all $n\geq 1$ and \(\omega \in \mathbb S^{1}\), we have 
    \begin{equation}\label{eq:X_n}
        (n+1)^{\frac{-1}{\ell_{1}^{-}}} \lesssim 1+ Z_{0,n}(\omega) \lesssim (n+1)^{\frac{-1}{\ell_{1}^{+}}}
    \end{equation}
    \begin{equation}\label{eq:Y_n}
        (n+1)^{\frac{-1}{\ell_{2}^{-}}} \lesssim 1-Y_{n,0}(\omega) \lesssim (n+1)^{\frac{-1}{\ell_{2}^{+}}} 
    \end{equation}
    \begin{equation}\label{eq:S_n}
         n^{\frac{-1}{\ell_{1}^{-}k_{2}^{-}}}\lesssim Y_{0,n}(\omega) \lesssim n^{\frac{-1}{\ell_{1}^{+}k_{2}^{+}}}
    \end{equation}
    \begin{equation}\label{eq:Z_n'}
         n^{\frac{-1}{\ell_{2}^{-}k_{1}^{-}}}\lesssim |Z_{n,0}(\omega)| \lesssim n^{\frac{-1}{\ell_{2}^{+}k_{1}^{+}}} 
    \end{equation}
If $\ell_{1}^{+} =0$ then
\begin{equation}\label{eq:ell1=0}
\begin{split}
   &\left(1+b_{1}^{+}+\delta\right)^{-n-1} \lesssim 1+Z_{0,n} \lesssim \big(1+b_{1}^{-}\big)^{-n-1}\\
    &\big(1+b_{1}^{+}+\delta \big)^{-\frac{n}{k_{2}^{-}}} \lesssim Y_{0,n} \lesssim \big(1+b_{1}^{-}\big)^{-\frac{n}{k_{2}^{+}}}
    \end{split}
\end{equation}
If \(\ell_{2}^{+}=0\) then the estimates in \eqref{eq:Y_n} and \eqref{eq:Z_n'} change similar to \eqref{eq:ell1=0}.
If \(\ell_{1}^{-}=0\) and \(\ell_{1}^{+} \neq 0\) then only the left side of \eqref{eq:X_n} and \eqref{eq:S_n} is replaced by \(\left(1+b_{1}^{+}+\delta\right)^{-n-1}\) and \(\left(1+b_{1}^{+}+\delta\right)^{-\frac{n}{k_{2}^{-}}}\) respectively. Similarly, if \(\ell_{2}^{-}=0\) and \(\ell_{2}^{+} \neq 0\), then the changes occur only on the left side of \eqref{eq:Y_n} and \eqref{eq:Z_n'}.
\end{lemma}

\begin{remark}
The asymptotic relations in \eqref{eq:X_n} depend on $(\ell_{1}^{\pm},b_{1}^{\pm},\delta)$, in \eqref{eq:Y_n} on $(b_{2}^{\pm},\ell_{2}^{\pm},\delta)$, in \eqref{eq:S_n} on $(\ell_{1}^{\pm},k_{2}^{\pm},b_{1}^{\pm},a_{2}^{\pm},\delta)$, and in  \eqref{eq:Z_n'} on $(k_{1}^{\pm},a_{1}^{\pm},b_{2}^{\pm},\ell_{2}^{\pm},\delta)$.
Notice that the dependence of the above inequalities on the corresponding parameters comes from their definition, for example, $Z_{0,n}$ is defined with the help of $T_{\omega}|_{I^{-}}$, which contains the parameter $\ell_{1}$ and $b_{1}$.
\end{remark}

\begin{proof}[Proof of Lemma \ref{lem:estimats_X}]
We provide proofs of \eqref{eq:X_n}, \eqref{eq:S_n}, \eqref{eq:ell1=0}, and the rest follow from an identical argument.  Moreover, it is sufficient to prove the estimates for \( n\geq \bar n \), since for the finitely many \( 1\leq n<\bar n \) all the quantities involved are bounded above and below by positive constants, recall Lemma~\ref{lem:graphs_are_disjoint} and Remarks \ref{rem:Ynk} and \ref{rem:Znk}, and so the estimates hold for all \( n\geq 1 \). Recall that \(Z_{0,n-1}(\omega) \leq P_{n}\), therefore, we focus on getting an upper bound for the sequence \(\{P_{n}\}\), which will lead us to our desired result.
For \(n \geq \bar{n}\) by definition of \(P_{n+1}\) we have \(P_{n+1}+b^{-}_{1}(1+P_{n+1})^{\ell_{1}^{+}+1}= P_{n}\), adding \(1\) to both sides and letting \(R_{n}=P_{n}+1\) we have 
\begin{equation*}
    R_{n+1}+b_{1}^{-}(R_{n+1})^{\ell_{1}^{+}+1}= R_{n}, \quad \text{therefore} \quad
    R_{n+1}^{-\ell_{1}^{+}}\left(1+b_{1}^{-}R_{n+1}^{\ell_{1}^{+}}\right)^{-\ell_{1}^{+}}=R_{n}^{-\ell_{1}^{+}}
\end{equation*}
Using binomial expansion and rearranging
\begin{equation*}
    R_{n+1}^{-\ell_{1}^{+}}\left(1-\ell_{1}^{+}b_{1}^{-}R_{n+1}^{\ell_{1}^{+}}+ O\left(R_{n+1}^{2\ell_{1}^{+}}\right)\right)= R_{n}^{-\ell_{1}^{+}},\quad \text{therefore} \quad R_{n+1}^{-\ell_{1}^{+}}-R_{n}^{-\ell_{1}^{+}} =  b_{1}^{-}\ell_{1}^{+} +o(1).
\end{equation*}
Applying the above relation recursively we get \(R^{-\ell_{1}^{+}}_{n+1}=n\ell_{1}^{+}b_{1}^{-}+o(n)\), that is \(R^{-\ell_{1}^{+}}_{n+1}=  n(\ell_{1}^{+}b_{1}^{-}+o(1))\), therefore, \(P_{n}+1 \lesssim   n^{\frac{-1}{\ell_{1}^{+}}}\). For the left-hand side of the inequality, we define \(T_{\ell_{1}^{-}}(x):=x+(b_{1}^{+}+\delta)(1+x)^{\ell_{1}^{-}+1}\) on \(U_{-1}\), and extend this map so that \(T_{\omega}(x)\leq T_{\ell_{1}^{-}}(x)\) for every \(\omega\in\mathbb{S}^{1}\) and \(x\leq T_{\ell_{1}^{-}}^{-1}(0)\). With this extension, the \(n\)-th point in the backward orbit of \(0\) under \(T_{\ell_{1}^{-}}\) is bounded above by \(Z_{0,n-1}(\omega)\) for every \(\omega\). Repeating the argument used above to estimate the backward orbit of \(0\) therefore gives the desired estimates.

Now, the first estimates in \eqref{eq:ell1=0}, if \(\ell_{1}^{+} =0\), we consider \(T_{\ell_{1}^{+}}(x):= x+b_{1}^{-}(1+x)\), with \(P_{0}=0\) and \(P_{n}= T_{\ell_{1}^{+}}(P_{n+1})\). Notice that \(T_{\ell_{1}^{+}}^{n+1}:[-1,P_{n+1}] \to [-1,0]\) is a bijection, furthermore, \(1=P_{0}+1= T_{\ell_{1}^{+}}^{n}(P_{n})-(-1)= T^{n+1}_{\ell_{1}^{+}}(P_{n+1})-(-1)\).
Notice \(T^{'}_{\ell_{1}^{+}}=1+b_{1}^{-}\), therefore by the mean value theorem we have, \(P_{n+1}+1\leq (1+b_{1}^{-})^{-n-1}\), hence the right side of the relation. For the left side, we use \(T_{\ell_{1}^{-}}(x) = x+(b_{1}^{+}+\delta)(1+x) \).

Now, let us prove \eqref{eq:S_n}, by the definition of \(Y_{0,n}\), we have \(T_{\omega}(Y_{0,n}(\omega))= Z_{0,n-1}(F(\omega))\), that is \(a_{2}(\omega)(Y_{0,n}(\omega))^{k_{2}(\omega)}+ \xi_{\omega}^{(2)}(Y_{0,n}(\omega))=1+Z_{0,n-1}(F(\omega))\). Furthermore, by assumption \eqref{eq:higherorder} we have \(0\leq \xi_{\omega}^{(2)}(Y_{0,n}(\omega))\leq \delta {(Y_{0,n}(\omega))^{k_{2}(\omega)}}\), therefore
\begin{equation}\label{eq:estemat1}
    \frac{1}{(a_{2}(\omega)+ \delta) }\left(1+Z_{0,n-1}(F(\omega))\right)\leq (Y_{0,n}(\omega))^{k_{2}(\omega)}\leq \frac{1}{a_{2}(\omega) }\left(1+Z_{0,n-1}(F(\omega))\right).
\end{equation}
Recall by definition \(a^{-}>0\), therefore using the estimates for \(Z_{0,n-1}\) from \eqref{eq:X_n} (respectively, from the first estimates in \eqref{eq:ell1=0}) in \eqref{eq:estemat1}, we have \eqref{eq:S_n} (respectively, the second estimates in \eqref{eq:ell1=0}).
\end{proof}

\section{Separability, Non-Singularity, and Expansivity}\label{sec:eventually_uniformly_expansion}
In Section \ref{sec:markov_partition_and_size_estemats} we proved the first property of Gibbs-Markov maps for the induced map $T_{\mathbf{H}}$, namely that  it is \emph{full branch} on the partition \( \mathscr{R}\). We now address the second and third properties,  \emph{non-singularity} and \emph{separability}, as well as the fact that it is an \emph{expanding} induced map in the sense of Definition \ref{def:expanding}.
We will deduce these  from a stronger expansivity property which uses the metric and differentiable  properties of our maps. The main technical result of this section is the following. 
\begin{proposition}\label{pro:expansion}
There exist  $ \Cr{Cveh}, \lambda, \tilde\lambda >1$  such that     for all $n,k,r$ and 
$(\omega, x) \in \mathscr{A}_{n,k,r}$
    \begin{equation}\label{eq:exp12}
        \bigg|\frac{\partial}{\partial x}T^{\tau(\omega,x)}_{\omega}(x) \bigg| >\widetilde{\lambda},
    \end{equation}
for all $(\omega_{1},x_{1}),(\omega_{2},x_{2}) \in \mathscr{A}_{n,k,r}$ and \(0\leq p< \tau:=n+k\),
\begin{equation}\label{eq:distortion4}
    |T^{\tau}_{\omega_{1}}(x_{1})-T^{\tau}_{\omega_{2}}(x_{2})| \geq \lambda |T_{\omega_{1}}^{p}(x_{1})-T_{\omega_{2}}^{p}(x_{2})|-2\Cr{Cveh}|F^{\tau}(\omega_{1})- F^{\tau}(\omega_{2})|,
\end{equation}
and  for all $(\omega_{1},x_{1}),(\omega_{2},x_{2}) \in \mathscr{A}_{n,k,r}$ and 
\( 1\leq m\leq s((\omega_{1},x_{1}),(\omega_{2}, x_{2}))\),
\begin{equation}\label{eq:eventual_expansion}
    d(T^{m}_{\mathbf{H}}(\omega_{1},x_{1}),T^{m}_{\mathbf{H}}(\omega_{2},x_{2})) \gtrsim \lambda^{m} d((\omega_{1},x_{1}),(\omega_{2},x_{2})).
\end{equation}
\end{proposition}

Recall that \( s((\omega_{1},x_{1}),(\omega_{2}, x_{2}))\) is the \emph{separation time} defined in \eqref{eq:separation_time}, which is the number of iterates of the induced map such that the images of the two points continue to fall inside the same element of the partition \( \mathscr{R}\). The estimate in \eqref{eq:eventual_expansion} is weaker than saying that the first iterate of the induced map is uniformly expanding  but it says that we do eventually have uniform expansion.

\begin{remark}\label{rem:vertical_expansion}
Notice that \eqref{eq:exp12} only implies the vertical expansion of  points with the \emph{same horizontal coordinates} and therefore does not  imply \eqref{eq:eventual_expansion}. To visualize this, choose two points in the inducing region with different vertical and horizontal coordinates, $(\omega_{1},x_{1}), (\omega_{2},x_{2}) \in \mathbf{H}$. If we iterate these two points, then each time the vertical coordinates are iterated under different fibre maps. Therefore, it is not guaranteed that even after the first application of the skew product $\widehat{F}$ they move far apart from each other. In fact, they might come close to each other vertically, and therefore once they return to the inducing domain, their separation is not guaranteed. As we shall see, the proof of \eqref{eq:exp12}  is relatively straightforward, see Lemma \ref{lem:vexpansion}, whereas the proof of \eqref{eq:eventual_expansion} is quite involved, requiring various auxiliary lemmas and  non-trivial geometric constructions  in Sections \ref{sec:techlem} and \ref{sec:proof_of_pro_expansion}. The inequality \eqref{eq:distortion4} is a key step in the proof of \eqref{eq:eventual_expansion} and  we prove it in Section \ref{sec:proof_of_pro_expansion}.
\end{remark}

Before proving Proposition \ref{pro:expansion}, we first show in the next three lemmas that it implies the required non-singularity, separability, and expansivity conditions.

\begin{lemma}\label{lem:nonsingularity}
    The induced map $T_{\mathbf{H}}$ is non-singular.
\end{lemma}

\begin{proof}
Notice that in this case the Jacobian $J_{T_{\mathbf{H}}}$ is given by the absolute value of the determinant of the differential of \(T_{\mathbf{H}}\), moreover the differential of the induced map $T_{\mathbf{H}}$ is a lower triangular matrix. 
\begin{equation}\label{eq:det_mat}
        DT_{\mathbf{H}}(\omega,x)= \begin{bmatrix}
(F^{n+k})^{'}(\omega)& 0 \\
\frac{\partial}{\partial \omega } T^{n+k}_{\omega}(x) & \frac{\partial}{\partial x } T^{n+k}_{\omega}(x)
\end{bmatrix}
\quad \text{and} \quad detDT_{\mathbf{H}}(\omega,x)= (F^{n+k})^{'}(\omega) \frac{\partial}{\partial x } T^{n+k}_{\omega}(x).
\end{equation}
Therefore, by \((A1)\) and \eqref{eq:exp12}, for all \(n,k,r\) and $(\omega,x) \in \mathscr{A}_{n,k,r}$ we have \( J_{T_{\mathbf{H}}}(\omega,x)= |detDT_{\mathbf{H}}(\omega,x)|\geq \theta^{n+k} \widetilde{\lambda}>0\). Hence we have non-singularity. 
\end{proof}

We now prove the separability property, namely that the sequence $\mathscr{R}_{n}= \vee_{k=0}^{n-1}T^{-k}_{\mathbf{H}}(\mathscr{R})$ generated by the partition $\mathscr{R}$   is a basis on the inducing domain $\mathbf{H}$. For a measurable set \(A \subset \mathbb{S}^{1} \times I\) we define the \emph{diameter} of \(A\) as
\( 
    diam(A):=\max\{d((\omega_{1},x_{1}),(\omega_{2},x_{2});  (\omega_{1},x_{1}),(\omega_{2},x_{2}) \in A \}
\) and the \emph{diameter of the partition} as \(diam(\mathscr{R}_{n}):=\max\{diam(\Omega_{n});  \ \Omega_{n}\in \mathscr{R}_{n}\}\).
 
\begin{lemma}\label{lem:basis}
    $\{\mathscr{R}_{n}\}_{n\in \mathbb{N}}$ is a basis on $\mathbf{H}.$
\end{lemma}

\begin{proof}
It is enough to show that \(diam(\mathscr{R}_{n}) \to 0\), see \cite[Corollary 2.5]{Alv2020} . Let \((\omega_{1},x_{1}), (\omega_{2},x_{2}) \in~\mathbf{H}\) be such that \((\omega_{1},x_{1}),(\omega_{2},x_{2}) \in \Omega_{n} \in \mathscr{R}_{n}\). Then the forward iterates of these points lies in a same element of the partition up to time \(n\geq 1\) and, therefore, by Proposition \ref{pro:expansion},
    \begin{equation*}
        d((\omega_{1},x_{1}),(\omega_{2},x_{2})) \lesssim \frac{1}{\lambda^{n}}d(T^{n}_{\mathbf{H}}(\omega_{1},x_{1}),T^{n}_{\mathbf{H}}(\omega_{2},x_{2})) \lesssim \frac{1}{\lambda^{n}}.
    \end{equation*}
Hence \(diam(\Omega_{n}) \to 0\) as \(n \to \infty\). 
\end{proof}

\begin{lemma}\label{lem:inducedexpanding}
The induced map \( T_{\mathbf{H}} \) is expanding.
\end{lemma}

\begin{proof}
Let \( \mathscr{A}_{n,k,r}\in\mathscr R \) and \( (\omega_{1},x_{1}),(\omega_{2},x_{2})\in\mathscr{A}_{n,k,r}\), let \( \tau=n+k \) be their common return time and let \( N:=s((\omega_{1},x_{1}),(\omega_{2},x_{2})) \) be their separation time. By the definition \eqref{eq:metricd} of the metric,
\begin{equation}\label{eq:dTH}
d(T_{\mathbf H}(\omega_{1},x_{1}),T_{\mathbf H}(\omega_{2},x_{2})) = |F^{\tau}(\omega_{1})-F^{\tau}(\omega_{2})| + |T^{\tau}_{\omega_{1}}(x_{1})-T^{\tau}_{\omega_{2}}(x_{2})|.
\end{equation}
We first verify the second condition in Definition \ref{def:expanding}, which is where \eqref{eq:distortion4} is used. For \( j=\tau \) there is nothing to prove, so let \( 0\leq j<\tau \). Rearranging \eqref{eq:distortion4} with \( p=j \), and using the fact that \( 2\Cr{Cveh}>1 \) and \eqref{eq:dTH}, we get
\[
|T^{j}_{\omega_{1}}(x_{1})-T^{j}_{\omega_{2}}(x_{2})| \leq \frac{1}{\lambda}\Big(|T^{\tau}_{\omega_{1}}(x_{1})-T^{\tau}_{\omega_{2}}(x_{2})| + 2\Cr{Cveh}|F^{\tau}(\omega_{1})-F^{\tau}(\omega_{2})|\Big) \leq \frac{2\Cr{Cveh}}{\lambda}\, d(T_{\mathbf H}(\omega_{1},x_{1}),T_{\mathbf H}(\omega_{2},x_{2})).
\]
Moreover \( \omega_{1},\omega_{2} \) lie in the same element of \( \mathscr P_{\tau} \), on which \( F^{\tau} \) is injective, and so by the uniform expansion of \( F \) and \eqref{eq:dTH} we have
\[
|F^{j}(\omega_{1})-F^{j}(\omega_{2})|\leq \theta^{-(\tau-j)}|F^{\tau}(\omega_{1})-F^{\tau}(\omega_{2})| \leq d(T_{\mathbf H}(\omega_{1},x_{1}),T_{\mathbf H}(\omega_{2},x_{2})).
\]
Adding the two inequalities we get
\[
d(\widehat F^{j}(\omega_{1},x_{1}),\widehat F^{j}(\omega_{2},x_{2})) \leq \Big(1+\frac{2\Cr{Cveh}}{\lambda}\Big) d(T_{\mathbf H}(\omega_{1},x_{1}),T_{\mathbf H}(\omega_{2},x_{2})),
\]
which is the second condition.

For the first condition, since the metric \( d \) is bounded, it is sufficient to consider \( N\geq 2 \). Then \( T_{\mathbf H}(\omega_{1},x_{1}) \) and \( T_{\mathbf H}(\omega_{2},x_{2}) \) lie in the same element of \( \mathscr R \) and their separation time is \( N-1 \), and so \eqref{eq:eventual_expansion} applied to these two points with \( m=N-1 \) gives
\[
d(T^{N}_{\mathbf H}(\omega_{1},x_{1}),T^{N}_{\mathbf H}(\omega_{2},x_{2})) \gtrsim \lambda^{N-1}\, d(T_{\mathbf H}(\omega_{1},x_{1}),T_{\mathbf H}(\omega_{2},x_{2})).
\]
Since the left hand side is bounded by the diameter of \( \mathbb S^{1}\times I \), this gives \( d(T_{\mathbf H}(\omega_{1},x_{1}),T_{\mathbf H}(\omega_{2},x_{2}))\lesssim \lambda^{-N} \), which is the first condition with \( \gamma=1/\lambda<1 \). Taking \( \Cr{Cexpanding} \) to be the larger of the two constants completes the proof.
\end{proof}

In Section \ref{sec:techlem} we state and prove three independent technical lemmas which will come together in Section \ref{sec:proof_of_pro_expansion} in the proof  of Proposition \ref{pro:expansion}.

\subsection{Three Technical Lemmas}\label{sec:techlem}

The proof of Proposition \ref{pro:expansion} needs three independent auxiliary technical lemmas which will also play a role in the distortion estimates in Section \ref{sec:distortion}. The first is a simple but crucial result which gives some expansion within the fibres in the ``vertical'' direction and immediately implies 
 \eqref{eq:exp12}. Recall from the construction that for all $(\omega,x) \in \mathscr{A}_{n,k,r}$ we have $(\omega,x) \in \Hat{Z}_{n}$ and  $(F^{n}(\omega), T^{n}_{\omega}(x)) \in \Hat{Y}_{0,k}$. We will therefore get  \eqref{eq:exp12} from independent expansion estimates in these two regions.

\begin{lemma}\label{lem:vexpansion}
There exists $\widetilde{\lambda}>1$ such that for all $n, k\geq 1$, \(0\leq i<n\), \(0\leq j<k\), $(\omega,x) \in \Hat{Z}_{n}$ and  $(\tilde \omega,\tilde x) \in \Hat{Y}_{0,k}$ we have 
\begin{equation}\label{eq:exp11}
    \big|\frac{\partial}{\partial x}T^{n-i}_{F^{i}(\omega)}(T^{i}_{\omega}(x)) \big| \geq \widetilde{\lambda}, \quad \text{and} \quad \big|\frac{\partial}{\partial x}T^{k-j}_{F^{j}(\tilde \omega)}(T^{j}_{\tilde \omega}(\tilde x)) \big| \geq \widetilde{\lambda}.
\end{equation}
In particular, we have \eqref{eq:exp12} and for all \( m\geq 1 \), \( (\omega,x)\in\Hat{Y}_{m,0} \) and \( (\tilde \omega,\tilde x)\in\Hat{Z}_{0,m} \) 
we have 
\begin{equation}\label{eq:vexpansion1}
 \frac{\partial}{\partial x}T^{m}_{\omega}(x)\geq\tilde\lambda
 \quad \text{ and } \quad 
  \frac{\partial}{\partial x}T^{m}_{\tilde \omega}(\tilde x)\geq\tilde\lambda.
  \end{equation}
\end{lemma}

The second technical result  is the ``\emph{Monotone Staircase Lemma}'',  the key technical step mentioned on page~\pageref{page:stair} of the introduction, which helps us to get the expansion and distortion estimates required for the full branch induced map. We will show that for all $n,k,r$ and for any two points $(\omega_{1},x_{1}),(\omega_{2},x_{2}) \in ~\mathscr{A}_{n,k,r}$, with $x_{1}>x_{2}$, there is a path between these two points in the form of a staircase in such a way that  all the \emph{vertices} of the staircase (but not necessarily all the edges) lie inside $\mathscr{A}_{n,k,r}$, see Figure \ref{fig:stairs_path}. 

\begin{lemma}[Monotone Staircase Lemma]\label{lem:stairs_between_points}
For all \(n\geq 1\) and for all pairs $(\omega_{1},x_{1}),(\omega_{2},x_{2})$ in the interior of $\Hat{Z}_{n}$ (and similarly for $\Hat{Y}_{n,0}, \Hat{Y}_{0,n}$, or $\Hat{Z}_{0,n}$)  there exists \( m \geq 1\) and a sequence of points 
\begin{equation}\label{eq:p_enum}
  (\omega_{1},x_{1})  :=(\omega^{0},x^{0}) ,(\omega^{1},x^{0}),(\omega^{1},x^{1}),(\omega^{2},x^{1}),(\omega^{2},x^{2}), \dots (\omega^{m},x^{m-1}), (\omega^{m},x^{m}):=(\omega_{2},x_{2})
\end{equation}
belonging to the interior of  $  \Hat{Z}_{n}$ (respectively, in $\Hat{Y}_{n,0}, \Hat{Y}_{0,n}$ ,or $\Hat{Z}_{0,n}$) such that, if \( \omega_{1} \neq \omega_{2}\) then  the sequence  \( \{\omega^{i}\} \) is strictly monotone  in the  interval \( [\omega_{1}, \omega_{2}]\), and if \( x_{1}\neq x_{2}\), then the sequence \( \{x^{i}\} \)  is  strictly monotone  in the interval  \( [x_{1}, x_{2}]\) (in particular \( x^{m-1}\neq x^{m})\).
\end{lemma}

\begin{remark}
Notice that it follows from the indexing of the points in \eqref{eq:p_enum} that consecutive points either have the same horizontal coordinate or the same vertical coordinate, and that these alternate. So joining each two consecutive points by a straight line segment produces a ``staircase''. The monotonicity requirement implies that this staircase is either ascending or descending (and does not go ``up and down'') and also the first edge, between \(  (\omega_{1},x_{1})  :=(\omega^{0},x^{0})  \) and \( (\omega^{1},x^{0})\), is horizontal, and the last edge, between \( (\omega^{m},x^{m-1})\) and \( (\omega^{m},x^{m}):=(\omega_{2},x_{2})\) is vertical. This last requirement is not conceptually so crucial but will be used in some of the applications of the lemma. 
\end{remark}

In the third lemma we consider two points \( (\omega_{1},x), (\omega_{2},x)\) belonging to  the same region \( \Hat{Z}_{n}\) or \( \hat Y_{0,n}\), with the same ``vertical'' coordinate \( x \) and  two different ``horizontal'' coordinates \( \omega_{1}, \omega_{2}\) inside the same element of the partition \( \mathscr P_{k}\), recall \eqref{eq:refinement}. We will show that their ``vertical'' distance after \( k \) iterates  is bounded by a constant factor of their ``horizontal'' distance after \( k \) iterates,

\begin{lemma}\label{lem:expansion_contol_by_horizontal}
There exists \(\Cl{Cveh}> 1\) such that for all $1\leq k\leq n$, and all pairs $(\omega_{1},x), (\omega_{2},x) \in \Hat{Z}_{n}$, and similarly for all pairs \((\omega_{1},x),(\omega_{2},x) \in \hat Y_{0,n}\), with $\omega_{1},\omega_{2} \in \Omega \in \mathscr{P}_{k}$ we have
   \begin{equation*}
       |T^{k}_{\omega_{1}}(x)-T^{k}_{\omega_{2}}(x)| \leq \Cr{Cveh} |F^{k}(\omega_{1})-F^{k}(\omega_{2})|.
   \end{equation*}
\end{lemma}

The three lemmas above will be used in Section \ref{sec:proof_of_pro_expansion} to prove Proposition \ref{pro:expansion}  and in Section \ref{sec:overview_of_the_proof} to prove Proposition \ref{pro:implication_of_distortion}. In the rest of this section we give the proofs of the lemmas.

\subsubsection{Vertical Expansion}

\begin{proof}[Proof of Lemma \ref{lem:vexpansion}]
Both \eqref{eq:exp12} and \eqref{eq:vexpansion1} follows from \eqref{eq:exp11}, indeed,  for all $(\omega,x) \in \mathscr{A}_{n,k,r}$ we have $(F^{n}(\omega), T^{n}_{\omega}(x)) \in \Hat{Y}_{0,k}$ and therefore, applying both the  inequalities in \eqref{eq:exp11} for \(i=j=0\), we get 
\[
\frac{\partial}{\partial x}T^{\tau(\omega,x)}_{\omega}(x)  = \frac{\partial}{\partial x}T^{n+k}_{\omega}(x)  = \frac{\partial}{\partial x}T^{k}_{F^{n}(\omega)}(T^{n}_{\omega}(x)) \frac{\partial}{\partial x}T^{n}_{\omega}(x)\geq \widetilde{\lambda}^{2}.
\]
For \eqref{eq:vexpansion1}, choose \(n>m\) then \(\widehat{F}^{n-m}: \hat{Z}_{n} \to \hat{Y}_{m,0}\) is surjective, in particular for every \((\omega,x) \in \hat{Y}_{m,0}\), we have some \((\hat \omega, \hat x) \in \hat{Z}_{n}\) such that \(\widehat{F}^{n-m}(\hat{\omega},\hat x) =(F^{n-m}(\hat \omega), T^{n-m}_{\hat \omega}(\hat x))=(\omega,x)\). Therefore, applying the first inequality in \eqref{eq:exp11} to the point \((\hat \omega, \hat x)\) and for \(i=n-m\), we get the first inequality in \eqref{eq:vexpansion1}. Similarly, the second inequality in \eqref{eq:vexpansion1} follows from the second inequality in \eqref{eq:exp11}.

It remains to verify \eqref{eq:exp11}. Let \(c^{-}:=\min\{Z_{0,1}(\omega); ~ \omega \in \mathbb S^{1}\}\). Since \(Z_{0,1}\) is continuous, its graph is compact and disjoint from \(\mathbb S^{1} \times \{-1\}\), therefore, \(c^{-}>-1\) which leads to $\widetilde{\lambda}_{1}:=\min\{\frac{\partial}{\partial x}T_{\omega}(c^{-}); ~\omega \in \mathbb{S}^{1}\}>1$. By assumption (A2), we have $\frac{\partial^{2}}{\partial x^{2}}T_{\omega}\geq 0$, therefore, $\frac{\partial}{\partial x}T_{\omega}$ increases in $(-1,0)$. In particular, for all $(\omega,x) \in \hat{Z}_{0,1} \cup \Hat{Z}_{0,0}$, $\frac{\partial}{\partial x}T_{\omega}(x)\geq \widetilde{\lambda}_{1}$. Similarly, let \(c^{+}:=\max\{Y_{1,0}(\omega); ~ \omega \in \mathbb S^{1}\}\) and $\widetilde{\lambda}_{2}:=\min\{\frac{\partial}{\partial x}T_{\omega}(c^{+}); ~ \omega \in \mathbb{S}^{1}\}>1$, then by (A2), for all \((\omega,x)\in \hat Y_{0,0} \cup \hat{Y}_{1,0}\), \(\frac{\partial}{\partial x}T_{\omega}(x)\geq \widetilde{\lambda}_{2}\). Let $\widetilde{\lambda}=\min\{\widetilde{\lambda}_{1},\widetilde{\lambda}_{2}\}$. Now, for the first part of \eqref{eq:exp11}, notice that for each \(0\leq i<n\), the point  \(\widehat{F}^{n-i-1}(F^{i}(\omega),T^{i}_{\omega}(x))= (F^{n-1}(\omega),T_{\omega}^{n-1}(x))\) is in \(\hat{Z}_{0,0}\) \((n=1)\) or \(\hat{Y}_{1,0}\). Therefore, by the above analysis, the final derivative in the fibre before time \(n\) is \(\frac{\partial}{\partial x}T_{F^{n-1}(\omega)}(T^{n-1}_{\omega}(x))\geq \widetilde{\lambda}\). Therefore, by the chain rule and the fact that the derivatives at the rest of all the iterates are greater than or equal to \(1\) (recall Remark \ref{rem:derivative_ge1}) we have the  conclusion. The second part follows by a similar argument.
\end{proof}

\subsubsection{The Staircase Lemma}

In this section we prove Lemma \ref{lem:stairs_between_points}. 
 Before giving the details we make a simple observation concerning the relation between the geometry of the construction and the indexing in \eqref{eq:p_enum}. Suppose that we have a path from \( (\omega_{1},x_{1})\) to \( (\omega_{2},x_{2})\) made up of alternating horizontal edges going to the right and vertical edges going down, whose vertices all lie in the interior of \( \Hat{Z}_{n}\), and whose first edge is horizontal and last edge is vertical, as required by the Lemma. Then the vertices of the path are exactly a sequence of the form \eqref{eq:p_enum} in which \( \{\omega^{i}\}\) is strictly increasing and \( \{x^{i}\}\) is strictly decreasing. It is therefore sufficient to construct such a path, the indexing follows automatically.

\begin{figure}[h]
    \centering
    \includegraphics[width=0.55\linewidth]{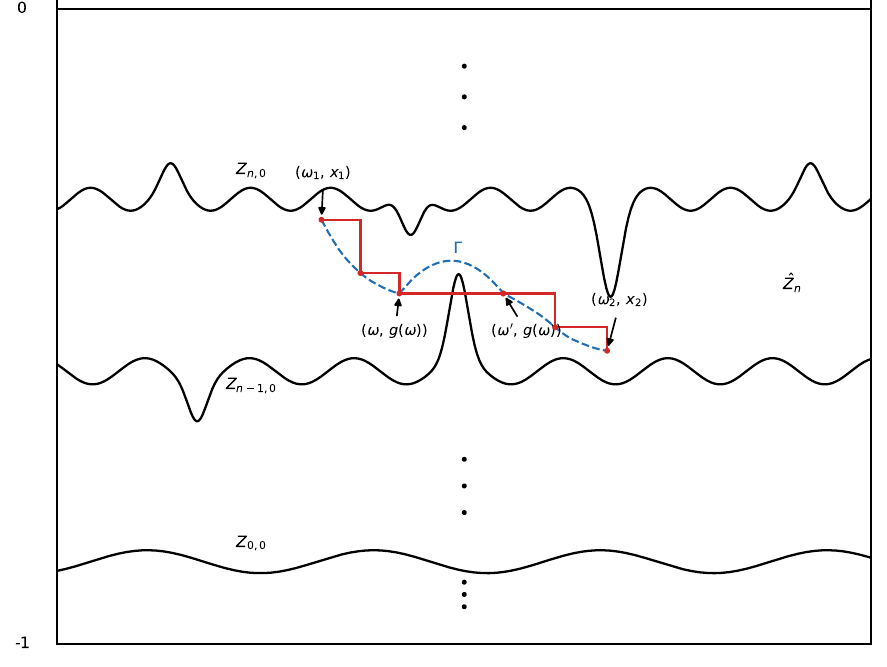}
    \caption{Connecting  two points through a Staircase}
    \label{fig:stairs_path}
\end{figure}

\begin{figure}[h]
    \centering
    \includegraphics[width=0.7\linewidth]{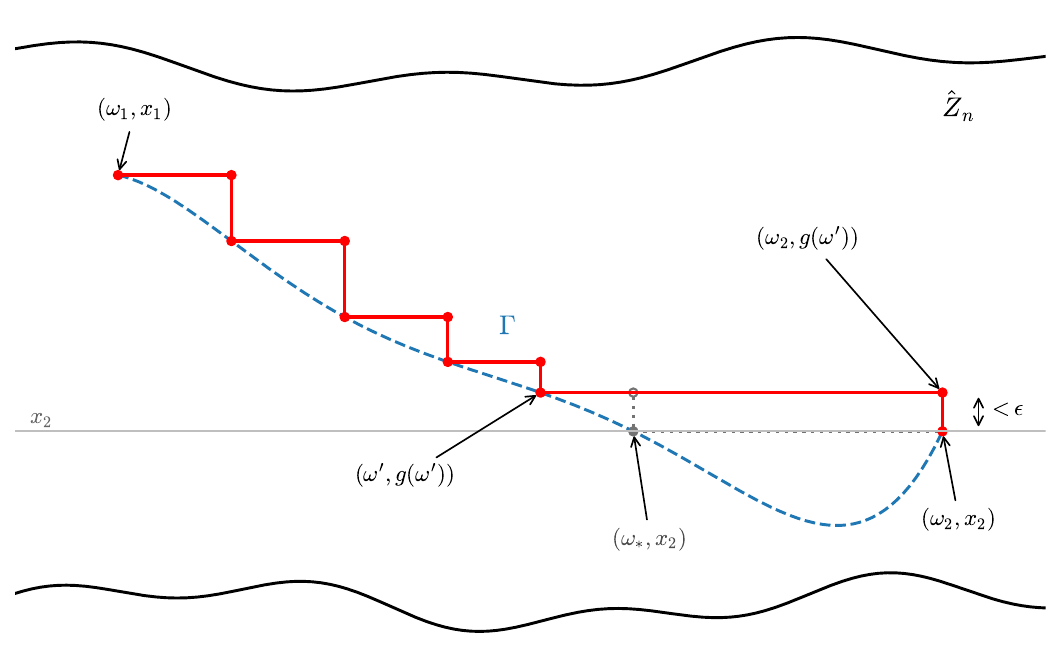}
    \caption{Cutting the corner when \( \bar\Gamma \) ends with a horizontal segment. }
    \label{fig:stairs_corner}
\end{figure}

\begin{proof}[Proof of Lemma \ref{lem:stairs_between_points}]
We give the proof for  $\Hat{Z}_{n}$, the same argument works for the other strips. If \(\omega_{1} = \omega_{2}\)  and/or \(x_{1}=x_{2} \) then the corresponding sequence is trivial. So we suppose that  
\[
\omega_{2} > \omega_{1}
\quad \text{ and } \quad  
x_{2} <  x_{1}.
\] 
The other cases all follow by exactly the same argument. 
Observe first of all that since the strip $\Hat{Z}_{n}$ is a connected region whose boundaries are two smooth curves we can find a smooth function   
\[
 g:[\omega_{1},\omega_{2}] \to I,  \text{ with } g(\omega_{1})=x_{1} \text{ and } g(\omega_{2})=x_{2},
 \] 
 with  at most a finite number of critical points,  and  whose graph $\Gamma$ lies in the interior of $\Hat{Z}_{n}$. In particular, there exists \( \epsilon > 0\) such that \( \Gamma \) has distance \( \geq \epsilon\) from the boundary of  $\Hat{Z}_{n}$ in the metric~\eqref{eq:metricd}.

We consider two cases, depending on whether \( g \) has critical points.

\vspace{3mm}
\noindent \emph{Case (1): \(g\) is strictly monotone decreasing (no critical points).}
Then we can find a strictly monotone  sequence  $\omega_{1}=\omega^{0} < \omega^{1} <  \dots <\omega^{k-1}< \omega^{k}= \omega_{2}$ such that $|\omega^{i}-\omega^{i+1}| < \epsilon$,  the corresponding sequence  $x^{i}= g(\omega^{i})$ is also strictly monotone, and therefore we get the required sequence \eqref{eq:p_enum}.  Indeed the points \( (\omega^{i},x^{i}) \) lie on \( \Gamma \) and the points \( (\omega^{i+1},x^{i}) \) are at distance \( |\omega^{i+1}-\omega^{i}|<\epsilon \) from \( (\omega^{i},x^{i})\in\Gamma \), and therefore all the points of \eqref{eq:p_enum} lie in the interior of \( \Hat{Z}_{n}\).

\medskip
\noindent \emph{Case (2): \(g\) is not strictly monotone decreasing. }
Then \( g \) has finitely many monotone branches, and applying the argument above to each of them would give a sequence of staircases going alternately ``down'' and ``up'' rather than a single monotone staircase. The idea is instead to follow \( \Gamma \) only while it is descending, and  whenever \( \Gamma \) turns back up, to ``wait'' at the current height, taking a single horizontal step over the whole excursion of \( \Gamma \) above this height and rejoining \( \Gamma \) at the point where it comes back down to the same height, see Figure~\ref{fig:stairs_path}. Notice that such a horizontal step may well leave the strip, but this is not a problem since we only require the \emph{vertices} of the staircase,  and not its edges, to lie in the interior of \( \Hat{Z}_{n}\), and both endpoints of the step lie on \( \Gamma \).

To do this, for \( \omega \in [\omega_{1}, \omega_{2}] \) we let 
\[
\bar g(\omega) := \max\big\{ \min\{ g(s) : \omega_{1}\leq s \leq \omega \}, \ x_{2} \big\}.
\]
Then \( \bar g \) is continuous and non-increasing, with \( \bar g(\omega) \geq x_{2}\), \( \bar g(\omega_{1}) = x_{1}\) and \( \bar g(\omega_{2}) = x_{2}\), and since \( g \) has finitely many critical points, its graph \( \bar\Gamma \) consists of finitely many arcs of \( \Gamma \) on which \( g \) is strictly decreasing, alternating with finitely many horizontal segments with endpoints on \( \Gamma \). We now build our path along \( \bar\Gamma \): on each descending arc we take the staircase given by \emph{Case} (1), and each horizontal segment we take as a single horizontal edge  (when a horizontal segment is followed by a descending arc, the segment and the first, horizontal, edge of the staircase on that arc are two consecutive collinear edges and we consider them as a single horizontal edge, so that horizontal and vertical edges alternate). Every vertex of this path either lies on \( \Gamma \) or, as in \emph{Case} (1), is within distance \( \epsilon \) of a point of \( \Gamma \), and therefore lies in the interior of \( \Hat{Z}_{n}\); the path moves only to the right and downwards, since \( \bar g \) is non-increasing. 

To conclude the proof, we just need to verify that the first edge is horizontal and the last edge is vertical. It is easy to see that the  first edge is horizontal, since both a horizontal segment and a staircase from \emph{Case} (1) begin with a horizontal edge. On the other hand, the last edge obtained from the construction above is not necessarily vertical. 
If \( \bar\Gamma \) ends with a descending arc then the last edge of the path is vertical, since the staircases of \emph{Case} (1) end with a vertical edge, and we are done. 
It just remains to consider the case in which \( \bar\Gamma \) ends with a horizontal segment \( [\omega_{*},\omega_{2}]\times\{x_{2}\} \), which happens if and only if \( g(\omega)=x_{2}\) for some \( \omega\in(\omega_{1},\omega_{2}) \), \( \omega_{*} \) being the smallest such \( \omega \), so that the path ends with a horizontal edge. In this case we modify the end of the path by ``cutting the corner'' at \( (\omega_{*},x_{2})\), see Figure \ref{fig:stairs_corner}. Notice first of all that this segment is preceded by a descending arc of \( \bar\Gamma \) (two horizontal segments are never adjacent, since \( \bar g \) is continuous and distinct segments lie at distinct heights, and \( \bar\Gamma \) is not a single horizontal segment since \( x_{1}>x_{2}\)), i.e. \( g \) is strictly decreasing on some interval \( [\omega_{a},\omega_{*}]\) with \( g(\omega_{*})=x_{2}\), and the staircase on this arc given by \emph{Case}~(1) has vertices \( (\omega^{i},g(\omega^{i}))\in\Gamma \) for \( \omega_{a}=\omega^{0}<\omega^{1}<\dots<\omega^{k}=\omega_{*} \). We make no assumption on \( g \) to the right of \( \omega_{*} \), where it may well go below \( x_{2}\), nor on whether \( \omega_{*} \) is a critical point of \( g \); we only use the fact that \( g \) is continuous and strictly decreasing on \( [\omega_{a},\omega_{*}] \), so that \( g(\omega)-x_{2} \) is positive and tends to \( 0 \) as \( \omega\to\omega_{*}\). Since in \emph{Case} (1) the points \( \omega^{i} \) can be chosen as close to each other as we like, we can choose \( \omega':=\omega^{k-1} \) so close to \( \omega_{*} \) that \( 0< g(\omega')-x_{2}<\epsilon\). We then replace the final part of the path after \( (\omega',g(\omega')\), which consists of a horizontal edge to \( (\omega_{*},g(\omega')) \), a vertical edge to \( (\omega_{*},x_{2})\) and the horizontal edge to \( (\omega_{2},x_{2}) \), by a single horizontal edge to \( (\omega_{2},g(\omega'))\) followed by a vertical edge down to \( (\omega_{2},x_{2})\). All the other vertices of the path are unchanged, and the only new vertex \( (\omega_{2},g(\omega')) \) lies directly above \( (\omega_{2},x_{2})\in\Gamma \) at distance \( g(\omega')-x_{2}<\epsilon \), and therefore lies in the interior of \( \Hat{Z}_{n}\). The path now ends with a vertical edge, which completes the proof.
\end{proof}

\subsubsection{Control of the vertical expansion by the horizontal}\label{sec:controle_of_VE_by_HE} 
In this section we prove Lemma \ref{lem:expansion_contol_by_horizontal}. We first state and prove two auxiliary sublemmas. 

\begin{sublemma}\label{lem:mean_value}
There exists \(\Cl{Cveh1} >0\) such that for all $x \in [-1,1]$ and  $\omega_{1},\omega_{2} \in \mathbb{S}^{1}$,
\begin{equation}\label{eq:mean_value11}
    |T_{\omega_{1}}(x)-T_{\omega_{2}}(x)| \leq \Cr{Cveh1} |\omega_{1}-\omega_{2}|.
    \end{equation}
\end{sublemma}

\begin{proof}
Notice that, for \(x=-1,0,1\) the left hand side vanishes and the result is trivial. 
So we just consider  \(x  \in (-1,0) \cup (0,1) \).  
To simplify the notation and emphasize that the point \( x \) is fixed and the dependence is just on \( \omega\), we let  $g(\omega) :=  T_{\omega}(x)$. 
Now for \( \omega_{1},\omega_{2} \in \mathbb{S}^{1}\),  let $h_{\omega_{1},\omega_{2}}:[0,1] \rightarrow \mathbb{S}^{1}$  be the parametrization of the smooth arc joining $\omega_{1}$ and $\omega_{2}$ on the circle.
Then we have, $g\circ h_{\omega_{1},\omega_{2}} :[0,1] \rightarrow I$ with $g\circ h_{\omega_{1},\omega_{2}}(0) = g(\omega_{1})$ and $g\circ h_{\omega_{1},\omega_{2}}(1) = g(\omega_{2})$, and  there exists a $\zeta \in [0,1]$ such that
\begin{equation}\label{eq:ind11}
    |g(\omega_{1})-g(\omega_{2})| \leq  |g^{'}(h_{\omega_{1},\omega_{2}}(\zeta))| \int_{0}^{1}|h_{\omega_{1},\omega_{2}}^{'}(t)| dt  \leq \max_{\omega}\left|\frac{\partial}{\partial \omega}T_{\omega}(x)\right|  |\omega_{1}-\omega_{2}|.
\end{equation}
It is therefore enough to show that \(|\frac{\partial}{\partial \omega}T_{\omega}(x)|\) is uniformly bounded. If $x \in {U_{+1}}$, then 
\begin{equation*}
    \frac{\partial}{\partial \omega} T_{\omega}(x) =  -b_{2}(\omega) \big(1-x\big)^{\ell_{2}(\omega)+1} \ln(1-x) \ell_{2}^{'}(\omega) -b_{2}^{'}(\omega)(1-x)^{\ell_{2}(\omega)+1} - \frac{\partial}{\partial \omega} \eta_{\omega}^{(2)}(x).
\end{equation*}
By assumption \eqref{eq:higherorder} and the fact that \(|\big(1-x\big)^{\ell_{2}(\omega)+1} \ln(1-x)|\leq \left(e(\ell_{2}^{-}+1)\right)^{-1}\), we have
\[
\Big|\frac{\partial}{\partial \omega} T_{\omega}(x)\Big|\leq \frac{(\delta+b_{2}^{+})\max\{|\ell_{2}^{'}(\omega)|\}}{e(\ell_{2}^{-}+1)}+\max\{|b_{2}^{'}(\omega)|\}\lesssim 1.
\]
A  similar argument gives a uniform upper bound for \(x \in U_{-1}\). Similarly for $x \in U_{-0}$, we have
\begin{equation*}
    \Big|\frac{\partial}{\partial \omega} T_{\omega}(x)\Big| \leq \frac{(\delta+a_{1}^{+}) \max\{|k_{1}^{'}(\omega)|\}}{ek_{1}^{-}}+\max\{|a_{1}^{'}(\omega)|\}\lesssim 1.
\end{equation*}
A similar bound for \(x \in U_{+0}\). If \(x \notin U_{\pm1}\cup U_{\pm0}\), then the map $(\omega,x) \to T_{\omega}(x)$ is $C^{2}$ on the closure of $\mathbb{S}^{1} \times (I\setminus (U_{\pm1}\cup U_{\pm0}))$ and therefore, it has uniformly bounded derivatives there. This completes the proof of the Sublemma.
\end{proof}
In the next sublemma we are starting with two points at the same vertical height, but after the first iteration we might lose this advantage, and the points might not remain at the same vertical height. We resolve this issue by using Lemma \ref{lem:stairs_between_points}.

\begin{sublemma}\label{lem:inductive_expansion_control}
There exist \( C>0 \), \( N\geq 1 \),  and a sequence \(\{\nu_i\}\),  such that \( \nu_{i}\leq C \) for all \( i \), and \( \nu_{i}\leq 1+b^{+}+\delta \) for  \( i\geq N \) (in particular at most \( N \) of the \( \nu_{i} \) exceed \( 1+b^{+}+\delta \)), such that for all \(n\geq 1\) and for all \(k\leq n\)  and  all $(\omega_{1},x), (\omega_{2},x) \in \Hat{Z}_{n}$, similarly for all $(\omega_{1},x), (\omega_{2},x)\in \Hat{Y}_{0,n}$, with $\omega_{1},\omega_{2} \in \Omega \in \mathscr{P}_{k}$ we have,
    \begin{equation}\label{eq:hexpansion_controle}
        \big|T^{k}_{\omega_{1}}(x)-T^{k}_{\omega_{2}}(x)\big| \leq \sum_{i=0}^{k-1}\nu_{0}\nu_{1}\cdots \nu_{i} |F^{k-i-1}(\omega_{1})-F^{k-i-1}(\omega_{2})|
    \end{equation}

\end{sublemma}

\begin{proof} 
We give the proof in case $(\omega_{1},x),(\omega_{2},x)\in \Hat{Z}_{n}$ and the other will follow from a similar reasoning. First of all notice that for \(k=1\) the statement follows from Lemma \ref{lem:mean_value}, therefore we consider the case \(k>1\). Now notice that $(F^{k-1}(\omega_{1}),T^{k-1}_{\omega_{1}}(x)),(F^{k-1}(\omega_{2}),T^{k-1}_{\omega_{2}}(x)) \in \Hat{Y}_{n-k+1,0}$. Let \(T^{k-1}_{\omega_{1}}(x)\geq T^{k-1}_{\omega_{2}}(x)\) and consider the sequence of points between $(F^{k-1}(\omega_{1}),T^{k-1}_{\omega_{1}}(x))$ and $(F^{k-1}(\omega_{2}),T^{k-1}_{\omega_{2}}(x))$ generated by Lemma \ref{lem:stairs_between_points}, that is for some \(m\geq 1\) we have
\[(\omega^{0},x^{0})= (F^{k-1}(\omega_{1}),T^{k-1}_{\omega_{1}}(x)),(\omega^{1},x^{0}),(\omega^{1},x^{1}), \dots, \newline (\omega^{m},x^{m-1}), (\omega^{m},x^{m})=(F^{k-1}(\omega_{2}),T^{k-1}_{\omega_{2}}(x)).\] 
Now take the left side of \eqref{eq:hexpansion_controle}, add and subtract $T_{\omega^{i}}(x^{i-j})$ for all $i=1, \dots m$ and $j=0,1$, 
\begin{equation}\label{eq:control101}
\begin{split}
    \big|T^{k}_{\omega_{1}}(x)-T^{k}_{\omega_{2}}(x)\big| =|T_{\omega^{0}}(x^{0})-T_{\omega^{m}}(x^{m})|= \big| \sum_{i=0}^{m-1}T_{\omega^{i}}(x^{i})-T_{\omega^{i+1}}(x^{i})+\sum_{j=1}^{m}T_{\omega^{j}}(x^{j-1})-T_{\omega^{j}}(x^{j})\big|.
\end{split}
\end{equation}
Let $I_{i}= |T_{\omega^{i}}(x^{i})-T_{\omega^{i+1}}(x^{i})|$ for all $i=0, \dots m-1$, and $J_{j}= |T_{\omega^{j}}(x^{j-1})-T_{\omega^{j}}(x^{j})|$ for all $j =1, \dots m$, therefore \eqref{eq:control101} leads to
\begin{equation}\label{eq:hexpansion_controle1}
    \big|T^{k}_{\omega_{1}}(x)-T^{k}_{\omega_{2}}(x)\big|\leq \sum_{i=0}^{m-1}I_{i}+\sum_{j=1}^{m}J_{j}.
\end{equation}
By Lemma \ref{lem:mean_value} we have \(I_{i}\leq \Cr{Cveh1}|\omega^{i}-\omega^{i+1}|\). Notice that by assumption (A2), the derivative of  \(T_{\omega^{j}}\) is decreasing in \([x^{j},x^{j-1}]\), hence the maximum is achieved at \(x^{j}\), also notice that \((\omega^{j},x^{j}) \in \hat Y_{n-k+1,0}\). Therefore, applying the Mean Value Theorem to each \(J_{j}\) and letting \(\nu_{1}:=\max\{|\frac{\partial}{\partial x}T_{\omega}(x)|; \quad (\omega,x) \in \hat Y_{n-k+1,0}\}\), inequality \eqref{eq:hexpansion_controle1} leads to 
\begin{equation}\label{eq:controle1}
     \big|T^{k}_{\omega_{1}}(x)-T^{k}_{\omega_{2}}(x)\big|\leq \Cr{Cveh1} \sum_{i=0}^{m-1}|\omega^{i}-\omega^{i+1}|+ \sum_{j=1}^{m}\nu_{1}|x^{j-1}-x^{j}|= \Cr{Cveh1}|\omega^{0}-\omega^{m}|+\nu_{1}|x^{0}-x^{m}|,
\end{equation}
where the last equality follows because \(\{\omega^{0},\omega^{1}, \dots, \omega^{m}\}\) and \(\{x^{0},x^{1}, \dots,x^{m}\}\) are in order, as given by Lemma \ref{lem:stairs_between_points}. Putting back the values of \(\omega^{0},\omega^{m},x^{0}\) and \(x^{m}\) in \eqref{eq:controle1} we have
\begin{equation}\label{eq:estemats1231}
    \big|T^{k}_{\omega_{1}}(x)-T^{k}_{\omega_{2}}(x)\big|\leq  \Cr{Cveh1}|F^{k-1}(\omega_{1})-F^{k-1}(\omega_{2})|+\nu_{1}|T^{k-1}_{\omega_{1}}(x)-T^{k-1}_{\omega_{2}}(x)|.
\end{equation}
The case \(T^{k-1}_{\omega_{1}}(x)\leq T^{k-1}_{\omega_{2}}(x)\) follows similarly. We also apply the above argument to the points \((F^{k-2}(\omega_{1}),T^{k-2}_{\omega_{1}}(x))\), \((F^{k-2}(\omega_{2}),T^{k-2}_{\omega_{2}}(x))\), we get the above estimates for $|T^{k-1}_{\omega_{1}}(x)-T^{k-1}_{\omega_{2}}(x)|$ with \(\nu_{2}\) as the maximum of the derivative over \(\hat{Y}_{n-k+2,0}\), we then put these estimates into \eqref{eq:estemats1231}. Repeating this process $k-2$ time, we find $\nu_{0}=\Cr{Cveh1}, \nu_{1}, \nu_{2}, \nu_{3}, \dots, \nu_{k-2}$ such that,
\begin{equation}\label{eq:control111}
    \big|T^{k}_{\omega_{1}}(x)-T^{k}_{\omega_{2}}(x)\big|\leq \sum_{i=0}^{k-2}\nu_{0}\nu_{1} \dots \nu_{i}|F^{k-1-i}(\omega_{1})-F^{k-1-i}(\omega_{2})|+\nu_{0} \dots \nu_{k-2}|T_{\omega_{1}}(x)-T_{\omega_{2}}(x)|.
\end{equation}
Using Lemma \ref{lem:mean_value} for the last term on the right of \eqref{eq:control111} then for $\nu_{k-1}=\Cr{Cveh1},$ we have
\begin{equation}\label{eq:hexpansion_control3}
    \big|T^{k}_{\omega_{1}}(x)-T^{k}_{\omega_{2}}(x)\big|\leq \sum_{i=0}^{k-1}\nu_{0}\nu_{1} \dots \nu_{i}|F^{k-1-i}(\omega_{1})-F^{k-1-i}(\omega_{2})|.
\end{equation}
For uniform boundedness, notice that by (A2), for all \(i=1, \dots,k-2\) we have \(\nu_{i}\leq \max\{|\frac{\partial}{\partial x}T_{\omega}(x)|; ~(\omega,x) \in \hat Y_{1,0}\}:=C\). It remains to verify that as \(n\) goes to infinity, then for most \(i\) we have, \(\nu_{i}<1+b^{+}+\delta\). Let \(T_{\ell^{-}_{2}}(x):=x-(b^{+}_{2}+\delta)(1-x)^{\ell^{-}_{2}+1} \), then for \(n\) sufficiently large \(\nu_{i}\leq \max\{|T_{\ell^{-}_{2}}^{'}(x)|;~(\omega,x) \in \hat Y_{n-k+i,0}\}\) for all \(i=1, \dots ,k-2\). Therefore, if \(\ell^{-}_{2}=0\) and/or \(\ell^{+}_{2} =0\), then \(\nu_{i} \leq T^{'}_{\ell^{-}_{2}}(x)=1+b_{2}^{+}+\delta\), and if $\ell^{+}_{2} \neq 0$ and \(\ell^{-}_{2}\neq 0\), then using the estimates from \eqref{eq:Y_n} we have
\begin{equation}\label{eq:most_i-old}
 \nu_{i} \leq 
     1+(b^{+}_{2}+\delta)(\ell^{-}_{2}+1)C^{'}(n-k+i)^{\frac{-\ell^{-}_{2}}{\ell^{+}_{2}}},  
\end{equation}
where \( C' \) is the constant in \eqref{eq:Y_n}. As \(n\) goes to infinity, the right side converges to \(1\), in particular, there exists \(N\) large enough such that \(\nu_{i} \leq 1+b^{+}+\delta\) for all \(i \geq N\). 
\end{proof}
We are now ready to complete the proof of Lemma \ref{lem:expansion_contol_by_horizontal}. We first make an elementary remark, which will also be used in other arguments below. 

\begin{remark}
We remark that the uniformly expanding map \(F\) has the bounded distortion property, in particular for all \(i\), \(\omega_{1}, \omega_{2} \in \Omega \in \mathscr{P}_{i}\), and \(j\leq i\)
\begin{equation}\label{eq:Fbd}
 \ln\frac{|\left(F^{j}\right)^{'}(\omega_{1})|}{|\left(F^{j}\right)^{'}(\omega_{2})|}\lesssim |F^{j}(\omega_{1})-F^{j}(\omega_{2})|\leq 1.
\end{equation}
\end{remark}

\begin{proof}[Proof of Lemma \ref{lem:expansion_contol_by_horizontal}]

First of all notice that by the Mean Value Theorem, \eqref{eq:Fbd} and assumption (A1),
\begin{equation}\label{eq:expansion_control4}
    \frac{|F^{k-i-1}(\omega_{1})-F^{k-i-1}(\omega_{2})|}{|F^{k}(\omega_{1})-F^{k}(\omega_{2})|}\lesssim \frac{|\omega_{1}-\omega_{2}|}{|F^{i+1}(\omega_{1})-F^{i+1}(\omega_{2})|} \leq \frac{1}{\theta^{i+1}}.
\end{equation}
Substituting \eqref{eq:expansion_control4} into  the right side of the inequality \eqref{eq:hexpansion_controle} we have
\begin{equation}\label{eq:seriesf}
  \sum_{i=0}^{k-1} \nu_{0}\nu_{1}\cdots\nu_{i}|F^{k-i-1}(\omega_{1})-F^{k-i-1}(\omega_{2})| \lesssim |F^{k}(\omega_{1})-F^{k}(\omega_{2})| \sum_{i=0}^{k-1} \frac{\nu_{0}\nu_{1}\cdots \nu_{i}}{\theta^{i+1}}.
\end{equation}
By Sublemma \ref{lem:inductive_expansion_control} we have \( \nu_{i}\leq C \) for all \( i \) and, by assumption (A5), \( \nu_{i}\leq 1+b^{+}+\delta<\theta \) for all but at most \( N \) values of \( i \), with \( C \) and \( N \) independent of \( n,k \) and of the points. Therefore \( \nu_{0}\cdots\nu_{i}\leq C^{N}(1+b^{+}+\delta)^{i+1} \), and so the sum on the right of \eqref{eq:seriesf} is bounded by \( C^{N}\sum_{i\geq 0}((1+b^{+}+\delta)/\theta)^{i+1}<\infty \), uniformly in \( n \) and \( k \), hence the required result follows. Notice that this is the only point in the proof where the second condition in \eqref{eq:expans_cond} is used, and only in the case \( \ell^{-}=0 \), since otherwise \( \nu_{i}\approx 1 \) for most \( i \) by \eqref{eq:most_i-old}, recall Remark \ref{rem:expans_cond}.
\end{proof}

\subsection{Proof of Proposition \ref{pro:expansion}}\label{sec:proof_of_pro_expansion}

We now use Lemmas \ref{lem:vexpansion}-\ref{lem:stairs_between_points} to prove Lemmas \ref{lem:lower_bound_VE1}  and  \ref{lem:inductive_step2} below, which we will then use to complete the proof of  Proposition \ref{pro:expansion}. 
Recall $\widetilde{\lambda}$ from Lemma \ref{lem:vexpansion} and choose 
\begin{equation}\label{eq:tilde_lambda}
    1<\lambda<\min\{\widetilde{\lambda},\theta\}.
\end{equation} 
We begin by proving this crucial expansion estimate. 

\begin{lemma}\label{lem:lower_bound_VE1}
    For all \((\omega_{1},x_{1}),(\omega_{2},x_{2}) \in \hat Z_{n}\), similarly for all  \((\omega_{1},x_{1}),(\omega_{2},x_{2}) \in \hat{Y}_{0,n}\), with \(\omega_{1},\omega_{2} \in \Omega \in \mathscr{P}_{n}\), and for all \(0\leq p< n\) we have
\begin{equation}\label{eq:distortion1}
    |T^{n}_{\omega_{1}}(x_{1})-T^{n}_{\omega_{2}}(x_{2})|\geq \lambda |T_{\omega_{1}}^{p}(x_{1})-T_{\omega_{2}}^{p}(x_{2})|-\Cr{Cveh}|F^{n}(\omega_{1})-F^{n}(\omega_{2})|.
\end{equation}
In particular, we have \eqref{eq:distortion4}.
\end{lemma}
\begin{proof}
It suffices to consider $(\omega_{1},x_{1}),(\omega_{2},x_{2})\in\hat{Z}_{n}$, as the case $(\omega_{1},x_{1}),(\omega_{2},x_{2})\in\hat{Y}_{0,n}$ follows by an analogous argument.  Notice that in this case  \((F^{p}(\omega_{1}),T_{\omega_{1}}^{p}(x_{1})), (F^{p}(\omega_{2}),T_{\omega_{2}}^{p}(x_{2})) \in \Hat{Y}_{n-p,0}\) or \(\hat{Z}_{n}\)(when \(p=0\)). Let $T^{p}_{\omega_{1}}(x_{1})\geq T_{\omega_{2}}^{p}(x_{2})$ and consider the sequence of points \eqref{eq:p_enum} between \((F^{p}(\omega_{1}),T_{\omega_{1}}^{p}(x_{1}))\),and \((F^{p}(\omega_{2}),T_{\omega_{2}}^{p}(x_{2}))\) given  by Lemma \ref{lem:stairs_between_points}.  That is  for some \(m\geq 1\)
\[(\omega^{0},x^{0})= (F^{p}(\omega_{1}),T^{p}_{\omega_{1}}(x_{1})),(\omega^{1},x^{0}),(\omega^{1},x^{1}), \dots, \newline (\omega^{m},x^{m-1}), (\omega^{m},x^{m})=(F^{p}(\omega_{2}),T^{p}_{\omega_{2}}(x_{2})).\]
All points lie in a same strip \(\Hat{Y}_{n-p,0}\) or \(\Hat{Z}_{n}\). Adding and subtracting $T^{n-p}_{\omega^{i}}(x^{i-j})$ for all $i =1, \dots , m$ and $j=0,1$  to  \(|T^{n-p}_{\omega^{0}}(x^{0})-T^{n-p}_{\omega^{m}}(x^{m})|\),  we have
\begin{equation}\label{eq:distortion811}
\begin{split}
    |T^{n-p}_{\omega^{0}}(x^{0})-T^{n-p}_{\omega^{m}}(x^{m})|&= \big|\sum_{j=1}^{m}\left(T_{\omega^{j}}^{n-p}(x^{j-1})-T_{\omega^{j}}^{n-p}(x^{j})\right)+ \sum_{i=0}^{m-1}\left(T_{\omega^{i}}^{n-p}(x^{i})-T_{\omega^{i+1}}^{n-p}(x^{i})\right)\big|\\ 
    &\geq \big|\sum_{j=1}^{m}\left(T_{\omega^{j}}^{n-p}(x^{j-1})-T_{\omega^{j}}^{n-p}(x^{j})\right)\big|-\big| \sum_{i=0}^{m-1}\left(T_{\omega^{i}}^{n-p}(x^{i})-T_{\omega^{i+1}}^{n-p}(x^{i})\right)\big|.
    \end{split}
\end{equation}
For $i=0, \dots , m-1$,  we let  $I_{i}:= |T^{n-p}_{\omega^{i}}(x^{i})-T^{n-p}_{\omega^{i+1}}(x^{i})|$,  then by  Lemma \ref{lem:expansion_contol_by_horizontal}, we have \(I_{i} \leq \Cr{Cveh}|F^{n-p}(\omega^{i})-F^{n-p}(\omega^{i+1})|\). For all  $j=1, \dots, m$, we let $J_{j}:=|T^{n-p}_{\omega^{j}}(x^{j-1})-T^{n-p}_{\omega^{j}}(x^{j})|=T^{n-p}_{\omega^{j}}(x^{j-1})-T^{n-p}_{\omega^{j}}(x^{j})$, where the equality follows because $x^{j-1}\geq x^{j}$, and $T^{n-p}_{\omega^{j}}$ is monotone. By \eqref{eq:vexpansion1}  we have \(J_{j}\geq \lambda |x^{j-1}-x^{j}|\). Therefore, \eqref{eq:distortion811} implies  
\begin{equation}\label{eq:bound_by_IiJj}
  |T^{n-p}_{\omega^{0}}(x^{0})-T^{n-p}_{\omega^{m}}(x^{m})| \geq \sum_{j=1}^{m}J_{j}- \sum_{i=0}^{m-1}I_{i}\geq \lambda\sum_{j=1}^{m}|x^{j-1}-x^{j}|-\Cr{Cveh}\sum_{i=0}^{m-1}|F^{n-p}(\omega^{i})-F^{n-p}(\omega^{i+1})|.
\end{equation}
Since \(x^{i}\)'s and \(\omega^{i}\)'s are in order and \(F^{n}\) is a homeomorphism on \(\Omega\), therefore \(\sum_{j=1}^{m}|x^{j-1}-x^{j}|=|x^{0}-x^{m}|\) and \(\sum_{i=0}^{m-1}|F^{n-p}(\omega^{i})-F^{n-p}(\omega^{i+1})|=|F^{n}(\omega^{0})-F^{n}(\omega^{m})|\). Putting this in \eqref{eq:bound_by_IiJj} and  replacing \(\omega^{0},x^{0},\omega^{m},x^{m}\) by their actual values we have \eqref{eq:distortion1}.

Now, to prove  \eqref{eq:distortion4}, since  \((F^{n}(\omega_{1}),T^{n}_{\omega_{1}}(x_{1})), (F^{n}(\omega_{2}),T^{n}_{\omega_{2}}(x_{2})) \in \hat{Y}_{0,k}\), with $F^{n}(\omega_{1}),F^{n}(\omega_{2}) \in \Omega \in \mathscr{P}_{k}$. Therefore, applying \eqref{eq:distortion1} for the points in \(\hat{Y}_{0,k}\), we get that for all \(0\leq s <k\),
\begin{equation}\label{eq:distortion21212}
\begin{split}
    |T^{k}_{F^{n}(\omega_{1})}(T^{n}_{\omega_{1}}(x_{1}))-T^{k}_{F^{n}(\omega_{2})}(T^{n}_{\omega_{2}}(x_{2}))| & \geq \lambda |T^{s}_{F^{n}(\omega_{1})}(T^{n}_{\omega_{1}}(x_{1}))-T^{s}_{F^{n}(\omega_{2})}(T^{n}_{\omega_{2}}(x_{2}))|\\
    &-\Cr{Cveh}|F^{k}(F^{n}(\omega_{1}))-F^{k}(F^{n}(\omega_{2}))|.
\end{split}
\end{equation}
In particular \eqref{eq:distortion4} follows for all \(n\leq p < \tau\). If \(p<n\), then using \eqref{eq:distortion21212} for \(s=0\) we have 
\begin{equation}\label{eq:distortion2}
|T^{\tau}_{\omega_{1}}(x_{1})-T^{\tau}_{\omega_{2}}(x_{2})| \geq \lambda|T^{n}_{\omega_{1}}(x_{1})-T^{n}_{\omega_{2}}(x_{2})|-\Cr{Cveh}|F^{\tau}(\omega_{1})-F^{\tau}(\omega_{2})|
\end{equation}
Use \eqref{eq:distortion1} in \eqref{eq:distortion2} we have
\begin{equation*}
    |T^{\tau}_{\omega_{1}}(x_{1})-T^{\tau}_{\omega_{2}}(x_{2})|  \geq \lambda^{2} |T_{\omega_{1}}^{p}(x_{1})-T_{\omega_{2}}^{p}(x_{2})|-\Cr{Cveh}(1+\frac{\lambda}{\theta^{k}})|F^{\tau}(\omega_{1})-F^{\tau}(\omega_{2})|.
\end{equation*}
This completes the proof of \eqref{eq:distortion4}. Notice that an identical argument works in case $x_{1} \leq x_{2}$, just by taking the sequence of points generated by Lemma \ref{lem:stairs_between_points} starting at $(\omega_{2},x_{2})$ and ending at $(\omega_{1},x_{1})$.
\end{proof}

Given two points \( (\omega_{1},x_{1}),(\omega_{2},x_{2}) \in \mathscr{A}_{n,k,r}\) with separation time \( s:= s((\omega_{1},x_{1}),(\omega_{2},x_{2}))\), for \( 0\leq i< s\) the points \( T^{i}_{\mathbf H}(\omega_{1},x_{1})\) and \( T^{i}_{\mathbf H}(\omega_{2},x_{2})\) lie in the same element of \( \mathscr{R}\) and therefore have the same return time, which we denote by
\[
\tau_{i}:=\tau(T^{i}_{\mathbf H}(\omega_{1},x_{1}))=\tau(T^{i}_{\mathbf H}(\omega_{2},x_{2})),
\]
and for \( 1\leq m\leq s\) we let
\[
t_{m}:= \tau_{0}+\tau_{1}+\dots+\tau_{m-1},
\]
so that \( T^{m}_{\mathbf H}(\omega_{j},x_{j})=\widehat F^{t_{m}}(\omega_{j},x_{j})\) for \( j=1,2\). Notice that \( t_{1}=\tau_{0}=n+k\) and \( t_{m+1}=t_{m}+\tau_{m}\). We also let
\[
\Cl{el}:= 2\Cr{Cveh}\sum_{i=0}^{\infty}\Big(\frac{\lambda}{\theta}\Big)^{i}
= \frac{2\Cr{Cveh}\theta}{\theta-\lambda},
\]
which is finite since \( \lambda<\theta\) by \eqref{eq:tilde_lambda}.

\begin{lemma}\label{lem:inductive_step2}
For all \( (\omega_{1},x_{1}),(\omega_{2},x_{2}) \in \mathscr{A}_{n,k,r}\) and all \( 1\leq m\leq s((\omega_{1},x_{1}),(\omega_{2},x_{2}))\) we have
\begin{equation}\label{eq:distortion13}
|T^{t_{m}}_{\omega_{1}}(x_{1})-T^{t_{m}}_{\omega_{2}}(x_{2})|
\geq \lambda^{m}|x_{1}-x_{2}| - \Cr{el}\,|F^{t_{m}}(\omega_{1})-F^{t_{m}}(\omega_{2})|.
\end{equation}
\end{lemma}

\begin{proof}
We define coefficients \( \Lambda_{1},\dots,\Lambda_{s}\) recursively by
\[
\Lambda_{1}:=1 \quad \text{ and } \quad \Lambda_{m+1}:= 1+ \lambda\theta^{-\tau_{m}}\Lambda_{m},
\]
and prove by induction on \( m\) that for all \( 1\leq m\leq s\) we have
\begin{equation}\label{eq:distortion81}
|T^{t_{m}}_{\omega_{1}}(x_{1})-T^{t_{m}}_{\omega_{2}}(x_{2})|
\geq \lambda^{m}|x_{1}-x_{2}| - 2\Cr{Cveh}\Lambda_{m}|F^{t_{m}}(\omega_{1})-F^{t_{m}}(\omega_{2})|.
\end{equation}
For \( m=1\) we have \( t_{1}=\tau_{0}=n+k\) and \( \Lambda_{1}=1\), and so \eqref{eq:distortion81} is exactly \eqref{eq:distortion4} with \( p=0\), which was already proved in Lemma \ref{lem:lower_bound_VE1}. Now suppose that \eqref{eq:distortion81} holds for some \( 1\leq m<s\). For \( j=1,2\) let \( \omega_{j}':=F^{t_{m}}(\omega_{j})\) and \( x_{j}':=T^{t_{m}}_{\omega_{j}}(x_{j})\), so that \( (\omega_{j}',x_{j}')=T^{m}_{\mathbf H}(\omega_{j},x_{j})\). Since \( m<s\), the two points \( (\omega_{1}',x_{1}'),(\omega_{2}',x_{2}')\) belong to the same element \( \mathscr{A}_{n',k',r'}\) of \( \mathscr R\), with \( n'+k'=\tau_{m}\), and therefore \eqref{eq:distortion4} with \( p=0\) applied to these two points gives
\[
|T^{\tau_{m}}_{\omega_{1}'}(x_{1}')-T^{\tau_{m}}_{\omega_{2}'}(x_{2}')|
\geq \lambda|x_{1}'-x_{2}'| - 2\Cr{Cveh}|F^{\tau_{m}}(\omega_{1}')-F^{\tau_{m}}(\omega_{2}')|.
\]
Since \( t_{m+1}=t_{m}+\tau_{m}\) we have \( T^{\tau_{m}}_{\omega_{j}'}(x_{j}')=T^{t_{m+1}}_{\omega_{j}}(x_{j})\) and \( F^{\tau_{m}}(\omega_{j}')=F^{t_{m+1}}(\omega_{j})\), and \( |x_{1}'-x_{2}'|\) is exactly the left hand side of \eqref{eq:distortion81}, and so the inductive hypothesis gives
\begin{equation}\label{eq:distortion9}
\begin{split}
|T^{t_{m+1}}_{\omega_{1}}(x_{1})-T^{t_{m+1}}_{\omega_{2}}(x_{2})|
\geq \lambda^{m+1}|x_{1}-x_{2}|
& - 2\Cr{Cveh}\lambda\Lambda_{m}|F^{t_{m}}(\omega_{1})-F^{t_{m}}(\omega_{2})| \\
& - 2\Cr{Cveh}|F^{t_{m+1}}(\omega_{1})-F^{t_{m+1}}(\omega_{2})|.
\end{split}
\end{equation}
Finally, \( \omega_{1}'\) and \( \omega_{2}'\) lie in the same element of \( \mathscr P_{\tau_{m}}\), and so by the uniform expansion of \( F\) we have
\[
|F^{t_{m}}(\omega_{1})-F^{t_{m}}(\omega_{2})| = |\omega_{1}'-\omega_{2}'|
\leq \theta^{-\tau_{m}}|F^{\tau_{m}}(\omega_{1}')-F^{\tau_{m}}(\omega_{2}')|
= \theta^{-\tau_{m}}|F^{t_{m+1}}(\omega_{1})-F^{t_{m+1}}(\omega_{2})|.
\]
Substituting this into \eqref{eq:distortion9} we get
\[
|T^{t_{m+1}}_{\omega_{1}}(x_{1})-T^{t_{m+1}}_{\omega_{2}}(x_{2})|
\geq \lambda^{m+1}|x_{1}-x_{2}|
- 2\Cr{Cveh}\big(1+\lambda\theta^{-\tau_{m}}\Lambda_{m}\big)|F^{t_{m+1}}(\omega_{1})-F^{t_{m+1}}(\omega_{2})|,
\]
which is \eqref{eq:distortion81} for \( m+1\) by the definition of \( \Lambda_{m+1}\). This completes the induction.

It just remains to bound the coefficients \( \Lambda_{m}\). Since \( \tau_{m}\geq 1\) for every \( m\), the recursive definition gives \( \Lambda_{m+1}\leq 1+(\lambda/\theta)\Lambda_{m}\), and therefore, by another immediate induction, \( \Lambda_{m}\leq \sum_{i=0}^{m-1}(\lambda/\theta)^{i} \leq \sum_{i=0}^{\infty}(\lambda/\theta)^{i}\) for all \( 1\leq m\leq s\). Substituting this into \eqref{eq:distortion81} gives \eqref{eq:distortion13} by the definition of \( \Cr{el}\).
\end{proof}

We can now prove Proposition \ref{pro:expansion}.

\begin{proof}[Proof of Proposition \ref{pro:expansion}]
The statements in \eqref{eq:exp12} and \eqref{eq:distortion4} have already been proved in Lemma \ref{lem:vexpansion} and \ref{lem:lower_bound_VE1} respectively, so we just need to prove \eqref{eq:eventual_expansion}. Let \( 1\leq m\leq s((\omega_{1},x_{1}),(\omega_{2}, x_{2}))\) and let \( t_{m}\) be as defined before Lemma \ref{lem:inductive_step2}.
Adding $\Cr{el}|F^{t_{m}}(\omega_{1})-F^{t_{m}}(\omega_{2})|$ to both sides of  \eqref{eq:distortion13} we have
\begin{equation}\label{eq:distortion1111}
\Cr{el}\big( |T^{t_{m}}_{\omega_{1}}(x_{1})-T^{t_{m}}_{\omega_{2}}(x_{2})| +|F^{t_{m}}(\omega_{1})- F^{t_{m}}(\omega_{2})|\big)\geq   {\lambda}^{m}|x_{1}-x_{2}|.
\end{equation}
By  definition of \( t_{m} \) we have   \( T^{m}_{\mathbf{H}}(\omega_{1},x_{1})=\widehat{F}^{t_{m}}(\omega_{1},x_{1})\) and \( T^{m}_{\mathbf{H}}(\omega_{2},x_{2})=\widehat{F}^{t_{m}}(\omega_{2},x_{2})\) and
\eqref{eq:distortion1111} gives
\begin{equation}\label{eq:distortion5}
     d(T^{m}_{\mathbf{H}}(\omega_{1},x_{1}),T^{m}_{\mathbf{H}}(\omega_{2},x_{2})) \geq \frac{\lambda^{m}}{\Cr{el}}|x_{1}-x_{2}|.
\end{equation}
Moreover, using the uniform expansivity of \( F \) and the fact that \( t_{m}\geq m \) and \( \lambda<\theta \), we also have
\begin{equation}\label{eq:distortion6}
     d(T^{m}_{\mathbf{H}}(\omega_{1},x_{1}),T^{m}_{\mathbf{H}}(\omega_{2},x_{2})) \geq \theta^{t_{m}}~|\omega_{1}-\omega_{2}| \geq \lambda^{m}|\omega_{1}-\omega_{2}|.
\end{equation}
Adding \eqref{eq:distortion5} and \eqref{eq:distortion6} we have the required result.
\end{proof}

\begin{remark}
Notice the crucial role of the Monotone Staircase Lemma \ref{lem:stairs_between_points} at two stages in this section: First in the proof of Sublemma \ref{lem:inductive_expansion_control} to estimate the precise bounds for \(\nu_{i}\), which are essential for the convergence of the series in \eqref{eq:seriesf}.  Second, it provides the factor \(\lambda\) in the lower bound for \(J_j\) in \eqref{eq:bound_by_IiJj}, which is needed to obtain the overall expansion. This is one of the places where this key innovation enables us to prove the expansion while still having the full branch property.
\end{remark}

\section{Distortion estimates}\label{sec:distortion}
To conclude that $T_{\mathbf{H}}$ is Gibbs-Markov, and thereby complete the proof of Proposition~\ref{pro:gibbs_markov_structure}, it remains to establish bounded distortion. Notice that in our case, the Jacobian of the induced map is the absolute value of the determinant of its differential. The distortion of the uniformly expanding base map $F$ is already controlled by \eqref{eq:Fbd}, so the main task is to control the distortion in the fibre direction. We proceed in two steps. First, we control the distortion of the iterates from the inducing region into a right neighbourhood of \(0\). Second, we control the distortion of the iterates from the right neighbourhood of \(0\) back to the inducing region. That is, we have the following result.

\begin{proposition}\label{pro:implication_of_distortion}
The following statements hold:
\begin{enumerate}
    \item For all $n \geq 1$, and $(\omega_{1}, x_{1}), (\omega_{2}, x_{2}) \in \Hat{Z}_{n}$, with $\omega_{1},\omega_{2} \in \Omega \in \mathscr{P}_{n}$, we have
   \begin{equation}\label{eq:dist11}
        \big| \ln(\frac{\partial}{\partial x}T^{n}_{\omega_{1}})(x_{1}) - \ln(\frac{\partial}{\partial x}T^{n}_{\omega_{2}})(x_{2}) \big| \lesssim  d(\widehat{F}^{n}(\omega_{1},x_{1}), \widehat{F}^{n}(\omega_{2},x_{2})).
    \end{equation}
   \item For all $k \geq 1$, and $(\omega_{1}, x_{1}), (\omega_{2}, x_{2}) \in \Hat{Y}_{0,k}$, with $\omega_{1},\omega_{2} \in \Omega \in \mathscr{P}_{k}$, we have
   \begin{equation}\label{eq:dist14}
        \big| \ln(\frac{\partial}{\partial x}T^{k}_{\omega_{1}})(x_{1}) - \ln(\frac{\partial}{\partial x}T^{k}_{\omega_{2}})(x_{2}) \big| \lesssim   d(\widehat{F}^{k}(\omega_{1},x_{1}), \widehat{F}^{k}(\omega_{2},x_{2})).
    \end{equation}
\end{enumerate}
\end{proposition}
In Section~\ref{sec:new_dist} we first show that Proposition \ref{pro:implication_of_distortion} implies the bounded distortion property~\eqref{eq:jacobian} and thus complete the proof of Proposition \ref{pro:gibbs_markov_structure}.  In Section \ref{sec:overview_of_the_proof} we reduce the  proof of Proposition~\ref{pro:implication_of_distortion} to that of two  further Propositions \ref{pro:hdistortion} and \ref{pro:vdistortion} which deal with the ``horizontal'' and ``vertical'' distortions respectively. Finally, in Sections \ref{sec:horizontle_distortion} and \ref{sec:vertical_distortion} we prove Propositions \ref{pro:hdistortion} and \ref{pro:vdistortion}.

\subsection{Proof of the Bounded Distortion condition  \eqref{eq:jacobian}}\label{sec:new_dist}
In this section we assume Proposition \ref{pro:implication_of_distortion} and prove the bounded distortion property~\eqref{eq:jacobian}, thus completing the proof of Proposition \ref{pro:gibbs_markov_structure}.
Since \(J_{\mathbf{H}}(\omega,x)=|det DT_{\mathbf{H}}(\omega,x)|\),   it is enough to prove that  for all \(n,k,r\) and $(\omega_{1},x_{1}), (\omega_{2},x_{2}) \in \mathscr{A}_{n,k,r}$ we have
   \begin{equation}\label{eq:distortion_gamma}
        \bigg|\ln\frac{\left|detDT_{\mathbf{H}}(\omega_{1},x_{1})\right|}{\left|detDT_{\mathbf{H}}(\omega_{2},x_{2})\right|}\bigg|
       \lesssim  d(T_{\mathbf{H}}(\omega_{1},x_{1}), T_{\mathbf{H}}(\omega_{2},x_{2})).
    \end{equation}
Indeed, then~\eqref{eq:jacobian} follows from the fact that  \( T_{\mathbf{H}}\) is expanding, recall Lemma \ref{lem:inducedexpanding},  and the first property of expanding maps, recall Definition \ref{def:expanding}. 
Notice first of all that we can split the distortion into two parts. 
Indeed  from \eqref{eq:det_mat} we have
\begin{equation}\label{eq:total1}
     \bigg|\ln\frac{\left|detDT_{\mathbf{H}}(\omega_{1},x_{1})\right|}{\left|detDT_{\mathbf{H}}(\omega_{2},x_{2})\right|}\bigg|
     \leq \left|\ln\frac{(F^{n+k})^{'}(\omega_{1})}{(F^{n+k})^{'}(\omega_{2})}\right|
     + \left|\ln \frac{\frac{\partial}{\partial x } T^{n+k}_{\omega_{1}}(x_{1})}{\frac{\partial}{\partial x } T^{n+k}_{\omega_{2}}(x_{2})}\right|.
\end{equation}
For the first term, the bound \eqref{eq:Fbd} and the definition of the metric in \eqref{eq:metricd} immediately give 
\begin{equation}\label{eq:distfirstterm}
\left|\ln\frac{(F^{n+k})^{'}(\omega_{1})}{(F^{n+k})^{'}(\omega_{2})}\right| 
\lesssim |F^{n+k}(\omega_{1}) - F^{n+k}(\omega_{2})| \leq d(T_{\mathbf{H}}(\omega_{1},x_{1}), T_{\mathbf{H}}(\omega_{2},x_{2})).
\end{equation}
The second term is the most problematic and we estimate it in a dedicated Lemma, the proof of which is exactly where we use Proposition \ref{pro:implication_of_distortion}.
\begin{lemma}
 For all \(n,k,r\) and $(\omega_{1},x_{1}), (\omega_{2},x_{2}) \in \mathscr{A}_{n,k,r}$ we have 
\begin{equation}\label{eq:dist4}
   \left|\ln \frac{\frac{\partial}{\partial x } T^{n+k}_{\omega_{1}}(x_{1})}{\frac{\partial}{\partial x } T^{n+k}_{\omega_{2}}(x_{2})}\right| \lesssim   d(T_{\mathbf{H}}(\omega_{1},x_{1}), T_{\mathbf{H}}(\omega_{2},x_{2})).
\end{equation}
\end{lemma}
\begin{proof}
By the chain rule we have 
\[
\frac{\partial}{\partial x } T^{n+k}_{\omega_{1}}(x_{1}) =  \frac{\partial}{\partial x } T^{n}_{\omega_{1}}(x_{1})  \frac{\partial}{\partial x } T^{k}_{F^{n}(\omega_{1})}(T^{n}_{\omega_{1}}(x_{1})) 
\text{ and } 
\frac{\partial}{\partial x } T^{n+k}_{\omega_{2}}(x_{2}) =  \frac{\partial}{\partial x } T^{n}_{\omega_{2}}(x_{2})  \frac{\partial}{\partial x } T^{k}_{F^{n}(\omega_{2})}(T^{n}_{\omega_{2}}(x_{2})) 
\]
and therefore, taking logs, the left hand side of \eqref{eq:dist4} is bounded by
\begin{equation} \label{eq:dist4a}
 \big|\ln \frac{\partial}{\partial x}T^{n}_{\omega_{1}}(x_{1})-\ln \frac{\partial}{\partial x}T^{n}_{\omega_{2}}(x_{2})\big|  +\big|\ln \frac{\partial}{\partial x}T^{k}_{F^{n}(\omega_{1})}(T^{n}_{\omega_{1}}(x_{1}))-\ln \frac{\partial}{\partial x}T^{k}_{F^{n}(\omega_{2})}(T^{n}_{\omega_{2}}(x_{2})) \big|,
\end{equation}
which are exactly the two terms bounded in \eqref{eq:dist11} and \eqref{eq:dist14} in Proposition \ref{pro:implication_of_distortion}. Indeed, 
notice that by the definition of the skew-product, 
\[
\widehat{F}^{n}(\omega_{1},x_{1}) = (F^{n}(\omega_{1}),T^{n}_{\omega_{1}}(x_{1}))
\quad \text{ and } \quad 
\widehat{F}^{n}(\omega_{2},x_{2}) = (F^{n}(\omega_{2}),T^{n}_{\omega_{2}}(x_{2}))
\]
and by the assumption that $(\omega_{1},x_{1}), (\omega_{2},x_{2}) \in \mathscr{A}_{n,k,r}$ it follows that 
\[
(F^{n}(\omega_{1}),T^{n}_{\omega_{1}}(x_{1})), (F^{n}(\omega_{2}),T^{n}_{\omega_{2}}(x_{2})) \in \hat{Y}_{0,k}
\quad \text{ and } \quad 
F^{n}(\omega_{1}), F^{n}(\omega_{2}) \in \Omega \in \mathscr{P}_{k}. 
\]
Therefore, using \eqref{eq:dist11} and \eqref{eq:dist14} to bound \eqref{eq:dist4a}, and applying this to \eqref{eq:dist4}, we get 
\begin{equation}\label{eq:dist141}
 \begin{aligned}
    \left|\ln \frac{\frac{\partial}{\partial x } T^{n+k}_{\omega_{1}}(x_{1})}{\frac{\partial}{\partial x } T^{n+k}_{\omega_{2}}(x_{2})}\right|  
   & \lesssim 
    d (\widehat{F}^{n}(\omega_{1},x_{1}),\widehat{F}^{n}(\omega_{2},x_{2}))+  d(\widehat{F}^{n+k}(\omega_{1},x_{1}), \widehat{F}^{n+k}(\omega_{2},x_{2})).
\\
& = 
    d (\widehat{F}^{n}(\omega_{1},x_{1}),\widehat{F}^{n}(\omega_{2},x_{2}))+   d(T_{\mathbf{H}}(\omega_{1},x_{1}), T_{\mathbf{H}}(\omega_{2},x_{2})).
\end{aligned}
\end{equation}
By the second property of expanding maps in Definition \ref{def:expanding} we have  \(     d (\widehat{F}^{n}(\omega_{1},x_{1}),\widehat{F}^{n}(\omega_{2},x_{2})) \lesssim d(T_{\mathbf{H}}(\omega_{1},x_{1}), T_{\mathbf{H}}(\omega_{2},x_{2}))\) and substituting this into \eqref{eq:dist141} completes the proof.  
\end{proof}

\begin{proof}[Proof of Bounded Distortion assuming Proposition \ref{pro:implication_of_distortion}]
Substituting \eqref{eq:distfirstterm} and \eqref{eq:dist4} into \eqref{eq:total1} we get~\eqref{eq:distortion_gamma} which, as mentioned above, implies the bounded distortion property~\eqref{eq:jacobian}.
\end{proof}

\subsection{Overview of the  proof of Proposition \ref{pro:implication_of_distortion}}\label{sec:overview_of_the_proof}
In this Section we prove Proposition \ref{pro:implication_of_distortion}   modulo the Monotone Staircase Lemma~\ref{lem:stairs_between_points} and two key technical Propositions \ref{pro:hdistortion} and~\ref{pro:vdistortion}. The key idea is to separate the distortion into a ``horizontal'' and a ``vertical'' component. To do this, notice first of all that 
by adding and subtracting $\ln \big(\frac{\partial}{\partial x}T^{n}_{\omega_{2}}(x_{1})\big)$, to the left side of \eqref{eq:dist11} we have  
\begin{equation}\label{eq:mix_ratio}
 \begin{split}
  \big| \ln(\frac{\partial}{\partial x}T^{n}_{\omega_{1}})(x_{1}) - \ln(\frac{\partial}{\partial x}T^{n}_{\omega_{2}})(x_{2}) \big|
   \leq & \big|\ln\big(\frac{\partial}{\partial x } T^{n}_{\omega_{1}}(x_{1})\big)- \ln\big(\frac{\partial}{\partial x } T^{n}_{\omega_{2}}(x_{1})\big)\big|
   \\
   & +\big|\ln\big(\frac{\partial}{\partial x } T^{n}_{\omega_{2}}(x_{1})\big)-\ln\big(\frac{\partial}{\partial x } T^{n}_{\omega_{2}}(x_{2})\big)\big|.
\end{split}
\end{equation}
Notice that the point $(\omega_{2},x_{1})$ might not be in the same strip in which $(\omega_{1},x_{1})$ and $(\omega_{2},x_{2})$ are, but we now have another advantage: in the first term on the right of~\eqref{eq:mix_ratio}, both the partial derivatives are evaluated at the same point $x_{1}$. In other words it is the distortion of two points with same vertical height. We refer to it as \emph{horizontal distortion} and one of the key results that we will use to bound this quantity is the following.

\begin{proposition}\label{pro:hdistortion}
For all  $n \geq 1$, and  $(\omega_{1}, x), (\omega_{2}, x) \in \Hat{Z}_{n}$, similarly for all $(\omega_{1},x), (\omega_{2},x) \in \Hat{Y}_{0,n}$, with $\omega_{1} ,\omega_{2} \in \Omega \in  \mathscr{P}_{n}$,
\begin{equation}\label{eq:ratio1} 
    \bigg|\ln( \frac{\partial}{\partial x}T^{n}_{\omega_{1}}(x))- \ln( \frac{\partial}{\partial x}T^{n}_{\omega_{2}}(x))\bigg| \lesssim |F^{n}(\omega_{1})-F^{n}(\omega_{2})|.
\end{equation}
\end{proposition}
Similarly, in the second term on the right of inequality \eqref{eq:mix_ratio}, we are differentiating the same map (in fact, the composition of maps). In other words, it is the distortion of points with the same horizontal position, we refer to it as \emph{vertical distortion} and one of the key results that we will use to bound this quantity is the following.

\begin{proposition}\label{pro:vdistortion}
For all $n \geq 1$, and $(\omega, x_{1}), (\omega, x_{2}) \in \Hat{Z}_{n}$, similarly for all  $(\omega,x_{1}), (\omega,x_{2}) \in \Hat{Y}_{0,n}$,
\begin{equation}\label{eq:distortion_ratioo}
        \bigg|\ln( \frac{\partial}{\partial x}T^{n}_{\omega}(x_{1}))- \ln( \frac{\partial}{\partial x}T^{n}_{\omega}(x_{2}))\bigg|\lesssim |T^{n}_{\omega}(x_{1})- T^{n}_{\omega}(x_{2})|.
    \end{equation}
\end{proposition}

\begin{remark}\label{rem:vdist12}
At first sight \eqref{eq:mix_ratio} seems to reduce Proposition~\ref{pro:implication_of_distortion} immediately to Propositions~\ref{pro:hdistortion} and~\ref{pro:vdistortion}. Indeed, if the auxiliary point \( (\omega_{2},x_{1}) \) belonged to \( \Hat{Z}_{n} \), then Proposition~\ref{pro:hdistortion} would apply to the points \( (\omega_{1},x_{1}), (\omega_{2},x_{1}) \) and bound the first term on the right of \eqref{eq:mix_ratio}, and Proposition~\ref{pro:vdistortion} would apply to the points \( (\omega_{2},x_{1}), (\omega_{2},x_{2}) \) and bound the second term.

The difficulty is that \( (\omega_{2},x_{1}) \) need not belong to \( \Hat{Z}_{n} \). The strip \( \Hat{Z}_{n} \) is bounded by the graphs of \( Z_{n-1,0} \) and \( Z_{n,0} \), which are not constant in \( \omega \), and so the horizontal segment joining \( (\omega_{1},x_{1}) \) to \( (\omega_{2},x_{1}) \) may leave \( \Hat{Z}_{n} \) and end in some other strip \( \Hat{Z}_{m} \) (with \( m>n \) if \( x_{1}>x_{2} \), say). Both propositions require the two points to lie in the same strip, and this is not merely a technical convenience: as we explain in Remarks~\ref{rem:hdistortion} and~\ref{rem:vdistortion} below, the constants in their proofs degenerate as the index of the strip tends to infinity, and so neither proposition extends to pairs of points in \( \Hat{Z}_{n} \) and \( \Hat{Z}_{m} \) with a constant independent of \( m \).

We overcome this using the Monotone Staircase Lemma~\ref{lem:stairs_between_points}. Instead of a single auxiliary point we join \( (\omega_{1},x_{1}) \) to \( (\omega_{2},x_{2}) \) by the finite chain of points \eqref{eq:p_enum}, all of which lie in \( \Hat{Z}_{n} \) and each of which differs from the next only in the horizontal or only in the vertical coordinate, so that Propositions~\ref{pro:hdistortion} and~\ref{pro:vdistortion} apply to each pair of consecutive points. The number of points in the chain depends on the two points and may be arbitrarily large, and so it is essential that the final bound does not depend on it. This is exactly where the monotonicity of the sequences \( \{\omega^{i}\} \) and \( \{x^{i}\} \) is used: the contributions of the horizontal steps add up exactly to \( |F^{n}(\omega_{1})-F^{n}(\omega_{2})| \), and the contributions of the vertical steps telescope to the vertical distance between the images of the two original points, up to horizontal correction terms which are controlled by Lemma~\ref{lem:expansion_contol_by_horizontal}. The details are given in the proof of Proposition~\ref{pro:implication_of_distortion} below. This is one of the key technical innovations of the paper, and allows us to obtain a full branch Gibbs-Markov induced map, in contrast with the induced map of \cite{Gou07}.
\end{remark}

In the rest of this section we show that Propositions \ref{pro:hdistortion} and \ref{pro:vdistortion} imply Proposition  \ref{pro:implication_of_distortion}. 

\begin{proof}[Proof of Proposition \ref{pro:implication_of_distortion}]
We only prove the first    estimate \eqref{eq:dist11} in the statement, \eqref{eq:dist14} follows by a completely  analogous argument. For simplicity, let $x_{1}\geq x_{2}$ and consider a sequence of points between \((\omega_{1},x_{1})\) and \((\omega_{2},x_{2})\) generated by Lemma \ref{lem:stairs_between_points}, i.e for some \(m\geq 1\), 
we have a sequence
\begin{equation}\label{eq:stairs}
(\omega_{1},x_{1}) =: (\omega^{0},x^{0}), (\omega^{1},x^{0}),(\omega^{1},x^{1}), \dots,   (\omega^{m}, x^{m-1}), (\omega^{m},x^{m}) := (\omega_{2},x_{2}).
\end{equation}
As discussed in Lemma \ref{lem:stairs_between_points}, these points form a monotone staircase, made up of horizontal and vertical lines segments, which connects the two points \((\omega_{1},x_{1}) \) and \((\omega_{2},x_{2}) \). We can therefore bound the distortion term in \eqref{eq:dist11}  by the sum of corresponding distortion terms between the endpoints of these segments.
For all $i = 0, \dots m-1$ and  $j= 1, \dots, m$, let 
\begin{equation*}
    I_{i}:= \big| \ln(\frac{\partial}{\partial x}T^{n}_{\omega^{i}})(x^{i}) - \ln(\frac{\partial}{\partial x}T^{n}_{\omega^{i+1}})(x^{i}) \big|, ~ \text{and} ~ J_{j}:=  \big| \ln(\frac{\partial}{\partial x}T^{n}_{\omega^{j}})(x^{j-1}) - \ln(\frac{\partial}{\partial x}T^{n}_{\omega^{j}})(x^{j}) \big|.
\end{equation*}
Notice that \( I_{0}\) is the distortion along the first horizontal segment joining the two points \( (\omega^{0}, x^{0}) \) and \( (\omega^{1}, x^{0} ) \)  and more generally \( I_{i}\) is the distortion along the horizontal segment joining the points \( (\omega^{i}, x^{i}) \) and \( (\omega^{i+1}, x^{i} ) \).  Similarly, \( J_{1}\) is the distortion along the first vertical segment joining the two points \( (\omega^{1}, x^{0}) \) and \( (\omega^{1}, x^{1} ) \)  and more generally \( J_{j}\) is the distortion along the vertical segment joining the points \( (\omega^{j}, x^{j-1}) \) and \( (\omega^{j}, x^{j} ) \). 
Therefore, 
adding and subtracting $\ln(\frac{\partial}{\partial x}T^{n}_{\omega^{i}})(x^{i-j})$ to the left side of \eqref{eq:dist11} for all $i=1, \dots, m$ and $j =0,1$  (analogously to what we did in \eqref{eq:distortion811}), by the standard triangle inequality
 we have
\begin{equation}\label{eq:distortion_implication}
        \big| \ln(\frac{\partial}{\partial x}T^{n}_{\omega_{1}})(x_{1}) - \ln(\frac{\partial}{\partial x}T^{n}_{\omega_{2}})(x_{2}) \big| \leq \sum_{i=0}^{m-1}I_{i}+ \sum_{j=1}^{m}J_{j}\toadd{.}
    \end{equation}
Splitting the distortion into these horizontal and vertical pieces puts us exactly into the setting of Proposition \ref{pro:hdistortion} and \ref{pro:vdistortion}. 
Indeed, by Proposition \ref{pro:hdistortion}, we have 
\begin{equation}\label{eq:distortion_implication1}
    \sum_{i=0}^{m-1} I_{i} \lesssim \sum_{i=0}^{m-1} |F^{n}(\omega^{i})-F^{n}(\omega^{i+1})| =|F^{n}(\omega_{1})-F^{n}(\omega_{2})|.
\end{equation}

For the vertical terms, notice that for each \( j \) the two points \( (\omega^{j},x^{j-1}) \) and \( (\omega^{j},x^{j}) \) lie on the same fibre and, by Monotone Staircase Lemma \ref{lem:stairs_between_points}, in \( \Hat{Z}_{n} \), and so Proposition \ref{pro:vdistortion} gives \( J_{j}\lesssim |T^{n}_{\omega^{j}}(x^{j-1})-T^{n}_{\omega^{j}}(x^{j})| \) with a constant which does not depend on \( j \) or \( m \). Moreover \( x^{j-1}\geq x^{j} \) and \( T^{n}_{\omega^{j}} \) is increasing, and so the quantity inside the absolute value is non-negative. Therefore
\begin{equation}\label{eq:distortion_improment2}
    \sum_{j=1}^{m}J_{j} \lesssim \sum_{j=1}^{m}\big(T^{n}_{\omega^{j}}(x^{j-1})-T^{n}_{\omega^{j}}(x^{j})\big).
\end{equation}
Adding and subtracting \(T^{n}_{\omega^{0}}(x^{0})\) and rearranging the right side of \eqref{eq:distortion_improment2} we get
\[
\sum_{j=1}^{m}\big(T^{n}_{\omega^{j}}(x^{j-1})-T^{n}_{\omega^{j}}(x^{j})\big)
= \big(T^{n}_{\omega^{0}}(x^{0})-T^{n}_{\omega^{m}}(x^{m})\big)
 - \sum_{i=0}^{m-1}\big(T^{n}_{\omega^{i}}(x^{i})-T^{n}_{\omega^{i+1}}(x^{i})\big).
\]
The first term on the right is the vertical distance between the images of the two original points. Each term in the second sum is the vertical displacement, after \( n \) iterates, of the two endpoints \( (\omega^{i},x^{i}) \) and \( (\omega^{i+1},x^{i}) \) of a horizontal segment of the staircase, which have the same vertical coordinate and \( \omega^{i},\omega^{i+1}\in\Omega \), and so by Lemma \ref{lem:expansion_contol_by_horizontal} it is bounded in absolute value by \( \Cr{Cveh}|F^{n}(\omega^{i})-F^{n}(\omega^{i+1})| \). Since the \( \omega^{i} \) are monotone and \( F^{n} \) is monotone on \( \Omega \), these bounds add up to \( \Cr{Cveh}|F^{n}(\omega_{1})-F^{n}(\omega_{2})| \), exactly as in \eqref{eq:distortion_implication1}. Taking absolute values we therefore get
\begin{equation}\label{eq:distortion_implication2}
    \sum_{j=1}^{m}J_{j} \lesssim |T^{n}_{\omega^{0}}(x^{0})-T^{n}_{\omega^{m}}(x^{m})|
    + \Cr{Cveh}|F^{n}(\omega_{1})-F^{n}(\omega_{2})|.
\end{equation}
Plugging \eqref{eq:distortion_implication1} and \eqref{eq:distortion_implication2} in \eqref{eq:distortion_implication} we have the required result \eqref{eq:dist11}. The case \(x_{1} <x_{2}\) follows similarly.
\end{proof}

It remains to prove the two technical Propositions \ref{pro:hdistortion} and \ref{pro:vdistortion}. The next two sections are dedicated to establish these results.

\subsection{Horizontal Distortion}\label{sec:horizontle_distortion}
In this subsection we prove Proposition \ref{pro:hdistortion}. We only prove the result for $(\omega_{1},x),(\omega_{2},x) \in \Hat{Z}_{n}$, and the case $(\omega_{1},x),(\omega_{2},x) \in \Hat{Y}_{0,n}$ follows by analogous reasoning. Notice that we can bound the horizontal distortion ratio as
    \begin{equation}\label{eq:hdistortion7}
    \begin{split}
   \bigg |\ln( \frac{\partial}{\partial x }T^{n}_{\omega_{1}}(x))- \ln( \frac{\partial}{\partial x }T^{n}_{\omega_{2}}(x))\bigg|
    &\leq \bigg|\ln (\frac{\partial}{\partial x}T_{\omega_{1}}(x))- \ln (\frac{\partial}{\partial x}T_{\omega_{2}}(x))\bigg|\\
    &+ \sum_{k=1}^{n-1}\bigg|\ln (\frac{\partial}{\partial x}T_{F^{k}(\omega_{1})}(T^k_{\omega_{1}}(x))) - \ln (\frac{\partial}{\partial x}T_{F^{k}(\omega_{2})}(T^{k}_{\omega_{2}}(x)))\bigg|.
 \end{split}
 \end{equation} 
The first term on the right corresponds to the distortion of the first iterate and we refer to it as \emph{distortion in the inducing region}, we bound this in  Lemma \ref{lem:hdistortion2} below. The second term corresponds to the iterates in the \(\mathbb S^{1} \times I^{+}\), we refer to it as \emph{distortion in the neutral neighbourhood}, the key step is to prove that the ``distortion function" is locally Lipschitz in the neighbourhood $\mathbb{S}^{1} \times U_{+1}$; this together with Lemma \ref{lem:expansion_contol_by_horizontal}, will lead us to the desired consequence.
The existence of the singular points makes the estimates in the inducing region weaker than the one in the neutral neighbourhood but good enough for the Proposition.

\subsubsection{Distortion in the inducing region}
We prove the following result.
\begin{lemma}\label{lem:hdistortion2}
For all $n\geq 1$ and for all $(\omega_{1},x),(\omega_{2},x) \in \Hat{Z}_{n}$  with $\omega_{1},\omega_{2} \in \Omega \in \mathscr{P}_{n}$, we have
\begin{equation}\label{eq:nothing}
     \bigg|\ln \frac{\partial}{\partial x}T_{\omega_{1}}(x)- \ln \frac{\partial}{\partial x}T_{\omega_{2}}(x)\bigg|\lesssim |F^{n}(\omega_{1})-F^{n}(\omega_{2})|
\end{equation}
\end{lemma}

\begin{proof}
We use a similar strategy as in Lemma \ref{lem:mean_value}. We let $g(\omega) = \ln \frac{\partial}{\partial x}T_{\omega}(x)$ and apply the Mean Value Theorem to \(g\) between \(\omega_1\) and \(\omega_2\). When \( n<\bar n\) the strip \( \Hat{Z}_{n}\) lies in a region in which the first and second derivatives of \(T_{\omega}(x)\) are uniformly bounded and, in particular the derivative of \(g\) is bounded and therefore we have 
\[
    \bigg|\ln \frac{\partial}{\partial x}T_{\omega_{1}}(x)- \ln \frac{\partial}{\partial x}T_{\omega_{2}}(x)\bigg| 
    =|g(\omega_1)-g(\omega_2)|
    \lesssim |\omega_{1}-\omega_{2}|
    \lesssim |F^{n}(\omega_{1})-F^{n}(\omega_{2})|.
\]
When \(n \geq \bar{n}\), then \(\hat Z_{n} \subset \mathbb S^{1} \times U_{-0}\) and therefore differentiating $g$ gives
\begin{equation}\label{eq:deriavtive_wrt_omega}
    g^{'}(\omega) = \frac{a^{'}_{1}(\omega)k_{1}(\omega)|x|^{k_{1}(\omega)-1}+a_{1}(\omega)k^{'}_{1}(\omega)|x|^{k_{1}(\omega)-1}+a_{1}(\omega)k_{1}(\omega)|x|^{k_{1}(\omega)-1}\ln|x|k^{'}_{1}(\omega)- \frac{\partial^{2}}{\partial \omega \partial x}\xi_{\omega}^{(1)}(x)}{a_{1}(\omega)k_{1}(\omega)|x|^{k_{1}(\omega)-1}-\frac{\partial}{\partial x}\xi_{\omega}^{(1)}(x)}.
\end{equation}
By assumption \eqref{eq:higherorder} we have
\begin{equation*}
    |g^{'}(\omega)| \leq   \frac{\left(|a^{'}_{1}(\omega)| k_{1}(\omega)+|k^{'}_{1}(\omega)|\big|1+\ln|x|k_{1}(\omega) \big|(a_{1}(\omega)+\delta)\right)|x|^{k_{1}(\omega)-1}}{a_{1}(\omega)k_{1}(\omega)|x|^{k_{1}(\omega)-1}}\lesssim |\ln|x||
\end{equation*}
and therefore
 \begin{equation}\label{eq:hdistortion1}
     \big|\ln \frac{\partial}{\partial x}T_{\omega_{1}}(x)- \ln \frac{\partial}{\partial x}T_{\omega_{2}}(x)\big|\lesssim |\ln |x|||\omega_{1}- \omega_{2}|.
    \end{equation}
Since \((\omega,x) \in \hat Z_{n}\), if $\ell_{2}^{-} \neq 0$ then by \eqref{eq:Z_n'} we have \(|\ln|x|\big| \lesssim \ln(n) \lesssim \theta^{n}\), similarly if $\ell_{2}^{-} =0$ then we have $\big|\ln|x|\big| \lesssim n\ln(1+b_{2}^{+}+ \delta)\lesssim \theta^{n}$. Therefore, \eqref{eq:hdistortion1} leads to
\begin{equation*}
    \bigg|\ln \frac{\partial}{\partial x}T_{\omega_{1}}(x)- \ln \frac{\partial}{\partial x}T_{\omega_{2}}(x)\bigg|\lesssim \theta^{n}|\omega_{1}-\omega_{2}|\lesssim|F^{n}(\omega_{1})-F^{n}(\omega_{2})|
\end{equation*}
\end{proof}

\subsubsection{Distortion in the neutral neighbourhood}
To simplify the notation in this section, we consider the function, $U: \mathbb{S}^{1}\times (0,1] \rightarrow [0, \infty)$ defined as
\begin{equation}\label{eq:fun_U}
    U(\omega,x) := \ln \frac{\partial}{\partial x}T_{\omega}(x).
\end{equation}
We prove that the map  \eqref{eq:fun_U}, is a Lipchitz function with the Lipchitz constant blowing up when we approach the top of the cylinder vertically, that is, on the set $\mathbb{S}^{1}\times \{1\}$. 

\begin{lemma}\label{lem:lipschitz}
For all $n \geq 1$  and  $(\omega_{1}, x_{1}), (\omega_{2}, x_{2}) \in \Hat{Y}_{n,0}$,
\begin{equation}\label{eq:hdist1}
    |U(\omega_{1},x_{1})-U(\omega_{2},x_{2})|\lesssim \frac{1}{(n+1)^{\frac{\ell_{2}^{-}-1}{\ell_{2}^{-}}}} d((\omega_{1},x_{1}),(\omega_{2}, x_{2})).
    \end{equation}
If $\ell_{2}^{-}=0$, then $\frac{1}{(n+1)^{\frac{\ell_{2}^{-}-1}{\ell_{2}^{-}}}}$ is replaced by $(1+b_{2}^{+}+\delta)^{n+1}$ in \eqref{eq:hdist1}.
\end{lemma}
\begin{proof}
Let $n\geq \bar{n}$, and $x_{1}<x_{2}$. Define \(h:[0,1] \rightarrow \mathbb{S}^{1} \times [0,1]\), as follows
    \begin{equation*}
    h(t) := (h_{\omega_{1},\omega_{2}}(t), x_{1}+ (x_{2}-x_{1})t),
    \end{equation*}
where $h_{\omega_{1},\omega_{2}}$ is a parametrization with constant speed of the shorter smooth arc joining \(\omega_{1}\) and \(\omega_{2}\), with $h_{\omega_{1},\omega_{2}}(0)=\omega_{1}$ and $h_{\omega_{1},\omega_{2}}(1) =\omega_{2}.$ Composing it with \(U\) we have \(U \circ h : [0,1] \rightarrow (0, \infty) \). Let \(\nabla\) be the gradient, then
\begin{equation}\label{eq:fun_U_MV}
    |U \circ h(0) -U \circ h(1)| \leq \sup_{0\leq t\leq 1}|\nabla U(h(t))| \int_{0}^{1}|h^{'}(\alpha)|d\alpha.
\end{equation}
If $\ell_{2}^{-} \neq 0$, then for all \((\omega,x) \in \mathbb{S}^{1}\times U_{+1}\),
\begin{equation*}
\bigg|\frac{\partial U}{\partial \omega}\bigg| \leq  \bigg| \frac{(b_{2}(\omega)+\delta)\ell_{2}^{'}(\omega)(1-x)^{\ell_{2}(\omega)}\big|\ln(1-x)(\ell_{2}(\omega)+1)+ 1\big|+b_{2}^{'}(\omega)(1-x)^{\ell_{2}(\omega)}(\ell_{2}(\omega)+1) }{1+b_{2}(\omega)(1-x)^{\ell_{2}(\omega)}(\ell_{2}(\omega)+1)}\bigg| \lesssim 1,
\end{equation*}
and
\begin{equation*}
\bigg|\frac{\partial U}{\partial x}\bigg|= \bigg| \frac{(b_{2}(\omega)+\delta)(1-x)^{\ell_{2}(\omega)-1}\ell_{2}(\omega)(\ell_{2}(\omega)+1)}{1+b_{2}(\omega)(1-x)^{\ell_{2}(\omega)}(\ell_{2}(\omega)+1)} \bigg| \lesssim (1-x)^{\ell_{2}^{-}-1}.
\end{equation*}
Using this in \eqref{eq:fun_U_MV}, noticing that \( (1-x)^{\ell_{2}^{-}-1} \) is increasing in \( x \) and that \( x\leq x_{2} \) along \( h \), so that the supremum is bounded by its value at height \( x_{2} \), and the fact that the integral is the arc-length of \( h \), which is at most \( |\omega_{1}-\omega_{2}|+|x_{1}-x_{2}| \), we have
\begin{equation*}
 |U \circ h(0) -U \circ h(1)| \lesssim (1-x_{2})^{\ell_{2}^{-}-1}d((\omega_{1}, x_{1}), (\omega_{2},x_{2}))
   \lesssim  \frac{1}{(n+1)^{\frac{\ell_{2}^{-}-1}{\ell_{2}^{-}}}}d((\omega_{1}, x_{1}), (\omega_{2},x_{2})).
\end{equation*}
In the last inequality, we have used the estimates from \eqref{eq:Y_n}. Similarly, if $\ell_{2}^{-}=0$, then both the partial derivatives of \(U\) can be bounded by \((1-x_{2})^{-1}\), therefore, by Lemma \ref{lem:estimats_X} we have
\begin{equation*}
     |U \circ h(0) -U \circ h(1)|\lesssim (1+b_{2}^{+}+\delta)^{n+1}d((\omega_{1}, x_{1}), (\omega_{2},x_{2})).
\end{equation*}
If $n<\bar{n}$, notice that uniform bounds on the partial derivatives of \(U\) are enough, which is the case here because in this case the strip \(\hat Y_{n,0}\) lies in a region in which the first and second partial derivatives of \(T_{\omega}(x)\) are uniformly bounded, in particular the partial derivatives of \(U\) are bounded. 
\end{proof}

\begin{lemma}\label{lem:hdistortion1}
For all $(\omega_{1}, x), (\omega_{2}, x) \in \Hat{Z}_{n}$, with $\omega_{1} ,\omega_{2} \in \Omega \in  \mathscr{P}_{n}$, and for all $k <n$ we have
     \begin{equation}\label{eq:rep}
       |U(\widehat{F}^{k}(\omega_{1}, x))-U(\widehat{F}^{k}(\omega_{2}, x))| \lesssim  \frac{1}{(n-k+1)^{\frac{\ell_{2}^{-}-1}{\ell_{2}^{-}}}\theta^{n-k}}|F^{n}(\omega_{1})-F^{n}(\omega_{2})|.
     \end{equation}
If $\ell_{2}^{-}=0,$ then  $\frac{1}{(n-k+1)^{\frac{\ell_{2}^{-}-1}{\ell_{2}^{-}}}}$ in \eqref{eq:rep} is replaced by $(1+b_{2}^{+}+\delta)^{n-k+1}$.
\end{lemma}
\begin{proof}
Notice that $\widehat{F}^{k}(\omega_{1}, x), \widehat{F}^{k}(\omega_{2}, x) \in \Hat{Y}_{n-k,0}$ . If $\ell_{2}^{-}\neq 0$ then using Lemma \ref{lem:lipschitz}, we get
\begin{equation}\label{eq:hdistortion4}
|U(\widehat{F}^{k}(\omega_{1}, x))-U(\widehat{F}^{k}(\omega_{2}, x))| \lesssim \frac{1}{(n-k+1)^{\frac{\ell_{2}^{-}-1}{\ell_{2}^{-}}}} d(\widehat{F}^{k}(\omega_{1},x), \widehat{F}^{k}(\omega_{2},x)).
\end{equation}
Recall the definition of \(d\) in \eqref{eq:metricd}, then by Lemma \ref{lem:expansion_contol_by_horizontal}, inequality \eqref{eq:Fbd} and assumption (A1) we have 
\begin{equation*}
    d(\widehat{F}^{k}(\omega_{1},x), \widehat{F}^{k}(\omega_{2},x))\lesssim |F^{k}(\omega_{1})-F^{k}(\omega_{2})|\lesssim \frac{1}{\theta^{n-k}}|F^{n}(\omega_{1})-F^{n}(\omega_{2})|.
\end{equation*}
Putting this in \eqref{eq:hdistortion4}, we get \eqref{eq:rep}. Similarly, if $\ell_{2}^{-}=0,$ then using Lemma \ref{lem:lipschitz}, the required result follows from a similar argument.
\end{proof}

\subsubsection{Proof of  Proposition \ref{pro:hdistortion}}
Finally, putting Lemma \ref{lem:hdistortion2} and \ref{lem:hdistortion1}, together we can conclude the Proposition.
\begin{proof}[Proof of the Proposition \ref{pro:hdistortion}]
For the first term on the right of  \eqref{eq:hdistortion7} we use Lemma \ref{lem:hdistortion2}. Therefore we only need to estimate the second term. Notice that,
\begin{equation}\label{eq:distortion8}
     \sum_{k=1}^{n-1}\bigg|\ln (\frac{\partial}{\partial x}T_{F^{k}(\omega_{1})}(T^k_{\omega_{1}}(x))) - \ln (\frac{\partial}{\partial x}T_{F^{k}(\omega_{2})}(T^{k}_{\omega_{2}}(x)))\bigg|
     = \sum_{k=1}^{n-1}|U(\widehat{F}^{k}(\omega_{1}, x))-U(\widehat{F}^{k}(\omega_{2}, x))|.
        \end{equation}
If $\ell_{2}^{-} \neq 0$, then by Lemma \ref{lem:hdistortion1},   
\begin{equation}\label{eq:vertical_distortion_1}
\begin{split}
  \sum_{k=1}^{n-1}|U(\widehat{F}^{k}(\omega_{1}, x))-U(\widehat{F}^{k}(\omega_{2}, x))|
    & \lesssim |F^{n}(\omega_{1})-F^{n}(\omega_{2})| \sum_{k=1}^{n-1} \frac{1}{(n-k+1)^{\frac{\ell_{2}^{-}-1}{\ell_{2}^{-}}}\theta^{n-k}}
    \\ &
   \leq  
   |F^{n}(\omega_{1})-F^{n}(\omega_{2})| 
   \sum_{k=1}^{\infty} (k+1)^{\frac{-\ell_{2}^{-}+1}{{\ell_{2}^{-}}}}\theta^{-k}<\infty
   \\ 
   & \lesssim |F^{n}(\omega_{1})-F^{n}(\omega_{2})|.
\end{split}
\end{equation}
If $\ell_{2}^{-}=0$, then
\[
 \sum_{k=1}^{n-1}|U(\widehat{F}^{k}(\omega_{1}, x))-U(\widehat{F}^{k}(\omega_{2}, x))|
     \lesssim |F^{n}(\omega_{1})-F^{n}(\omega_{2})| \sum_{k=1}^{\infty} \frac{(1+b_{2}^{+}+\delta)^{k+1}}{\theta^{k}}\lesssim  |F^{n}(\omega_{1})-F^{n}(\omega_{2})|,
\]
where the boundedness of the series on the right follows from (A5), in particular \((1+b_{2}^{+}+\delta)<\theta\). Notice that this is the only point in the proof where the second condition in \eqref{eq:expans_cond} is used, and only in the case \( \ell_{2}^{-}=0 \) (and similarly \( \ell_{1}^{-}=0 \) in the case of \( \Hat{Y}_{0,n} \)), recall Remark \ref{rem:expans_cond}. Therefore, using the above estimates and \eqref{eq:nothing} in \eqref{eq:hdistortion7} we have \eqref{eq:ratio1}.
\end{proof}

\begin{remark}\label{rem:hdistortion}
Recall that in Remark~\ref{rem:vdist12} we discussed the difficulty of making sure that the auxiliary point \( (\omega_{2},x_{1}) \) lies in the same strip as the others. In fact this is not strictly necessary for our arguments but if \( (\omega_{2},x_{1}) \) lies in some strip \( \hat Z_m\) it is necessary to have a uniform bound on \( m-n\) (which we achieve by making sure that \(m=n\)). We discuss briefly here why a uniform bound for \( m-n\) is necessary, in particular why it  cannot be allowed to lie in a strip \( \Hat{Z}_{m} \) with \( m\gg n \). 

In the proof of the Proposition \ref{pro:hdistortion}, the distortion accumulated in the neutral neighbourhood is controlled by the sum in \eqref{eq:vertical_distortion_1}, and this sum is bounded because of the interplay between two facts. First of all, by Lemma~\ref{lem:lipschitz}, the Lipschitz constant of \( U \) on \( \Hat{Y}_{j,0} \) grows with \( j \), like \( (j+1)^{(1-\ell_{2}^{-})/\ell_{2}^{-}} \) if \( \ell_{2}^{-}\neq 0 \) and like \( (1+b_{2}^{+}+\delta)^{j+1} \) if \( \ell_{2}^{-}=0 \). Second, by Lemma~\ref{lem:hdistortion1}, the \( k \)-th iterates of the two points are closer than their \( n \)-th iterates by the exponential factor \( \theta^{-(n-k)} \) coming from the uniform expansion of \( F \). The two facts are tied together because the \( k \)-th iterates of both points lie in the same strip \( \Hat{Y}_{n-k,0} \), and so the Lipschitz constant which applies at time \( k \) is \( (n-k+1)^{(1-\ell_{2}^{-})/\ell_{2}^{-}} \), which is dominated by \( \theta^{-(n-k)} \) (and if \( \ell_{2}^{-}=0 \) it is \( (1+b_{2}^{+}+\delta)^{n-k+1} \), which is dominated by \( \theta^{-(n-k)} \) precisely because of condition (A5)).

If instead \( (\omega_{1},x_{1})\in\Hat{Z}_{n} \) and \( (\omega_{2},x_{1})\in\Hat{Z}_{m} \) with \( m>n \), the \( k \)-th iterates lie in the different strips \( \Hat{Y}_{n-k,0} \) and \( \Hat{Y}_{m-k,0} \), and the only Lipschitz constant available for \( U \) on the region between them is that of \( \Hat{Y}_{m-k,0} \). The sum to be controlled would then be \( \sum_{k=1}^{n-1}(m-k+1)^{(1-\ell_{2}^{-})/\ell_{2}^{-}}\theta^{-(n-k)} \), whose last term alone (\( k=n-1 \)) is of order \( (m-n+2)^{(1-\ell_{2}^{-})/\ell_{2}^{-}} \) and is therefore unbounded in \( m \). No bound uniform in \( m \) is possible, unless \(m-n\) is a known fixed number (zero in our case).
\end{remark}

\subsection{Vertical distortion}\label{sec:vertical_distortion}
In this section, we prove Proposition \ref{pro:vdistortion}. We give the proof in the case $(\omega,x_{1}),(\omega,x_{2}) \in \Hat{Z}_{n}$,  the case $(\omega,x_{1}),(\omega,x_{2}) \in \Hat{Y}_{0,n}$ follows by an analogous reasoning.
In Section \ref{sec:Uni_bound_1} we prove the following  uniform bound on the distortion ratio in the neutral neighbourhood $\mathbb{S}^{1}\times U_{+1}$. 
\begin{proposition}\label{pro:vdistortion1}
   For all $n \geq 1$ and for all $(\omega,x_{1}),(\omega,x_{2}) \in \Hat{Y}_{n,0}\cup\Hat{Y}_{n+1,0}$,
    \begin{equation*}
        \big| \ln \big(\frac{\partial}{\partial x}T^{n} _{\omega}(x_{1}) \big)-\ln \big(\frac{\partial}{\partial x}T^{n}_{\omega}(x_{2}) \big) \big|\lesssim 1.
    \end{equation*}
\end{proposition}

In Section \ref{sec:Uni_bound_2} we prove the following uniform bound on the distortion ratio in the singular neighbourhood $\mathbb{S}^{1}\times U_{-0}$.
\begin{proposition}\label{pro:vdistortion2}
  For all $n \geq 1$ and  for all $(\omega,x_{1}),(\omega,x_{2}) \in \Hat{Z}_{n}$,
    \begin{equation}\label{eq:vdistoetion123}
        \big| \ln\big(\frac{\partial}{\partial x}T_{\omega}(x_{1}) \big)-\ln\big(\frac{\partial}{\partial x}T_{\omega}(x_{2}) \big)\big| \lesssim 1
    \end{equation}
\end{proposition}
Combining Proposition \ref{pro:vdistortion1} and \ref{pro:vdistortion2}, we get a uniform bound for the distortion ratio in \eqref{eq:distortion_ratioo}. However we need more than that and  in Section \ref{sec:imp_uni_bound} we improve this  bound to prove Proposition~\ref{pro:vdistortion}.

\subsubsection{Koebe Principle}\label{sec:kobe_p}
In the proof of both the above Propositions \ref{pro:vdistortion1} and \ref{pro:vdistortion2}, we use the Koebe Principle, and so  we recall some terminology and techniques from the literature, for details, we refer to \cite{dMvS93}. For a $C^{3}$ interval map $g:I \subset \mathbb{R} \rightarrow \mathbb{R}$ the \emph{Schwarzian derivative} at a point \(x \in I\) is given by
\begin{equation}
    Sg(x) = \frac{g^{'''}(x)}{g^{'}(x)}-\frac{3}{2} \Big( \frac{g^{''}(x)}{g^{'}(x)} \Big)^{2}, \quad \text{if} \quad g^{'}(x)\neq 0.
\end{equation}
The Schwarzian derivative acts on the composition of two functions $f$ and $g$ as follows.
\begin{equation}\label{eq:schwarzian_on_composition}
    S(g \circ f) = (Sg \circ f) (f^{'})^{2}+Sf.
\end{equation}
For intervals $U \subset V$, we say that $V$ is an $\epsilon$-\emph{scaled} neighbourhood of $U$ if both components of $V\setminus U$ have length at least $\epsilon |U|$.

The \emph{Koebe Principle} states that if $J \subset J^{'}$ are intervals and $g: J^{'} \rightarrow g(J^{'})$ is a $C^{3}$ diffeomorphism with $Sg <0$ such that $g(J^{'})$ contains an $\epsilon$-scaled neighbourhood of $g(J)$, then there exists a constant $\kappa$, depending only on $\epsilon$, such that for all \( x,y\in J \)
\begin{equation}\label{eq:koebe_principle2}
   |\ln g^{'}(x)-\ln g^{'}(y)| \leq \kappa \frac{|x-y|}{|J|}.
\end{equation}
We will apply  the Koebe Principle to intervals of a particular kind. For $n > \bar{n}$ we consider
\begin{equation}\label{eq:notation1}
    J(\omega):=[Y_{n-1,0}(\omega),Y_{n+1,0}(\omega)], \quad \text{and} \quad J^{'}(\omega):= [Y_{n-2,0}(\omega), 1].
\end{equation}
Notice that,
\begin{equation}\label{eq:schrr}
  J(\omega)\subset  J^{'}(\omega) = [Y_{n-2,0}(\omega),Y_{n-1,0}(\omega)]\cup J(\omega)\cup  \bigcup_{k\geq n+2}[Y_{k-1,0}(\omega),Y_{k,0}(\omega)].
\end{equation}

\begin{lemma}\label{lem:epsolon_scal_negh}
Let \(J\) and \(J^{'}\) be as defined in \eqref{eq:notation1}, then there exists an $\epsilon>0$, such that for all $n >\bar{n}+1$ and $\omega \in \mathbb{S}^{1}$, the interval $T^{n-\bar{n}-1}_{\omega}(J^{'}(\omega))$ contains an $\epsilon$-scale neighbourhood of $T^{n-\bar{n}-1}_{\omega}(J(\omega))$.
\end{lemma}
\begin{proof}
Recall Section \ref{sec:Auxi_partition}, $\Hat{Y}_{n,0}$ is the region between the graph of $Y_{n-1,0}$ and $Y_{n,0}$, moreover \(\widehat{F}^{n-\bar{n}-1}:\Hat{Y}_{n,0} \to \Hat{Y}_{\bar{n}+1,0} \quad \text{and} \quad \widehat{F}^{n-\bar{n}-1}:\Hat{Y}_{n-2,0} \to \Hat{Y}_{\bar{n}-1,0}\) are surjective, in particular, for $\widetilde{\omega}=F^{n-\bar{n}-1}(\omega),$ and \(i\geq n-2\), the following is a bijection,
\begin{equation*}
    T^{n-\bar{n}-1}_{\omega}: [Y_{i,0}(\omega),Y_{i+1,0}(\omega)] \to [Y_{i-n+\bar{n}+1,0}(\widetilde{\omega}), Y_{i-n+\bar{n}+2,0}(\widetilde{\omega})].
\end{equation*}
Therefore, by \eqref{eq:schrr} we have
\begin{equation*}
    T^{n-\bar{n}-1}_{\omega}(J^{'}(\omega))=[Y_{\bar{n}-1,0}(\widetilde{\omega}), Y_{\bar{n},0}(\widetilde{\omega})] \cup T^{n-\bar{n}-1}_{\omega}(J(\omega)) \cup [Y_{\bar{n}+2,0}(\widetilde{\omega}),1]
\end{equation*}
By Lemma \ref{lem:graphs_are_disjoint}, the graphs of $Y_{\bar{n}-1,0}$ and $Y_{\bar{n},0}$ are disjoint, they are also compact, therefore, there exists $\delta_{1}>0$ such that \(\min\{|Y_{\bar{n}-1,0}(\omega)-Y_{\bar{n},0}(\omega)|, \quad \omega \in \mathbb{S}^{1} \}\geq \delta_{1}\).
By the continuity of \( Y_{\bar n+2,0} \) and Lemma \ref{lem:graphs_are_disjoint} there exists a constant \(\delta_{2}>0\), such that $|1-Y_{\bar{n}+2,0}|\geq  \delta_{2}$, we let \(\Delta=\min\{\delta_{1},\delta_{2}\}\). Furthermore, notice that we can fix $\delta_{3}$ such that, \(\max\{|T^{n-\bar{n}-1}_{\omega}(J(\omega))|; \quad \omega \in \mathbb{S}^{1}\}= \max\{|Y_{\bar{n},0}(\omega)-Y_{\bar{n}+2,0}(\omega)|, \quad\omega \in \mathbb{S}^{1}\} \leq \delta_{3}\).
Therefore, for all $n >\bar{n}+1$ and  $\omega \in \mathbb{S}^{1}$ the neighbourhoods on the left and right of $T_{\omega}^{n-\bar{n}-1}(J)$ inside $T^{n-\bar{n}-1}_{\omega}(J^{'}(\omega))$ have size at least $\epsilon \left|T^{n-\bar{n}-1}_{\omega}(J(\omega))\right|,$ with \(\epsilon= \Delta/\delta_{3}\), see figure \ref{fig:epsilon_scale_nbhd}.
\begin{figure}[h]
    \centering
\includegraphics[width=0.85\linewidth]{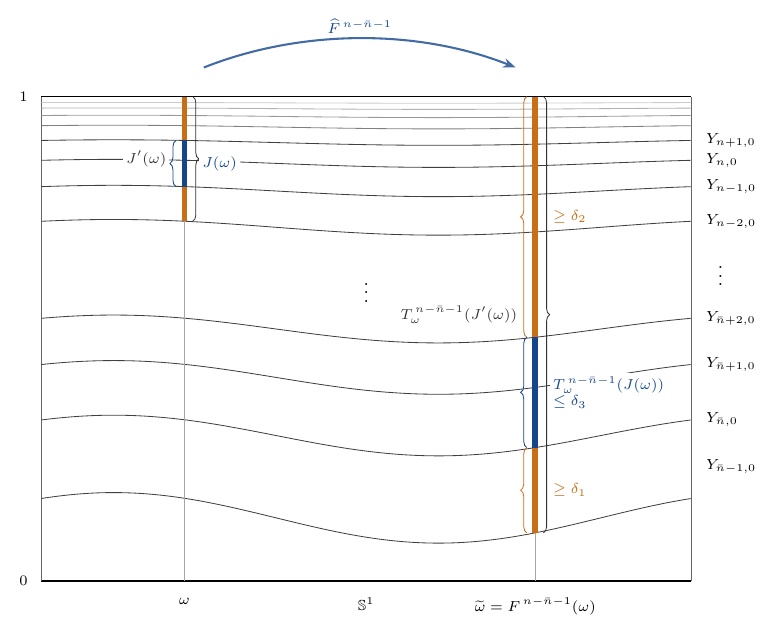}
    \caption{Illustration of the proof of Lemma \ref{lem:epsolon_scal_negh}. The intervals \( J(\omega)\subset J'(\omega) \) on the fibre over \( \omega \) and their images on the fibre over \( \widetilde\omega=F^{n-\bar n-1}(\omega) \). The image of \( J'(\omega) \) contains the two strips between \( Y_{\bar n-1,0} \) and \( Y_{\bar n,0} \) and between \( Y_{\bar n+2,0} \) and \( 1 \), whose widths are bounded below by \( \delta_{1} \) and \( \delta_{2} \) uniformly in \( \omega \), on either side of the image of \( J(\omega) \), whose length is at most \( \delta_{3} \).}
    \label{fig:epsilon_scale_nbhd}
\end{figure}
\end{proof}

\begin{remark}
    Notice that \(\Delta\) depends only on the skew product \(\widehat{F}\), in particular on \(\ell_{2},b, \delta\) and \(\bar{n}\). Therefore, \(\epsilon\) also depends only on \(\widehat{F}\).
\end{remark}

\begin{lemma}\label{lem:neg_shwarzian}
   For all $\omega \in \mathbb{S}^{1}$ and $n > \bar{n}+1$, if \(\ell_{2}(F^{i}(\omega))\neq 0\) and/or $\eta^{(2)}_{F^{i}(\omega)}$ is not identically zero for some $i < n-\bar{n}-1$, then the map
    \begin{equation}\label{eq:koebe1}   
    T^{n-\bar{n}-1}_{\omega}:  [Y_{n-2,0}(\omega),1] \to [Y_{\bar{n}-1,0}(F^{n-\bar{n}-1}(\omega)),1],
    \end{equation}
has negative Schwarzian derivatives.
\end{lemma}
\begin{proof}
Recalling the definition of \(T_{\omega}\) in \(U_{+1}\), we have
\begin{equation}\label{eq:schwar1}
\frac{ \frac{\partial^{3}}{\partial^{3} x}T_{\omega}(x)}{\frac{\partial}{\partial x} T_{\omega}(x)}= \frac{b_{2}(\omega)(\ell_{2}(\omega)+1)\ell_{2}(\omega)(\ell_{2}(\omega)-1)(1-x)^{\ell_{2}(\omega)-2}- \frac{\partial ^{3}}{\partial x^{3}}\eta_{\omega}^{(2)}(x)}{1+b_{2}(\omega)(\ell_{2}(\omega)+1)(1-x)^{\ell_{2}(\omega)}-\frac{\partial}{\partial x}\eta_{\omega}^{(2)}(x)},
\end{equation}
and
\begin{equation}\label{eq:schwar2}
    \frac{\frac{\partial^{2}}{\partial^{2} x}T_{\omega}(x)}{ \frac{\partial}{\partial x} T_{\omega}(x)}= \frac{ -b_{2}(\omega)(\ell_{2}(\omega)+1)\ell_{2}(\omega)(1-x)^{\ell_{2}(\omega)-1}-\frac{\partial^{2}}{\partial x^{2}}\eta_{\omega}^{(2)}(x)}{1+b_{2}(\omega)(\ell_{2}(\omega)+1)(1-x)^{\ell_{2}(\omega)}-\frac{\partial}{\partial x}\eta_{\omega}^{(2)}(x)}.
\end{equation}
Notice that using the above two expressions we have \(S(T_{\omega})(x)\leq 0\), we claim that the inequality is strict if \(\ell_{2}(\omega) \neq 0\) and/or \(\eta_{\omega}^{(2)}\neq 0\). First of all notice that, by (A5), in particular \eqref{eq:higherorder}, we have \(\frac{\partial}{\partial x}\eta_{\omega}^{(2)}(x) \leq \delta (\ell_{2}(\omega)+1)(1-x)^{\ell_{2}(\omega)}\) with \(\delta<b_{2}^{-}\), therefore, the denominator in the expressions \eqref{eq:schwar1} and \eqref{eq:schwar2} is positive.  Now if \(\ell_{2}(\omega)\neq 0\) then \(b_{2}(\omega)(\ell_{2}(\omega)+1)\ell_{2}(\omega)(\ell_{2}(\omega)-1)(1-x)^{\ell_{2}(\omega)-2} < 0\), and by (A4), in particular \eqref{eq:neg_schwar_eta}, \(\frac{\partial^{3}}{\partial x^{3}}\eta_{\omega}^{(2)}(x) \geq   0\). Similarly, if \(\ell_{2}(\omega)=0\) and \(\eta_{\omega}^{(2)} \neq 0\), then again by (A4), \(\frac{\partial^{3}}{\partial x^{3}}\eta_{\omega}^{(2)}(x) > 0\). From both of these observations we conclude that if \(\ell_{2}(\omega) \neq 0\) and/or \(\eta_{\omega}^{(2)}\neq 0\), then the numerator in \eqref{eq:schwar1} is negative, therefore the quantity in \eqref{eq:schwar1} is negative and  the claim follows.

We now use the claim to prove the Lemma. For all $0\leq j< n- \bar{n}-1$, and $\omega \in \mathbb{S}^{1}$, we have $T_{F^{j}(\omega)}([Y_{n-j-2,0}(F^{j}(\omega)),1])=[Y_{n-j-3,0}(F^{j+1}(\omega)),1]$ and $[Y_{n-j-2,0}(F^{j}(\omega)),1]\subset U_{+1}$. Therefore, $S(T_{F^{j}(\omega)})\leq 0$, and by the claim it is negative for \(j=i\). Recalling the action of the Schwarzian derivative on the composition \eqref{eq:schwarzian_on_composition}, we conclude that $S(T^{n-\bar{n}-1}_{\omega}) <0$.
\end{proof}

\subsubsection{Uniform bound on the distortion ratio in the neutral neighbourhood}\label{sec:Uni_bound_1}
We now prove a strong uniform distortion bound which in particular implies Proposition \ref{pro:vdistortion1}. 

\begin{proposition}\label{cla:strongdist}
For all \( n\geq 1 \) and \((\omega,x_{1}),(\omega,x_{2}) \in \Hat{Y}_{n,0} \cup \Hat{Y}_{n+1,0}\),
\begin{equation}\label{eq:strongdist}
     \big| \ln \big(\frac{\partial}{\partial x}T^{n} _{\omega}(x_{1}) \big)-\ln \big(\frac{\partial}{\partial x}T^{n}_{\omega}(x_{2}) \big) \big| \lesssim  \frac{|x_{1}-x_{2}|}{|Y_{n-1,0}(\omega) - Y_{n+1,0}(\omega)|}+  |T_{\omega}^{n}(x_{1})-T^{n}_{\omega}(x_{2})|
     \lesssim 1
\end{equation}
\end{proposition}
For the proof, notice that in the preceding subsection we have ensured that the map defined in \eqref{eq:koebe1} satisfies all the conditions required to implement the Koebe principle.

\begin{proof}[Proof of Proposition \ref{cla:strongdist}]

Let \(J\) and \(J^{'}\) be as defined in \eqref{eq:notation1}, notice that for  $(\omega,x_{1}),(\omega,x_{2}) \in \Hat{Y}_{n,0} \cup \Hat{Y}_{n+1,0}$, we have $(\omega,x_{1}),(\omega,x_{2}) \in \{\omega\} \times J(\omega) \subset \Hat{Y}_{n,0} \cup \hat{Y}_{n+1,0}$. 

Consider the case when \(n>\bar{n}+1\). Let $\ell_{2}(F^{i}(\omega)) \neq 0$, and/or \(\eta^{(2)}_{F^{i}(\omega)} \neq 0 \) for some $i < n-\bar{n}-1$. By Lemma \ref{lem:epsolon_scal_negh} we have that $T^{n-\bar{n}-1}_{\omega}(J^{'}(\omega))$ contains $\epsilon$-scaled neighbourhood of $T^{n-\bar{n}-1}_{\omega}(J(\omega))$, and by Lemma \ref{lem:neg_shwarzian} we have $S(T^{n-\bar{n}-1}_{\omega})<0$. Apply the Koebe Principle, particularly \eqref{eq:koebe_principle2}, to the map $T^{n-\bar{n}-1}_{\omega}: J^{'}(\omega) \to [Y_{\bar{n}-1,0},1]$, we have that there exists $\kappa(\epsilon)$, with $\epsilon>0$, and only depends on \( \epsilon \), and therefore only on \( \widehat F \), such that for all $x_{1},x_{2} \in J$ 
\begin{equation}\label{eq:distortion_in_right1}
    \big| \ln \big(\frac{\partial}{\partial x}T^{n-\bar{n}-1} _{\omega}(x_{1}) \big)-\ln \big(\frac{\partial}{\partial x}T^{n-\bar{n}-1}_{\omega}(x_{2}) \big) \big| \leq \kappa(\epsilon) \frac{|x_{1}-x_{2}|}{|J(\omega)|}.
\end{equation}
If $\ell_{2}(F^{i}(\omega)) = 0$ and \(\eta_{F^{i}(\omega)}=0\) for all $i <n-\bar{n}-1$, then the map $T^{n-\bar{n}-1}_{\omega}: J^{'}(\omega) \to [Y_{\bar{n}-1,0},1]$ is an affine map, and for all $ x \in J^{'}(\omega)$, we have $\frac{\partial}{\partial x}T^{n-\bar{n}-1}_{\omega}(x)= \left(1+b_{2}(F^{n-\bar{n}-2}(\omega))\right)\left(1+b_{2}(F^{n-\bar{n}-3}(\omega))\right) \cdot \cdot \cdot (1+b_{2}(\omega))$. Therefore,
\begin{equation*}
     \big| \ln \big(\frac{\partial}{\partial x}T^{n-\bar{n}-1} _{\omega}(x_{1}) \big)-\ln \big(\frac{\partial}{\partial x}T^{n-\bar{n}-1}_{\omega}(x_{2}) \big) \big|=0.
\end{equation*}
Now, for the remaining iterates, notice that in \(\cup_{1 \leq i\leq \bar{n}+2} \hat{Y}_{i,0}\) the map \(T_{\omega}(x)\) have uniformly bounded first and second partial derivatives with respect to \(x\). Moreover by Remark \ref{rem:derivative_ge1} the first partial derivative is bounded below by \(1\). Let $\widetilde{\omega}= F^{n-\bar{n}-1}(\omega)$, then by Mean Value Theorem
\begin{equation}\label{eq:distortion_in_right2}
\begin{split}
    \big| \ln\big(\frac{\partial}{\partial x}T^{\bar{n}+1}_{\widetilde{\omega}}\big(T_{\omega}^{n-\bar{n}-1}(x_{1})\big)\big)- \ln\big(\frac{\partial}{\partial x}T^{\bar{n}+1}_{\widetilde{\omega}}\big(T_{\omega}^{n-\bar{n}-1}(x_{2})\big)\big)\big| \lesssim |T_{\omega}^{n}(x_{1})-T^{n}_{\omega}(x_{2})|.
    \end{split}
\end{equation}
Finally, notice that
\begin{equation}\label{eq:final}
    \begin{split}
        \big| \ln \big(\frac{\partial}{\partial x}T^{n} _{\omega}(x_{1}) \big)-& \ln \big(\frac{\partial}{\partial x}T^{n}_{\omega}(x_{2}) \big) \big| \leq \big| \ln \big(\frac{\partial}{\partial x}T^{n-\bar{n}-1} _{\omega}(x_{1}) \big)-\ln \big(\frac{\partial}{\partial x}T^{n-\bar{n}-1}_{\omega}(x_{2}) \big) \big|\\
        &+\big| \ln\big(\frac{\partial}{\partial x}T^{\bar{n}+1}_{\widetilde{\omega}}\big(T_{\omega}^{n-\bar{n}-1}(x_{1})\big)\big)- \ln\big(\frac{\partial}{\partial x}T^{\bar{n}+1}_{\widetilde{\omega}}\big(T_{\omega}^{n-\bar{n}-1}(x_{2})\big)\big)\big|.
        \end{split}
\end{equation}
Putting the estimates from inequalities \eqref{eq:distortion_in_right1} and \eqref{eq:distortion_in_right2} in \eqref{eq:final} we have the required result. If \(n\leq \bar{n}+1\), then all the \(n\) iterates are in \(\cup_{1\leq i\leq \bar{n}+2}\hat{Y}_{i,0}\). Then again, in this region the partial derivatives are uniformly bounded and we conclude similarly  to~\eqref{eq:distortion_in_right2}. 
\end{proof}
As a consequence of Proposition \ref{pro:vdistortion1}, and inequality \eqref{eq:strongdist} we have the following Corollary.
\begin{corollary}\label{cor:improment2}
   For all $(\omega,x_{1}),(\omega,x_{2}) \in \Hat{Z}_{n}$,
 \begin{equation}\label{eq:improved_second_iterat}
         \big| \ln \big(\frac{\partial}{\partial x}T^{n-1} _{F(\omega)}(T_{\omega}(x_{1})) \big)-\ln \big(\frac{\partial}{\partial x}T^{n-1}_{F(\omega)}(T_{\omega}(x_{2})) \big) \big| 
     \lesssim |T^{n}_{\omega}(x_{1})-T^{n}_{\omega}(x_{2})|.
     \end{equation}
\end{corollary}
\begin{proof}
Notice that \(\left(F(\omega),T_{\omega}(x_{1})\right), \left(F(\omega),T_{\omega}(x_{2})\right) \in \Hat{Y}_{n-1,0}\cup\Hat{Y}_{n,0}\), therefore,  \eqref{eq:strongdist} leads to 
\begin{equation}\label{eq:imp1}
        \big| \ln \big(\frac{\partial}{\partial x}T^{n-1} _{F(\omega)}(T_{\omega}(x_{1})) \big)-\ln \big(\frac{\partial}{\partial x}T^{n-1}_{F(\omega)}(T_{\omega}(x_{2})) \big) \big|  \lesssim \frac{|T_{\omega}(x_{1})-T_{\omega}(x_{2})|}{|J(F(\omega))|}+ |T_{\omega}^{n}(x_{1})-T^{n}_{\omega}(x_{2})|,
\end{equation}
where \(J(F(\omega)) :=[Y_{n-2,0}(F(\omega)),Y_{n,0}(F(\omega))]\). By the Mean Value Theorem there are  \(\alpha_{1}, \alpha_{2} \in J(F(\omega))\) such that
\begin{equation}\label{eq:imp2}
    \frac{|T_{\omega}(x_{1})-T_{\omega}(x_{2})|}{|J(F(\omega))|} =\frac{\frac{\partial}{\partial x}T^{n-1}_{F(\omega)}(\alpha_{2})}{\frac{\partial}{\partial x}T^{n-1}_{F(\omega)}(\alpha_{1})} \frac{|T^{n}_{\omega}(x_{1})-T^{n}_{\omega}(x_{2})|}{|T^{n-1}_{F(\omega)}(J(F(\omega)))|} \lesssim \frac{|T^{n}_{\omega}(x_{1})-T^{n}_{\omega}(x_{2})|}{|T^{n-1}_{F(\omega)}(J(F(\omega))|},
\end{equation}
where the last inequality above follows by noticing that $(F(\omega),\alpha_{1}),(F(\omega),\alpha_{2}) \in \Hat{Y}_{n-1,0}\cup\Hat{Y}_{n,0}$, and then using Proposition \ref{pro:vdistortion1}.
Notice that $T^{n-1}_{F(\omega)}(J(F(\omega)))=[0,Y_{1,0}(F^{n}(\omega))]$, and the graph of \(Y_{1,0}\) is at a positive distance from \(0\), therefore the denominator on the right side of \eqref{eq:imp2} is uniformly bounded away from \(0\). Combining \eqref{eq:imp2} and \eqref{eq:imp1}, we have the required result.
\end{proof}
\subsubsection{Proof of Proposition \ref{pro:vdistortion2}}\label{sec:Uni_bound_2}
In this  section we prove Proposition \ref{pro:vdistortion2}. Notice that the Schwarzian derivative of the maps in the fibres is not negative in the neighbourhood $U_{-0}$; therefore, we cannot use the Koebe principle directly. Recall Section \ref{sec:markov_partition}, the strips $\Hat{Z}_{n}$ are the union of countably many strips $\{\Hat{Z}_{n,k}\}_{k\in \mathbb{N}}$, and each strip $\Hat{Z}_{n,k}$ is a pull back of the strip $\Hat{Y}_{n-1,k}$, in particular $\Hat{Z}_{n}$ is a pull back of $\Hat{Y}_{n-1,0}$. Therefore, first in Lemma \ref{lem:verticaDl11} and \ref{lem:verticaDl12} we establish a bound depending on $\{Y_{n,0}\}$ and then we use the Koebe Principle to get the desired result.

\begin{lemma}\label{lem:verticaDl11}
  For all $n\geq 1$, and $(\omega,x_{1}),(\omega,x_{2}) \in \Hat{Z}_{n}$, 
    \begin{equation}\label{eq:verticaD111}
        \big| \ln\big(\frac{\partial}{\partial x}T_{\omega}(x_{1}) \big)-\ln\big(\frac{\partial}{\partial x}T_{\omega}(x_{2}) \big)\big| \lesssim  \frac{|x_{1}-x_{2}|}{|Z_{n,0}(\omega)|} \lesssim \frac{|Z_{n,0}(\omega)-Z_{n-1,0}(\omega)|}{|Z_{n,0}(\omega)|}.
    \end{equation}
    \end{lemma}
\begin{proof}
By the Mean Value Theorem, for some \(x \in [x_{1},x_{2}]\), we have
\begin{equation}\label{eq:dist3}
        \big|\ln\big( \frac{\partial}{\partial x}T_{\omega}(x_{1})\big)- \ln\big(\frac{\partial}{\partial x}T_{\omega}(x_{2})\big) \big| \leq \left|\frac{\frac{\partial^{2}}{\partial x^{2}}T_{\omega}(x)}{\frac{\partial}{\partial x }T_{\omega}(x)}\right| |x_{1}-x_{2}|.
  \end{equation}
If \( n\geq\bar n \), so that \( x\in U_{-0} \) and \( T_{\omega} \) has the explicit form \eqref{eq:definition_of_Tomega}, therefore
\begin{equation*}  
\left|\frac{\frac{\partial^{2}}{\partial x^{2}}T_{\omega}(x)}{\frac{\partial}{\partial x }T_{\omega}(x)}\right| \leq  \frac{(a_{1}(\omega)+\delta)|k_{1}(\omega)-1|k_{1}(\omega)|x|^{k_{1}(\omega)-2}}{a_{1}(\omega)k_{1}(\omega)|x|^{k_{1}(\omega)-1}} \lesssim \frac{1}{|x|}.
\end{equation*}
Putting this in \eqref{eq:dist3} and noticing that \(x \in [Z_{n-1,0}(\omega),Z_{n,0}(\omega)]\), we have \eqref{eq:verticaD111} for \( n\geq\bar n \). If \( n<\bar n \), then \(\hat Z_{n}\) is in a region in which the partial derivatives of \(T_{\omega}(x)\) are uniformly bounded, and \( \frac{\partial}{\partial x}T_{\omega}\geq ~1 \), recall Remark \ref{rem:derivative_ge1}. Therefore, using \eqref{eq:dist3} then the left hand side of \eqref{eq:verticaD111} is \( \lesssim|x_{1}-x_{2}| \), hence \eqref{eq:verticaD111} follows in this case well.
  \end{proof}

\begin{lemma}\label{lem:verticaDl12}
For all $n\geq  \bar{n}+1$, and $\omega \in \mathbb{S}^{1}$,
     \begin{equation}\label{eq:121}
         \frac{|Z_{n,0}(\omega)-Z_{n-1,0}(\omega)|}{|Z_{n,0}(\omega)|}\lesssim 1+\bigg(1+\frac{|Y_{n-1,0}(F(\omega))-Y_{n-2,0}(F(\omega))|}{|Y_{n,0}(F(\omega))-Y_{n-1,0}(F(\omega))|}\bigg)^{\frac{1}{k_{1}(\omega)}}
     \end{equation}
\end{lemma}
\begin{proof}
Recall that \(T_{\omega}(Z_{i,0}(\omega))=Y_{i-1,0}(F(\omega))\). Using the definition of \(T_{\omega}|_{U_{-0}}\), rearranging and then using assumption (A5), particularly \eqref{eq:higherorder}, we have that for all $i\geq \bar{n}$ 
\begin{equation*}
    |Z_{i,0}(\omega)|\leq \left(\frac{|1-Y_{i-1,0}(F(\omega))|}{a_{1}(\omega)}\right)^{\frac{1}{k_{1}(\omega)}} \quad \text{and} \quad  |Z_{i,0}(\omega)|\geq \left(\frac{|1-Y_{i-1,0}(F(\omega))|}{a_{1}(\omega)+\delta}\right)^{\frac{1}{k_{1}(\omega)}}.
\end{equation*}
Using this on the left side of \eqref{eq:121} for $i=n$ and $n-1$, we have
\begin{equation}\label{eq:Zn2}
    \frac{|Z_{n,0}(\omega)-Z_{n-1,0}(\omega)|}{|Z_{n,0}(\omega)|}\leq 1+ \frac{|Z_{n-1,0}(\omega)|}{|Z_{n,0}(\omega)|}\lesssim 1+ \left(\frac{|1-Y_{n-2,0}(F(\omega))|}{|1-Y_{n-1,0}(F(\omega))|}\right)^{\frac{1}{k_{1}(\omega)}}.
\end{equation} 
By adding and subtracting $Y_{n-1,0}(F(\omega))$ to the numerator of the ratio on the right of \eqref{eq:Zn2} we have,
\begin{equation*}
     \frac{|1-Y_{n-2,0}(F(\omega))|}{|1-Y_{n-1,0}(F(\omega))|} \leq 1+ \frac{|Y_{n-1,0}(F(\omega))-Y_{n-2,0}(F(\omega))|}{|1-Y_{n-1,0}(F(\omega))|}\leq 1+  \frac{|Y_{n-1,0}(F(\omega))-Y_{n-2,0}(F(\omega))|}{|Y_{n,0}(F(\omega))-Y_{n-1,0}(F(\omega))|},
\end{equation*}
where the last inequality follows because, by Lemma \ref{lem:graphs_are_disjoint}, the sequence $\{Y_{n,0}(F(\omega))\}$ is monotone increasing and converges $1$, therefore $|1-Y_{n-1,0}(F(\omega))|\geq |Y_{n,0}(F(\omega))-Y_{n-1,0}(F(\omega))|$. Thus, putting this in \eqref{eq:Zn2}, we get the required result.
\end{proof}

\begin{lemma}\label{lem:verticaDl13}
For all $n>  1$, and $\omega \in \mathbb{S}^{1}$ ,
    \begin{equation}\label{eq:1212}
         \frac{|Y_{n-1,0}(F(\omega))-Y_{n-2,0}(F(\omega))|}{|Y_{n,0}(F(\omega))-Y_{n-1,0}(F(\omega))|} \lesssim 1.
    \end{equation}
\end{lemma}
\begin{proof}
Notice that by the repeated application of \eqref{eq:Xn}, we have \( T^{n-2}_{F(\omega)}(Y_{n-2,0}(F(\omega))) = Y_{0,0}(F^{n-1}(\omega))\), \(T^{n-2}_{F(\omega)}(Y_{n-1,0}(F(\omega))) = Y_{1,0}(F^{n-1}(\omega))\),
    and \(T^{n-2}_{F(\omega)}(Y_{n,0}(F(\omega))) =Y_{2,0}(F^{n-1}(\omega))\). Therefore, there exists $\alpha_{1},\alpha_{2} \in [Y_{n-2,0}(F(\omega)),Y_{n,0}(F(\omega))]$, such that
\begin{equation}\label{eq:above}
    \frac{|Y_{n-1,0}(F(\omega))-Y_{n-2,0}(F(\omega))|}{|Y_{n,0}(F(\omega))-Y_{n-1,0}(F(\omega))|}=\frac{\frac{\partial}{\partial x} T^{n-2}_{F(\omega)}(\alpha_{2})}{\frac{\partial}{\partial x} T^{n-2}_{F(\omega)}(\alpha_{1})}\frac{|Y_{1,0}(F^{n-1}(\omega))-Y_{0,0}(F^{n-1}(\omega))|}{|Y_{2,0}(F^{n-1}(\omega))-Y_{1,0}(F^{n-1}(\omega))|}.
\end{equation}
Notice that \((F(\omega),\alpha_{1}),(F(\omega),\alpha_{2}) \in \hat{Y}_{n-1,0} \cup \hat{Y}_{n,0}\), therefore, by Proposition \ref{pro:vdistortion1}, the ratio 
\begin{equation}\label{eq:lemma33}
\frac{|\frac{\partial}{\partial x} T^{n-1}_{F(\omega)}(\alpha_{2})|}{|\frac{\partial}{\partial x} T^{n-1}_{F(\omega)}(\alpha_{1})|}=\frac{\left|\frac{\partial}{\partial x} T_{F^{n-1}(\omega)}(T^{n-2}_{F(\omega)}(\alpha_{2}))\right|}{\left|\frac{\partial}{\partial x} T_{F^{n-1}(\omega)}(T^{n-2}_{F(\omega)}(\alpha_{1}))\right|}\frac{\left|\frac{\partial}{\partial x} T^{n-2}_{F(\omega)}(\alpha_{2})\right|}{\left|\frac{\partial}{\partial x} T^{n-2}_{F(\omega)}(\alpha_{1})\right|}
\end{equation}
is uniformly bounded. Since \((F^{n-1}(\omega),T^{n-2}_{F(\omega)}(\alpha_{2})), (F^{n-1}(\omega),T^{n-2}_{F(\omega)}(\alpha_{1})) \in \hat{Y}_{1,0} \cup \hat{Y}_{2,0}\) and the map \(T_{F^{n-1}(\omega)}(x)\) is \(C^{2}\) in the region \(\hat{Y}_{1,0} \cup \hat{Y}_{2,0}\), therefore, the first ratio on the right of the equation \eqref{eq:lemma33} is uniformly bounded and hence the second as well. By Lemma \ref{lem:graphs_are_disjoint} and the continuity of $Y_{n,0}$, the difference \(|Y_{2,0}(\omega)-Y_{1,0}(\omega)|\) is uniformly bounded away from \(0\). Using these two observations in \eqref{eq:above} we have the required result.
\end{proof}
\begin{corollary}\label{cor:bounded_ratio}
For all \(n\geq 1\) and \(\omega \in \mathbb{S}^{1},\)
     \begin{equation}\label{eq:rem1}
        \frac{|Z_{n,0}(\omega)-Z_{n-1,0}(\omega)|}{|Z_{n,0}(\omega)|} \lesssim 1.
    \end{equation}
\end{corollary}
\begin{proof}
Using \eqref{eq:1212} in \eqref{eq:121}, then for all $n \geq \bar{n}+2$, the statement follows.
For \(n<\bar{n}+2\), notice that the numerator is bounded by \(1\) and recall Remark \ref{rem:Znk}, we have \(|Z_{n,0}(\omega)|> |Z_{\bar{n}+2,0}(\omega)|\), and the latter is uniformly bounded away from \(0\), hence \eqref{eq:rem1} follows.
\end{proof}

\begin{proof}[Proof of Proposition \ref{pro:vdistortion2}]
By putting \eqref{eq:rem1} in \eqref{eq:verticaD111}, we get the required result. 
\end{proof}

\subsubsection{Proof of Proposition \ref{pro:vdistortion}}\label{sec:imp_uni_bound}
In this subsection, we give a proof of Proposition \ref{pro:vdistortion}. First of all, by using Proposition \ref{pro:vdistortion1} and \ref{pro:vdistortion2}, we improve the uniform bound in \eqref{eq:vdistoetion123}, then we combine this improved bound with Corollary~\ref{cor:improment2}.

\begin{lemma}\label{lem:improment1}
For all $(\omega,x_{1}),(\omega,x_{2}) \in \Hat{Z}_{n}$,
   \begin{equation}\label{eq:distortion_in_right3}
 \big|\ln\big( \frac{\partial}{\partial x}T_{\omega}(x_{1})\big)-\ln\big(\frac{\partial}{\partial x}T_{\omega}(x_{2})\big) \big| \lesssim |T^{n}_{\omega}(x_{1})-T^{n}_{\omega}(x_{2})|.
\end{equation}
\end{lemma}
\begin{proof}
By the repeated application of \eqref{eq:Z_Y}, we have \(T^{n}_{\omega}: [Z_{n-1,0}(\omega),Z_{n,0}(\omega)] \rightarrow [0,Y_{0,0}(F^{n}(\omega))]\),
is a bijection. Therefore, for some $\alpha_{1},\alpha_{2} \in [Z_{n-1,0}(\omega),Z_{n,0}(\omega)]$,
\begin{equation}\label{eq:dist1}
    \frac{|x_{1}-x_{2}|}{|Z_{n-1,0}(\omega)-Z_{n,0}(\omega)|} =\frac{|T^{n}_{\omega}(x_{1})- T^{n}_{\omega}(x_{2})|}{|Y_{0,0}(F^{n}(\omega))|}\frac{\frac{\partial}{\partial x}T^{n-1}_{F(\omega)}(T_{\omega}(\alpha_{1}))}{\frac{\partial}{\partial x}T^{n-1}_{F(\omega)}(T_{\omega}(\alpha_{2}))} \frac{\frac{\partial}{\partial x}T_{\omega}(\alpha_{1})}{\frac{\partial}{\partial x}T_{\omega}(\alpha_{2})} .
\end{equation}
Notice that $(\omega,\alpha_{1}),(\omega,\alpha_{2}) \in \Hat{Z}_{n}$, and $(F(\omega),T_{\omega}(\alpha_{1})),(F(\omega),T_{\omega}(\alpha_{2})) \in  \Hat{Y}_{n-1,0}\cup\Hat{Y}_{n,0}$, therefore, by Proposition \ref{pro:vdistortion1} and \ref{pro:vdistortion2}, the ratios of the partial derivatives in \eqref{eq:dist1} are bounded above. Moreover the graph of \(Y_{0,0}\) is  at a positive distance from \(0\), therefore \eqref{eq:dist1} leads to
\begin{equation}
   |x_{1}-x_{2}| \lesssim |T^{n}_{\omega}(x_{1})- T^{n}_{\omega}(x_{2})||Z_{n-1,0}(\omega)-Z_{n,0}(\omega)|,
\end{equation}
using this in the first inequality in \eqref{eq:verticaD111} we have
\begin{equation}\label{eq:improment2}
     \big|\ln\big( \frac{\partial}{\partial x}T_{\omega}(x_{1})\big)-\ln\big(\frac{\partial}{\partial x}T_{\omega}(x_{2})\big)\big|
     \lesssim |T^{n}_{\omega}(x_{1})- T^{n}_{\omega}(x_{2})|\frac{|Z_{n-1,0}(\omega)-Z_{n,0}(\omega)|}{|Z_{n,0}(\omega)|}.
\end{equation}
Using Corollary \ref{cor:bounded_ratio} in \eqref{eq:improment2}, we have the required result.
\end{proof}

\begin{proof}[Proof of Proposition \ref{pro:vdistortion}]
Using the chain rule with \eqref{eq:distortion_in_right3} and \eqref{eq:improved_second_iterat}, we have the required result.
\end{proof}
\begin{remark}\label{rem:vdistortion}
We now give the analogous explanation for the vertical distortion, recall Remark~\ref{rem:vdist12}. Suppose that \( (\omega_{2},x_{2})\in\Hat{Z}_{n} \) but \( (\omega_{2},x_{1})\in\Hat{Z}_{m} \) for some \( m\gg n \), and that we try to prove Proposition~\ref{pro:vdistortion} for this pair. After one iterate we have \( \widehat{F}(\omega_{2},x_{2})\in\Hat{Y}_{n-1,0} \) and \( \widehat{F}(\omega_{2},x_{1})\in\Hat{Y}_{m-1,0} \), and so the argument of Corollary~\ref{cor:improment2} would require the uniform bound of Proposition~\ref{pro:vdistortion1} for pairs of points in two strips whose indices differ by \( m-n \), which can be arbitrarily large, rather than in two adjacent strips. In the notation of Proposition~\ref{pro:vdistortion1} this means replacing the interval \( J(\omega) \) in \eqref{eq:notation1} by \( [Y_{n-1,0}(\omega),Y_{m,0}(\omega)] \). Proposition~\ref{pro:vdistortion1} is proved by applying the Koebe Principle to the intervals \( J(\omega)\subset J^{'}(\omega) \), and the constant it produces depends only on the \( \epsilon \) of Lemma~\ref{lem:epsolon_scal_negh}. With the larger interval \( J(\omega) \), the right hand component of \( T^{n-\bar{n}-1}_{\omega}(J^{'}(\omega))\setminus T^{n-\bar{n}-1}_{\omega}(J(\omega)) \) would be \( [Y_{m-n+\bar{n}+1,0}(\widetilde{\omega}),1] \), whose length tends to zero as \( m\to\infty \) since \( Y_{j,0}\to 1 \). Thus, \( T^{n-\bar{n}-1}_{\omega}(J^{'}(\omega)) \) would not contain an \( \epsilon \)-scaled neighbourhood of \( T^{n-\bar{n}-1}_{\omega}(J(\omega)) \) for any \( \epsilon \) independent of \( m \), and the Koebe Principle would give no uniform bound. This is why Proposition~\ref{pro:vdistortion1} is stated for two adjacent strips only.
\end{remark}

\section{Tail Estimates}\label{sec:tail} 
In this section, we prove Proposition \ref{pro:tail_lebesgue}. Recall \(d\hat{\mu}/d m\) is bounded above and below, therefore, it suffices to prove \eqref{eq:tail_lebesgue1} and \eqref{eq:tail_lebesgue2} for \(m\).
Recall  first of all that 
\[
\mathbf{H} = \bigcup_{i\geq 1} \bigcup_{k>1}\Hat{Z}_{i,k}
\quad\text{ and } \quad
 \{\tau >n \} =\bigcup_{i\geq 1} \bigcup_{k>n-i}\Hat{Z}_{i,k}.
\]
For the purposes of our estimate it is useful to divide the right hand side into disjoint subsets with \( i=1\) and \(i>n\) grouped together and \( 2\leq i\leq n\) grouped together. Thus we write 
\[
\mathcal{R}_{n}:= \big(\bigcup_{k\geq n} \Hat{Z}_{1,k}\big) 
\cup
\big(\bigcup_{i>n} \bigcup_{k\geq 1}\Hat{Z}_{i,k} \big) ,\quad \text{and}\quad \mathcal{S}_{n}:=\bigcup_{i=2}^{n}\bigcup_{k>n-i}\Hat{Z}_{i,k}.
\]
We then clearly have 
\begin{equation}\label{eq:tail0}
    \{\tau >n \} = \mathcal{R}_{n}\cup\mathcal{S}_{n}.
\end{equation}

\begin{lemma}\label{lem:tail1}
For all $n\geq 1,$ if \(\ell_{1}^{+} \neq 0\) and \(\ell_{2}^{+} \neq 0\), then
\begin{equation}\label{eq:ufc330}
m(\mathcal{R}_{n} ) 
\lesssim n^{\frac{-1}{\ell_{2}^{+}k_{1}^{+}}}+  n^{\frac{-1}{\ell_{1}^{+}k_{2}^{+}}} \lesssim n^{\frac{-1}{\beta}},
\end{equation}
and if $\ell_{1}^{+} =0$ and/or $\ell_{2}^{+}=0$ then the term $n^{\frac{-1}{\ell_{1}^{+}k_{2}^{+}}}$ and/or $n^{\frac{-1}{\ell_{2}^{+}k_{1}^{+}}}$ is replaced by $(1+b^{-}_{1})^{\frac{-n}{k_{2}^{+}}}$ and/or $(1+b^{-}_{2})^{\frac{-n}{k_{1}^{+}}}$ respectively. 
\end{lemma}
\begin{proof}
Notice that,
\begin{equation*}
    \bigcup_{i>n} \bigcup_{k\geq 1}\Hat{Z}_{i,k} = \bigcup_{\omega \in \mathbb{S}^{1}} \{\omega\} \times [Z_{n,0}(\omega),0], \quad \text{therefore}\quad m\left(\bigcup_{i>n} \bigcup_{k\geq 1}\Hat{Z}_{i,k}\right) \leq \max_{\omega}|Z_{n,0}(\omega)|.
\end{equation*}
If $\ell_{2}^{+} \neq 0$, then by \eqref{eq:Z_n'} 
\begin{equation}\label{eq:estem1}
    m(\bigcup_{i>n} \bigcup_{k\geq 1}\Hat{Z}_{i,k}) \lesssim n^{\frac{-1}{\ell_{2}^{+}k_{1}^{+}}}.
\end{equation}
Similarly,
\begin{equation}\label{eq:estem2}
    \bigcup_{k\geq n} \Hat{Z}_{1,k}= \bigcup_{\omega \in \mathbb{S}^{1}}\{\omega\} \times [Z_{0,0}(\omega),Z_{1,n-1}(\omega)].
\end{equation}
By \eqref{eq:Z_Y}, we have $T_{\omega}([Z_{0,0}(\omega),Z_{1,n-1}(\omega)])= [0,Y_{0,n-1}(F(\omega))]$, and by Remark \ref{rem:derivative_ge1}, \(\frac{\partial}{\partial x}T_{\omega}\geq 1\) in $[Z_{0,0}(\omega),Z_{1,n-1}(\omega)]$. Therefore, if \(\ell_{1}^{+}\neq 0\), then by Mean Value theorem and \eqref{eq:S_n}  we have 
\begin{equation*}
    |Z_{0,0}(\omega)-Z_{1,n-1}(\omega)|\leq |Y_{0,n-1}(F(\omega))|\lesssim n^{\frac{-1}{\ell_{1}^{+}k_{2}^{+}}}.
\end{equation*}
Combining this with \eqref{eq:estem2} we have,
\begin{equation}\label{eq:estem4}
    m\left(\bigcup_{k\geq n} \Hat{Z}_{1,k}\right) \lesssim n^{\frac{-1}{\ell_{1}^{+}k_{2}^{+}}}.
\end{equation}
Combining \eqref{eq:estem1} and \eqref{eq:estem4} we have the required result. Clearly, if $\ell_{2}^{+} =0$ and/or \(\ell_{1}^{+}=0\), by Lemma~\ref{lem:estimats_X}, the right side of \eqref{eq:estem1} and/or \eqref{eq:estem4} are replaced by \((1+b_{2}^{-})^{\frac{-n}{k^{+}_{1}}}\) and/or \((1+b_{1}^{-})^{\frac{-n}{k^{+}_{2}}}\) respectively\todelete{ }, and hence the stated conclusion follows.
\end{proof}
It remains to measure \(\mathcal{S}_{n}\). First of all we get the following. 
\begin{lemma}\label{lem:tail1000}
   For all $n \geq 2$, 
 \begin{equation}\label{eq:bound_for_A}      
 m(\mathcal{S}_{n})\lesssim \sum_{i=2}^{n}   \int_{0}^{1}|Z_{i-1,0}(\omega)||Y_{0,n-i}(F^{i}(\omega))|~d\omega
\end{equation}
\end{lemma}
\begin{proof}
First of all recall Remark \ref{rem:Znk}, the strip \(\hat{Z}_{i,k}\) is the region between the graph of \(Z_{i,k}\) and \(Z_{i,k-1}\), therefore
\begin{equation}\label{eq:tail10}
m \big(\Hat{Z}_{i,k}\big) = \int_{0}^{1}  |Z_{i,k}(\omega)-Z_{i,k-1}(\omega)| ~d\omega.
\end{equation}
By Remark \ref{rem:Zn} and \eqref{eq:Z_Y}, we have  \(\{\omega\} \times [Z_{i,k}(\omega),Z_{i,k-1}(\omega)] \subset \hat{Z}_{i}\). Therefore, by Mean Value Theorem and Proposition \ref{pro:vdistortion} we have 
\begin{equation*}
    \frac{|Z_{i,k}(\omega)-Z_{i,k-1}(\omega)|}{|Z_{i-1,0}(\omega)-Z_{i,0}(\omega)|} \lesssim \frac{|T^{i}_{\omega}([Z_{i,k}(\omega),Z_{i,k-1}(\omega)]|}{|T^{i}_{\omega}([Z_{i-1,0}(\omega),Z_{i,0}(\omega)]|} = \frac{|Y_{0,k-1}(F^{i}(\omega))-Y_{0,k}(F^{i}(\omega))|}{|Y_{0,0}(F^{i}(\omega))|}.
\end{equation*}
Notice that the graph of \(Y_{0,0}\) is bounded away from \(0\), therefore using the above in \eqref{eq:tail10} we have
\begin{equation*}
  m \big(\Hat{Z}_{i,k}\big)  \lesssim \int_{0}^{1}|Y_{0,k}(F^{i}(\omega))-Y_{0,k-1}(F^{i}(\omega))| |Z_{i-1,0}(\omega)-Z_{i,0}(\omega)| ~d\omega.
\end{equation*}
Combining this with the definition of \(\mathcal{S}_{n}\), we have
\begin{equation}\label{eq:ff}
 m(\mathcal{S}_{n}) \lesssim \int_{0}^{1} \sum_{i=2}^{n}  |Z_{i-1,0}(\omega)-Z_{i,0}(\omega)|\sum_{k >n-i} |Y_{0,k}(F^{i}(\omega))-Y_{0,k-1}(F^{i}(\omega))| ~d\omega.
\end{equation}
By the pointwise monotonicity of these sequences of functions, recall Remark \ref{rem:Znk} and \ref{rem:Ynk} we have \(|Z_{i-1,0}(\omega)-Z_{i,0}(\omega)|\leq |Z_{i-1,0}(\omega)|\) and $\sum_{k> n-i}  |Y_{0,k}(F^{i}(\omega))-Y_{0,k-1}(F^{i}(\omega))| =|Y_{0,n-i}(F^{i}(\omega))|$, using these in \eqref{eq:ff} we have \eqref{eq:bound_for_A}.
\end{proof}

\begin{lemma}\label{lem:tail4-old}
For all $n\geq 2$, if \(\ell_{1}^{+} \neq 0\) or \(\ell_{2}^{+} \neq 0\) then
    \begin{equation}\label{eq:ufc1-old}
        \mathcal{S}_{n}\lesssim n^{\frac{-1}{\beta}},
    \end{equation}
and if $\ell_{1}^{+}=\ell_{2}^{+}=0$, then
\begin{equation}\label{eq:ufc2-old}
    \mathcal{S}_{n} \lesssim  n(1+b^{-})^{\frac{-n}{k^{+}}}.
\end{equation}
\end{lemma}
\begin{proof}
Let $\ell_{1}^{+} \neq 0$ and $\ell_{2}^{+} \neq 0$. Therefore, putting the estimates from \eqref{eq:S_n} and~\eqref{eq:Z_n'} in \eqref{eq:bound_for_A} we have
\begin{equation}\label{eq:estimats_for_A3-old}
      \mathcal{S}_{n} \lesssim  \int_{0}^{1} \sum_{i=2}^{n-1}|Z_{i-1,0}(\omega)||Y_{0,n-i}(F^{i}(\omega))| +|Z_{n-1}(\omega)||Y_{0,0}(\omega)| \lesssim  \sum_{i=2}^{n-1} (i-1)^{\frac{-1}{\ell_{2}^{+}k_{1}^{+}}}(n-i)^{\frac{-1}{\ell_{1}^{+}k_{2}^{+}}}+(n-1)^{\frac{-1}{\ell_{2}^{+}k_{1}^{+}}}.
\end{equation}
Notice that
\begin{equation*}
\begin{split}
   \sum_{i=2}^{n-1} (i-1)^{\frac{-1}{\ell_{2}^{+}k_{1}^{+}}}(n-i)^{\frac{-1}{\ell_{1}^{+}k_{2}^{+}}} &\leq \sum_{i=1}^{n-2} \left(\frac{1}{i(n-i-1)}\right)^{\frac{1}{\beta}}  \\
    &\leq  2^{\frac{1}{\beta}}\left(\sum_{i=1}^{\lfloor\frac{n-1}{2}\rfloor}\left(\frac{1}{i(n-1)}\right)^{\frac{1}{\beta}}+\sum_{i= \lfloor\frac{n-1}{2}\rfloor+1}^{n-2}\left(\frac{1}{(n-1)(n-i-1)}\right)^{\frac{1}{\beta}}\right)\\
    &\lesssim \frac{2^{\frac{1}{\beta}+1}}{n^{\frac{1}{\beta}}} \sum_{i\geq 1} \left(\frac{1}{i}\right)^{\frac{1}{\beta}}\lesssim \frac{1}{n^{\frac{1}{\beta}}}
\end{split}
\end{equation*}
Using this in \eqref{eq:estimats_for_A3-old} we have \eqref{eq:ufc1-old}.
If $\ell_{1}^{+}\neq 0$ and $\ell_{2}^{+}=0$ (or \(\ell_{1}^{+}=0\) and \(\ell_{2}^{+} \neq 0\)), then by Lemma \ref{lem:estimats_X} we have exponential estimates \(|Z_{i,0}(\omega)|\lesssim (1+b_{2}^{-})^{\frac{-i+1}{k^{+}_{1}}} \lesssim (i-1)^{\frac{-1}{\ell_{1}^{+}k_{2}^{+}}}\), and hence we have the same conclusion as in the previous case.
Finally, if $\ell^{+}_{1}=\ell_{2}^{+}=0$. Recall \(b^{-}\) as the minimum of \(b_{1}^{-}\) and \(b_{2}^{-}\), then by Lemma \ref{lem:estimats_X}
\begin{equation}\label{eq:neg_part-old}
 \mathcal{S}_{n} \lesssim \sum_{i=2}^{n-1}(1+b^{-}_{2})^{\frac{-i+1}{k^{+}_{1}}}(1+b^{-}_{1})^{\frac{-n+i}{k^{+}_{2}}}\leq  \sum_{i=2}^{n-1}(1+b^{-})^{\frac{-n+1}{k^{+}}} \leq (n-3)(1+b^{-})^{\frac{-n+1}{k^{+}}}.
\end{equation}
This completes the proof.
\end{proof}
\begin{proof}[Proof of Proposition \ref{pro:tail_lebesgue}]
Recall that it is enough to estimate the tail with respect to \(m\). If $\ell_{1}^{+} \neq 0$ and $\ell_{2}^{+} \neq 0$. Using the estimates from \eqref{eq:ufc330} and \eqref{eq:ufc1-old}, the equation \eqref{eq:tail0} leads to \(m(\{\tau>n\} \leq m(\mathcal{R}_{n}) +m(\mathcal{S}_{n}) \lesssim n^{-1/\beta}\).
If $\ell_{1}^{+}=0$ and $\ell_{2}^{+} \neq 0$, then \(m(\mathcal{R}_{n}) \lesssim  (1+b_{1}^{-})^{\frac{-n}{k_{2}^{+}}} \lesssim n^{-1/\beta}\), therefore we have the same conclusion as the previous case. The case $\ell_{1}^{+} \neq 0$ and $\ell_{2}^{+} =0$ follows similarly. This completes the proof of \eqref{eq:tail_lebesgue1}.
Finally, if $\ell_{1}^{+} = \ell_{2}^{+}=0$, then by the second part of Lemma \ref{lem:tail1} and inequality \eqref{eq:ufc2-old} we have
\begin{equation*}
     m\big(\{\tau>n\}\big) \lesssim  (1+b_{2}^{-})^{\frac{-n}{k_{1}^{+}}}+ (1+b_{1}^{-})^{\frac{-n}{k_{2}^{+}}}+ n(1+b^{-})^{\frac{-n+1}{k^{+}}}\lesssim n(1+b^{-})^{\frac{-n+1}{k^{+}}}  \lesssim e^{-\frac{\ln(1+b^{-})}{2k^{+}} n}.
\end{equation*}
Therefore we have \eqref{eq:tail_lebesgue2}.
\end{proof}

\paragraph{Data availability.}
Data sharing is not applicable to this article as no data sets were generated or analysed during the current study.

\paragraph{Conflict of interest.}
The authors declare that they have no conflict of interest.

\makeatletter
\renewcommand{\@openbib@code}{\setlength{\itemsep}{0pt}\setlength{\parsep}{0pt}\setlength{\parskip}{0pt}}
\makeatother
\AtBeginEnvironment{thebibliography}{\small}
\newcommand{\etalchar}[1]{$^{#1}$}
\providecommand{\bysame}{\leavevmode\hbox to3em{\hrulefill}\thinspace}
\providecommand{\MR}{\relax\ifhmode\unskip\space\fi MR }
\providecommand{\MRhref}[2]{%
  \href{http://www.ams.org/mathscinet-getitem?mr=#1}{#2}
}
\providecommand{\href}[2]{#2}

\end{document}